\documentclass{article}
\usepackage{euscript,amssymb,amsfonts,amsmath,amsthm,mathtools}
\usepackage{multirow}
\usepackage[%
paperheight=11in,
paperwidth=8.25in,
textheight=234mm,%
left=18mm,%
top=2.5cm,%
right=18mm,%
headsep=0.3cm,%
headheight=14pt,%
footskip=0.5cm]{geometry}

\usepackage[utf8x]{inputenc}
\usepackage[T2A]{fontenc}
\usepackage[english]{babel}
\usepackage{boxedminipage}
\usepackage{indentfirst}
\usepackage{afterpage}
\usepackage{paralist}
\usepackage{url}
\usepackage{textcomp}
\usepackage{cite}
\usepackage{mdwlist}
\usepackage{bm}
\usepackage{graphicx}
\usepackage{float}
\usepackage[font=small,labelsep=period,bf]{caption}
\usepackage[usenames,dvipsnames]{color}

\usepackage{multirow}
\usepackage[dvipsnames]{xcolor}

\usepackage{soul}

\newcommand{\Prob}{\mathsf{{P}}}        			
\newcommand{\D}{\mathsf{{D}}}					
\newcommand{\E}{\mathsf{{E}}}					

\newcommand{\R}{\mathbb{R}}               
\newcommand{\N}{\mathbb{N}}               
\renewcommand{\C}{\mathbb{C}}               

\DeclareMathOperator{\sgn}{sgn}

\newtheorem{theorem}{Theorem}
\newtheorem{lemma}{Lemma}
\newtheorem{corollary}{Corollary}

\newtheorem{remark}{Remark}

\newcommand{\eps}{\varepsilon}
\newcommand{\abs}[1]{\left\vert#1\right\vert}
\newcommand{\I}{{\bf1}} 
\newcommand{\lyap}{L_{2+\delta,n}}
\newcommand{\scE}{{\textsc{e}}}
\newcommand{\scR}{{\textsc{r}}}
\newcommand{\esfr}{\widehat L_{\scE,n}}
\newcommand{\es}{L_{\scE,n}}
\newcommand{\rozfr}{L_{\scR,n}}
\newcommand{\roz}{L_{\scR,n}}
\newcommand{\lind}{L_n}
\newcommand{\osiM}{\Lambda_n}
\newcommand{\essM}{M_n}
\newcommand{\essC}{C_\scE}
\newcommand{\rozC}{C_\scR}
\newcommand{\qt}{\sigma^2}
\newcommand{\aexess}{\essC^*}
\newcommand{\aexroz}{\rozC^*}
\newcommand{\gopt}{\gamma_*}

\newcommand{\Sfn}{\overline{f}_n(t)}

\newcommand{\gmmalow}{\Upsilon}
\newcommand{\gmmahigh}{\Gamma}

\date{}

\title{Esseen--Rozovskii type estimates for the rate of convergence\\ in the Lindeberg theorem\footnote{Research supported by the Russian Foundation for Basic Research (projects 15-07-02984-a and 16-31-60110-mol\_a\_dk) and by the grant of the President of Russia No.~MD-2116.2017.1.}}

\author{Ruslan Gabdullin$^1$\!\!, Vladimir Makarenko$^1$\!\!, and Irina Shevtsova$^{1,2}$\!\!}

\begin{document}
\maketitle

\begin{abstract}
We present structural improvements of Esseen's (1969) and Rozovskii's (1974) estimates for the rate of convergence in the Lindeberg theorem and also compute the appearing absolute constants. We introduce the asymptotically exact constants in the constructed inequalities and obtain upper bounds for them. We analyze the values of Esseen's, Rozovskii's, and Lyapunov's fractions, compare them pairwise and provide some extremal distributions. As an auxiliary statement, we prove a sharp inequality for the quadratic tails of an arbitrary distribution (with finite second order moment) and its convolutional symmetrization.
\end{abstract}

\footnotetext[1]{Faculty of Computational Mathematics and Cybernetics,  Lomonosov Moscow State University, Moscow, Russia}

\footnotetext[2]{Institute of Informatics Problems of Federal Research Center ``Computer Science and Control,'' Russian Academy of Sciences, Moscow, Russia; e-mail: \url{ishevtsova@cs.msu.ru}}

\section{Introduction}

Let $X_1,X_2,\ldots,X_n$ be independent random variables (r.v.'s) on
a certain probability space $(\Omega,\mathcal A,\Prob)$ with
distribution functions (d.f.'s) $F_k(x)\coloneqq\Prob(X_k<x),$
$k=1,\ldots,n,$ and such that
\begin{equation}\label{EDX1}
\E X_k = 0,\quad \sigma_k^2\coloneqq\E X_k^2<\infty,\quad B_n^2\coloneqq\sum_{k=1}^n\sigma_k^2>0.
\end{equation}
Denote
$$
S_n=\sum_{k=1}^n X_k,\quad
\overline F_n(x)=\Prob(S_n <xB_n),\quad \Phi(x)=\frac1{\sqrt{2\pi}}\int_{-\infty}^{x} e^{-t^2/2}dt,
$$
$$
\Delta_n=\Delta_n(F_1,\ldots,F_n)\coloneqq\sup_{x\in\R}\abs{\overline F_n(x)-\Phi(x)},\ n\in\N,
$$
$$
\qt_k(z)=\E X_k^2\I(|X_k|\ge z),\quad \lind(z)\coloneqq\frac1{B_n^2}\sum_{k=1}^n\qt_k(zB_n)=\frac1{B_n^2}\sum_{k=1}^n\E X_k^2\I(|X_k|\ge zB_n),\quad z\ge0,
$$
so that $\sup_{z\ge0}\qt_k(z)=\qt_k(0)=\sigma_k^2$ for all $k=1,\ldots,n$ and $L_n(0)=1.$

In 1922 Lindeberg~\cite{Lindeberg1922} proved that $\Delta_n\to0,$
if the \textit{Lindeberg fraction} $\lind$ satisfies the condition
$$
\lim_{n\to\infty}\lind(z)=0\quad \text{for every } z>0.\eqno(L)
$$
In 1935 Feller~\cite{Feller1935} completed Lindeberg's theorem by
proving the necessity of condition $(L)$ for the CLT to hold, if the
random summands are asymptotically negligible in the sense that
$$
\lim\limits_{n\to\infty}B_n^{-2}\max\limits_{1\le k\le n}\sigma_k^2=0.\eqno(F)
$$
Condition $(F)$ is called the \textit{Feller condition}.

Lindeberg's theorem yields the celebrated Lyapunov
theorem~\cite{Lyapounov1901}, according to which $\Delta_n\to0$, if
the \textit{Lyapunov fraction} defined below tends to zero:
$$
\lyap\coloneqq\frac1{B_n^{2+\delta}} \sum_{k=1}^n\E|X_k|^{2+\delta}\to0, \quad n\to\infty, \text{ for certain } \delta>0.
$$
In the case of independent and identically distributed (i.i.d.)
random summands the Lyapunov fraction
$$
\lyap=\frac{\E|X_1|^{2+\delta}}{\sigma_1^{2+\delta}n^{\delta/2}}
$$
is of order $O(n^{-\delta/2})$ as $n\to\infty$ for every fixed
distribution $F_1$ of the random summands.

The numerical demonstration to the Lyapunov theorem is given by the
inequality
\begin{equation}\label{B-E(delta)}
\Delta_n\le C_{BE}(\delta)\cdot \lyap,\quad 0<\delta\le1,
\end{equation}
$C_{BE}(\delta)$ being  some absolute constants for every
$\delta\in(0,1]$, which was proved  by Lyapunov himself in the same
paper~\cite{Lyapounov1901} for $0<\delta<1$, and, 40 years later,
independently by Berry~\cite{Berry1941} in the i.i.d. case  and
Esseen~\cite{Esseen1942} in the general situation for $\delta=1$.
Inequality~\eqref{B-E(delta)} in the i.i.d. case with $\delta=1$
takes the form
\begin{equation}\label{B-E}
\Delta_n\le \frac{C_{BE}(1)}{\sqrt{n}}\cdot\frac{\E|X_1|^3}{\sigma_1^3}
\end{equation}
and establishes the exact order $O(n^{-1/2})$ of the rate of
convergence in the CLT for distributions with finite third-order
moments. The history of evaluation of the constants
$C_{BE}(\delta)$, especially, for $\delta=1$, as well as the
discussion of the structural improvements of the Berry--Esseen
inequality  is worthy of a particular narration for which we refer
to the papers~\cite{KorolevShevtsova2009, KorolevShevtsova2010SAJ,
Shevtsova2013Inf, Shevtsova2014DAN}. In these papers some structural
improvements of inequality~\eqref{B-E(delta)} are also presented.

The numerical demonstration to the Lindeberg theorem was given in
1966 by Osipov's inequality~\cite{Osipov1966} in the form
\begin{equation}\label{OsipovIneq}
\Delta_n\le C(\osiM(\eps)+\lind(\eps))\quad\text{for every }\eps>0,
\end{equation}
where
$$
\osiM(\eps)\coloneqq  \frac1{B_n^3}\sum_{k=1}^n\E|X_k|^3\I(|X_k|<\eps B_n),\quad \eps>0,
$$
and $C$ is an absolute constant. The quantity $\osiM(\eps)+\lind(\eps)$ is called the \textit{Osipov fraction} with the parameter $\eps>0$. Since $\osiM(\eps)\le \eps B_n^{-2} \sum_{k=1}^n \E
X_k^2\I(|X_k|\le\eps B_n)\le\eps$, inequality~\eqref{OsipovIneq} yields the bound
$$
\Delta_n\le C(\eps+\lind(\eps)) \quad\text{for every }\eps>0,
$$
and hence, the Lindeberg theorem. Thus, in the case of
asymptotically negligible random summands satisfying condition
$(F)$, by the Feller theorem,  the left-hand and right-hand sides
of~\eqref{OsipovIneq} either tend or do not tend to zero
simultaneously. In other words, inequality~\eqref{OsipovIneq}
relates the \textit{criteria} of convergence with the \textit{rate}
of convergence in the classical CLT and hence, according to
Zolotarev's terminology~\cite{Zolotarev1986}, it can be called a
\textit{natural} convergence rate estimate.

It is easy to see that inequality~\eqref{OsipovIneq} with
\textit{some} $\eps>0$ yields the Berry--Esseen
inequality~\eqref{B-E(delta)} for \textit{every} $\delta\in(0,1]$,
since
$$
\osiM(\eps)+\lind(\eps)=\sum_{k=1}^n\bigg(\frac{\E|X_k|^3\I(|X_k|<\eps B_n)}{B_n^3} +\frac{\E X_k^2\I(|X_k|\ge\eps B_n)}{B_n^2} \bigg)\le \sum_{k=1}^n\frac{\E|X_k|^{2+\delta}}{B_n^{2+\delta}}=\lyap.
$$
However, inequality~\eqref{OsipovIneq} with $\eps=1$
\begin{multline}\label{OsipovIneqEps=1}
\Delta_n\le C(\osiM(1)+\lind(1))
=\frac{C}{B_n^3}\sum_{k=1}^n\left(\E|X_k|^3\I(|X_k|< B_n)+B_n\E
X_k^2\I(|X_k|\ge B_n)\right)=
\\
=\frac{C}{B_n^2}\sum_{k=1}^n\E X_k^2\min\left\{1,\frac{|X_k|}{B_n}\right\}=C\int_0^{1}\lind(z)dz,
\end{multline}
is of great interest, because $\eps=1$ minimizes the right-hand side
of~\eqref{OsipovIneq}. Indeed, to find the minimizer one can
observe, following the outline of the reasoning
in~\cite{Paditz1986}, that for every $\eps>0$
$$
\osiM(\eps )+\lind(\eps )= \sum_{k=1}^n\bigg[\E\Big|\frac{X_k}{B_n}\Big|^3\I\Big(\Big|\frac{X_k}{B_n}\Big|<\eps\Big) +\E\Big(\frac{X_k}{B_n}\Big)^2 \I\Big(\Big|\frac{X_k}{B_n}\Big|\ge\eps\Big)\bigg] 
$$
and for every r.v. $X\coloneqq X_k/B_n,$ $k=1\ldots,n,$ with  $A$
being  arbitrary Borel subset of $\R$, we have
$$
\E|X|^3\I(|X|<1)+ \E X^2\I(|X|\ge1)= \E X^2\min\{ 1,|X|\}(\I(X\in A) +\I(X\notin A))\le
$$
$$
\le \E|X|^3\I(X \in A) + \E X^2\I(X\notin A)
$$
(here the optimality of the set $A=[-1,1]$ was first noted,
probably, by Loh~\cite{Loh1975} whose hardly available thesis we
cite by~\cite{BarbourHall1984}; however, in~\cite{Paditz1986}, this
fact was proved independently). Now the claim follows by taking
$A\coloneqq[-\eps,\eps]$.

Inequality~\eqref{OsipovIneqEps=1} has an interesting history. First
of all we note that it is a trivial corollary to earlier and
formally more general results of Katz~\cite{Katz1963} and
Petrov~\cite{Petrov1965}, but this connection remained unnoticed for
a long time, even by Petrov himself (for example, see~\cite[Ch.\,5,
\S\,3, Theorems~5 and~8]{Petrov1972}). Furthermore,
inequality~\eqref{OsipovIneqEps=1} was re-proved formally in a more
general, but in fact, in an equivalent form, by
Feller~\cite{Feller1968} with $C=6$ and by
Paditz~\cite{Paditz1980,Paditz1984} in the i.i.d. case with
$C=4.77$. Two years later Paditz~\cite{Paditz1986} announced a
sharper bound $C\le3.51$ without a complete proof. The works of
Feller and Paditz did not refer to the above mentioned work of
Osipov. Inequality~\eqref{OsipovIneq} was also re-proved by Barbour
and Hall~\cite{BarbourHall1984}  who cited Feller's result, but
applied the new Tikhomirov--Stein method and obtained a rougher
bound $C\le18$. Chen and Shao~\cite{ChenShao2001} who cited only
Feller's work re-proved~\eqref{OsipovIneq} with $C=4.1$. Finally,
Korolev and his disciples successively improved the upper bounds for
$C$ to $2.011$ in~\cite{KorolevPopov2011,KorolevPopov2012} and to
$1.87$ in~\cite{KorolevDorofeyeva2017}.

Thus, inequalities~\eqref{OsipovIneq} and~\eqref{OsipovIneqEps=1}
attracted much attention. However, as far back as in 1969
Esseen~\cite{Esseen1969} proved the bound
\begin{equation}\label{EsseenIneqSupZ>0}
\Delta_n\le \frac{C_1}{B_n^3}\sum_{k=1}^n\sup_{z>0}\left\{\abs{\E X_k^3\I(|X_k|<z)} + z\E X_k^2\I(|X_k|\ge z)\right\} = C_1\cdot\esfr^3(\infty),
\end{equation}
where
$$
\esfr^3(\eps)\coloneqq\frac{1}{B_n^3}\sum_{k=1}^n\sup_{0<z<\eps B_n}\left\{\abs{\mu_k(z)} + z\qt_k(z)\right\}, \quad \mu_k(z)\coloneqq\E X_k^3\I(|X_k|<z),
\quad \eps>0,\ z\ge0,
$$
and, as its corollary, by the use of the traditional truncation
techniques, deduced that
\begin{equation}\label{EsseenIneqSupZbounded}
\Delta_n\le C_2\cdot\esfr^3(1),
\end{equation}
with $C_1$ and $C_2$ being some absolute constants. Obviously,
$\esfr^3(1)\le\esfr^3(\infty)$. Moreover, due to the left-continuity
of the functions $\mu_k(z)$, $\qt_k(z)$, $k=1,\ldots,n,$ for $z>0$,
the least upper bound in the definition of $\esfr^3(\eps)$ can be
replaced by the one over the set $z\in(0,\eps B_n].$

Comparing inequalities~\eqref{EsseenIneqSupZ>0}
and~\eqref{EsseenIneqSupZbounded}  with~\eqref{OsipovIneqEps=1} we
observe, first, that the fractions $\esfr^3(\infty)$ and
$\esfr^3(1)$ here depend not on the \textit{absolute}, but on the
\textit{algebraic} third-order truncated moments, which may vanish,
for example, in the case of symmetrically distributed random
summands. Second, Esseen's fraction $\esfr^3(1)$ never exceeds the
Osipov fraction $\osiM(\eps)+\lind(\eps)$ with $\eps=1$ and hence
with arbitrary $\eps>0.$ Indeed, for every random summand
$X\coloneqq X_k$, $k=1,\ldots,n$,  we have
$$
\abs{\E X^3\I(|X|< z)} + z\E X^2\I(|X|\ge z)\le \E|X|^3\I(|X|< z)+ z\E X^2\I(|X|\ge z)\eqqcolon h(z),
$$
with $h(z)$ being monotonically increasing with respect to $z\ge0$,
due to
$$
h(u)-h(v)=\E\left(|X|^3-vX^2\right)\I(v\le|X|<u)+(u-v)\E X^2\I(|X|\ge u)\ge0,\quad u\ge v\ge0,
$$
and, hence, writing $h_k(z)$ for the analogous function of the
distribution of $X_k$, $k=1,\ldots,n,$ we get
$$
\esfr^3(1) \le\frac{1}{B_n^3}\sum_{k=1}^n \sup_{0<z\le B_n}h_k(z) =\frac{1}{B_n^3}\sum_{k=1}^nh_k(B_n)= \osiM(1)+\lind(1) =\inf_{\eps>0}\big\{\osiM(\eps)+\lind(\eps)\big\}.
$$
This fact implies, in turn, that $C\le C_2,$ in particular,
inequality~\eqref{EsseenIneqSupZbounded} is also a natural
convergence rate estimate in the Lindeberg theorem.

Third, Esseen's inequality~\eqref{EsseenIneqSupZbounded} not only
relates the criteria of convergence with the rate of convergence as
Osipov's inequality does, but also provides a numerical
demonstration to the \textit{criteria of the rate of convergence to
be} of the traditional order  $\mathcal O(n^{-1/2})$, $n\to\infty$,
provided by the Berry--Esseen inequality for distributions with
finite third-order moments in the i.i.d. case. In fact, the
condition on the finiteness of the third-order moments can be
relaxed. Namely, in 1966 Ibragimov~\cite{Ibragimov1966} proved that
in the i.i.d. case $\Delta_n=\mathcal O(n^{-1/2})$ as $n\to\infty$
if and only if
$$
\mu_1(z)=\mathcal O(1),\quad z\qt_1(z)=\mathcal O(1),\quad z\to\infty.
$$
It is easy to see that Ibragimov's condition is weaker than the
condition $\E|X_1|^3<\infty$, and Esseen's
inequality~\eqref{EsseenIneqSupZbounded} in the i.i.d. case
$$
\Delta_n\le \frac{C_2}{\sigma_1^3\sqrt{n}}\sup_{0<z\le \sigma_1\sqrt{n}} \left\{\abs{\mu_1(z)} + z\qt_1(z)\right\}
$$
trivially yields the ``if'' part of Ibragimov's criteria.

In 1974 Rozovskii~\cite{Rozovskii1974} generalized Ibragimov's
theorem and proved another estimate involving algebraic truncated
third-order moments
$$
\essM(z)\coloneqq\frac1{B_n^3}\sum_{k=1}^n\mu_k(zB_n)=\frac1{B_n^3}\sum_{k=1}^n\E X_k^3\I(|X_k|<zB_n),\quad z>0,
$$
in the form
\begin{equation}\label{RozovskyIneq}
\Delta_n\le C_3\cdot\rozfr^3,
\end{equation}
where $C_3$ is an absolute constant and
$$
\rozfr^3\coloneqq \abs{M_n(1)} + \sup_{0<z\le1}z\lind(z)=\frac{1}{B_n^3}\bigg(\bigg|\sum_{k=1}^n\mu_k(B_n)\bigg| + \sup_{0<z\le B_n}z\sum_{k=1}^n\qt_k(z)\bigg).
$$
At first glance, Rozovskii's inequality~\eqref{RozovskyIneq} is more
favorable than Esseen's inequalities~\eqref{EsseenIneqSupZ>0}
and~\eqref{EsseenIneqSupZbounded}. Indeed, the right-hand side
of~\eqref{RozovskyIneq} is always finite, while the right-hand side
of~\eqref{EsseenIneqSupZ>0} may be infinite. In~\eqref{RozovskyIneq}
the first term may vanish not only in the symmetric, but also in the
non-symmetric case, for example, for even $n$ if the appearing
truncated third-order moments $\mu_k(B_n)$, $k=1,\ldots,n,$ have the
same absolute value, but alternating signs. One more advantage of
inequality~\eqref{RozovskyIneq} over~\eqref{EsseenIneqSupZ>0}
and~\eqref{EsseenIneqSupZbounded}  is that the values of $\mu_k(z)$
are used in~\eqref{RozovskyIneq} only in one and the same point
$z=B_n$ for all $k=1,\ldots,n,$ while the right-hand sides
of~\eqref{EsseenIneqSupZ>0} and~\eqref{EsseenIneqSupZbounded}
require the information on $\mu_k(z)$ in every point of the interval
$z\in(0,B_n].$ However, a deeper analysis (see
Theorem~\ref{ThEssRozOsipovLyapFracComparison}\,(ii) below) shows
that Rozovskii's fraction $\rozfr^3$ may take greater values than
each of Esseen's fractions $\esfr^3(1)$, $\esfr^3(\infty)$ and even
greater than Lyapunov's fraction $L_{3,n}$, while Esseen's fractions
always satisfy
\begin{equation}\label{Fractions:Esseen<=Osipov<=Lyap}
\esfr^3(\infty)\ \le\ L_{3,n},\qquad \esfr^3(1)\ \le\ \osiM(1)+\lind(1)\ \le\ \lyap,\quad \delta\in(0,1].
\end{equation}
Thus the choice between inequalities~\eqref{EsseenIneqSupZbounded}
and~\eqref{RozovskyIneq} depends not only on the concrete values of
the fractions $\esfr^3$ and $\rozfr^3$, but also on the values of
the appearing absolute constants $C_2$ and~$C_3$. However, the
values of these constants, as well as of $C_1,$ remain unknown.

In the present paper, we compute upper bounds for  the constants
$C_1,C_2,$ and $C_3$. Moreover, we prove new natural convergence
rate estimates in the Lindeberg--Feller theorem generalizing
Esseen's~\eqref{EsseenIneqSupZ>0}, \eqref{EsseenIneqSupZbounded} and
Rozovskii's~\eqref{RozovskyIneq} inequalities, and provide an
explicit analytical representation with an algorithm of evaluation
of the appearing constants. Namely, we introduce a
\textit{truncation} parameter $\eps>0$ and a \textit{balancing}
parameter $\gamma>0$ and denote
$$
\es^3(\eps,\gamma) \coloneqq \sup\limits_{0 < z\le \eps} \left\{ \gamma \abs{\essM(z)} + z\lind(z)\right\}= \frac1{B_n^3}\sup_{0 < z\le\eps B_n} \bigg\{ \gamma\bigg|\sum_{k=1}^n\mu_k(z)\bigg|+ z\sum_{k=1}^n\qt_k(z)\bigg\},
$$
$$
\roz^3(\eps,\gamma)\coloneqq \gamma \abs{\essM(\eps)} + \sup_{0<z\le\eps}z\lind(z) = \frac1{B_n^3}\bigg\{\gamma\bigg|\sum_{k=1}^n\mu_k(\eps B_n)\bigg|+ \sup_{0 < z\le\eps B_n} z\sum_{k=1}^n\qt_k(z)\bigg\},
$$
where the least upper bounds with respect to $0<z\le\ldots$ can be
replaced by those over the open sets $0<z<\ldots,$ due to the
left-continuity of $\mu_k(\cdot)$ and $\qt_k(\cdot)$. It is easy to
see that
$$
\es^3(\eps,1)=\sup_{0<z\le\eps}\{\abs{\essM(z)}+z\lind(z)\} \le\sup_{0<z\le\eps}\{\osiM(z)+z\lind(z)\} =\osiM(\eps)+\eps\lind(\eps),
$$
\begin{multline}\label{es^3(eps,1)<=lyap}
\sup_{\gamma\in(0,1]}\es^3(\eps,\gamma)\le \frac1{B_n^3}\sum_{k=1}^n\sup_{0 < z\le\eps B_n} \Big[ \E|X_k|^3\I(|X_k|<z)+ z\E X_k^2\I(|X_k|\ge z)\Big]
\\
\le \frac1{B_n^3}\sum_{k=1}^n\sup_{0 < z\le\eps B_n}z^{1-\delta} \E|X_k|^{2+\delta}= \eps^{1-\delta}\lyap\quad\text{for every}\quad \eps>0,\ \delta\in(0,1],
\end{multline}
moreover, $\roz^3(1,1)=\rozfr^3,$ $\es^3(1,1)\le\esfr^3(1)$ with equality sign,
for example, in the i.i.d. case, and $\es^3(\eps,\gamma)=\roz^3(\eps,\gamma)$ in
the symmetric case.

\begin{theorem}\label{ThEsseenRozovskiiIneq}
Under the above assumptions, for every $\gamma>0$
\begin{eqnarray}\label{MainEsseenIneq}
\Delta_n \le \essC(\eps,\gamma) \cdot \es^3(\eps,\gamma),\quad \eps\in(0,\infty],
\\
\label{MainRozovskyIneq}
\Delta_n \le \rozC(\eps,\gamma) \cdot\rozfr^3(\eps,\gamma),\quad \eps\in(0,\infty),
\end{eqnarray}
where $\essC(\eps,\gamma)$ and $\rozC(\eps,\gamma)$ depend only on
the arguments in the brackets, take finite values for every
$\gamma>0$ and $\eps$ specified above and can be computed for every
$\eps,\gamma$ under consideration by an algorithm provided below in
the proof. Both $\essC(\eps,\gamma)$ and $\rozC(\eps,\gamma)$ are
monotonically decreasing with respect to $\gamma>0,$
$\essC(\eps,\gamma)$ is also monotonically decreasing with respect
to $\eps>0.$ In particular,
$$
\essC(\infty,1)\le \max\big\{\essC(\infty,0.97),\ \essC(4.35,1)\big\}\le 2.66, \quad \essC(1,1)\le\essC(1,0.72) \le 2.73,
$$
$$
\essC(\infty,\infty) \le \max\big\{\essC(\infty,1.43),\ \essC(4,1.62),\ \essC(2.74,3),\ \essC(2.56,\infty)\big\}\le2.65,
$$
$$
\rozC(1,1)\le\sup_{\gamma\ge\gopt}\rozC(1,\gamma) \le 2.73,\quad \inf_{\eps>0,\,\gamma>0}\rozC(\eps,\gamma)\le \inf_{\eps>0}\sup_{\gamma\ge\gopt}\rozC(\eps,\gamma)\le \rozC(2.12,\gopt)\le2.66,
$$
$$
\sup_{\gamma\ge0.4} \max\big\{\rozC(2.63,\gamma),\ \rozC(1.76,\gamma)\big\} \le2.70,
\quad
\sup_{\gamma\ge0.2} \max\big\{\rozC(5.39,\gamma),\ \rozC(1.21,\gamma)\big\} \le2.87,
$$
where
$$
\gopt\coloneqq\frac{1}{\sqrt{6\varkappa}} =0.5599\ldots,\quad \varkappa := x^{-2}\sqrt{\vphantom{\tfrac12}(\cos x-1+x^2/2)^2+(\sin x-x)^2}\,\Big\vert_{x=x_0} = 0.531551\ldots,
$$
and $x_0=5.487414\ldots$ is the unique root of the equation $8(\cos
x - 1) + 8x\sin x  - 4x^2\cos x-x^3\sin x = 0$ on the interval
$(\pi,2\pi).$
\end{theorem}

\begin{remark}
Within the method used: (i) $\rozC(1,\gamma)$ does not depend on
$\gamma$ for $\gamma\ge\gopt$; (ii) further increase of
$\gamma\ge0.73$ does not reduce the constructed upper bound $2.73$
for $\essC(1,\gamma)$ by more than $0.01$; (iii) the same concerns
the presented upper bounds for $\essC(\infty,1)$ and
$\essC(\infty,\infty)$.
\end{remark}

\begin{remark}
Within the method used, $\essC(1,1) =\rozC(1,1)$, due to
Remark~\ref{Rem:p_E=p_R,|t|<=t_0}.
\end{remark}

\begin{figure}[h!]
\includegraphics[width=0.5\textwidth]{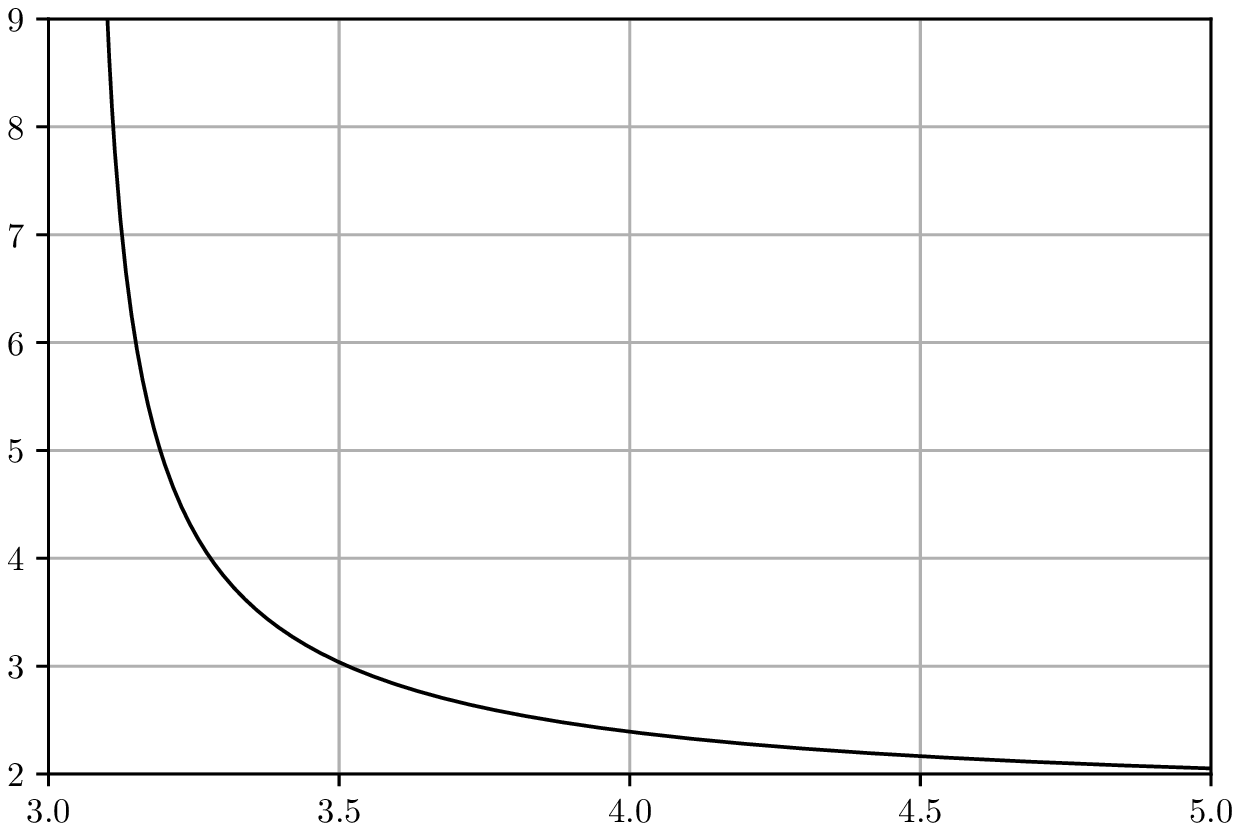}
\includegraphics[width=0.5\textwidth]{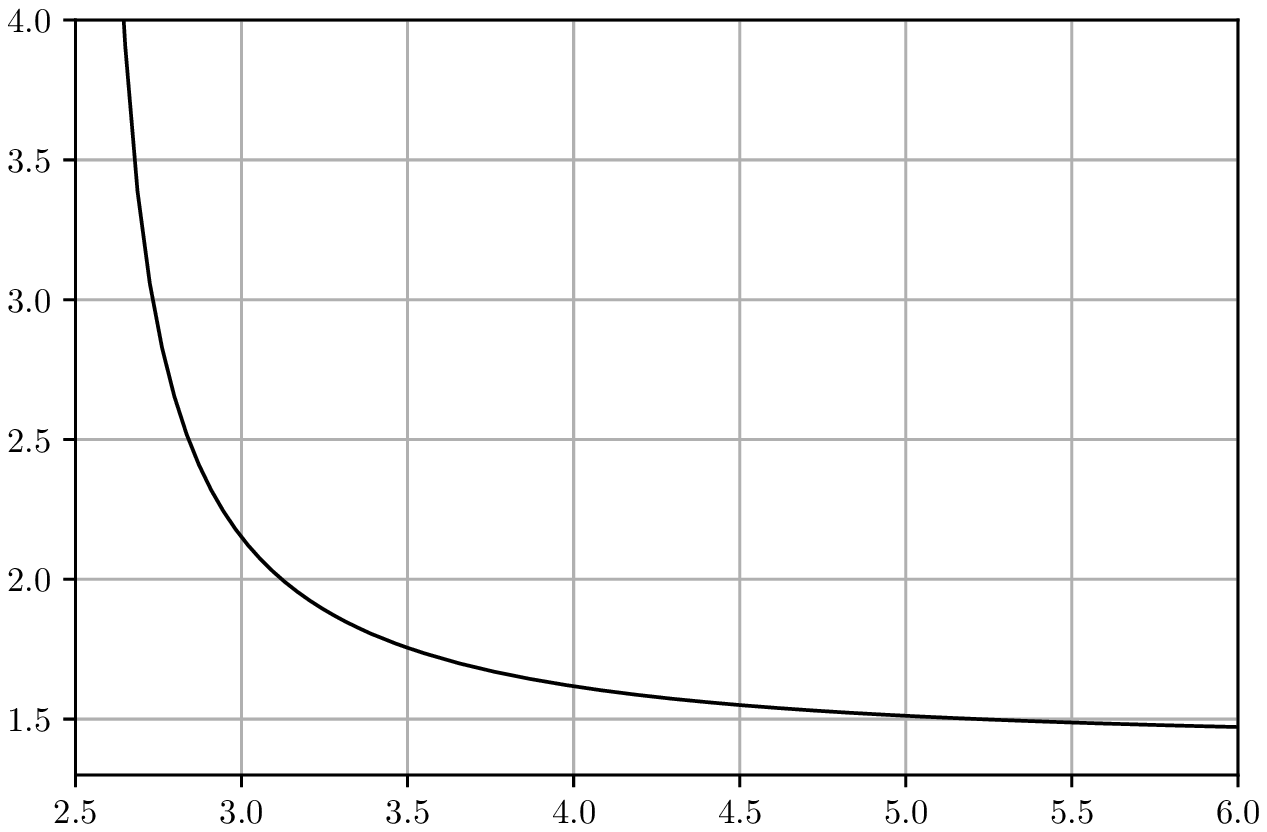}
    \caption{Level curves $\gamma=\gamma(\eps)$ for the upper bounds to the asymptotically exact ``constant'' $\aexess(\eps,\gamma)$ defined in~\eqref{aexCessDef} and to the absolute $\essC(\eps,\gamma)$ ``constant'' in the Esseen-type inequality~\eqref{MainEsseenIneq}. Left: $\big\{(\eps,\gamma)\colon\widehat\aexess(\eps,\gamma)=1.72\big\}$ with $\widehat\aexess(\eps,\gamma)\coloneqq C_0(\eps,\gamma,0+)$ defined in~\eqref{aexess(eps,gamma)<=}. Right: $\big\{(\eps,\gamma)\colon\max\limits_{L_0\le L\le L_1}C_1(\eps,\gamma,L)=2.65\big\}$ with $C_1(\eps,\gamma,L)$ defined in~\eqref{C1(L)def}.}
    \label{Fig:aex+absEssPlotGamma(Eps)}
\end{figure}

The plot of the  level curve $\gamma=\gamma(\eps)$ delivering the
constant value $2.65$ to $\essC(\eps,\gamma)$ is given on
Fig.\,\ref{Fig:aex+absEssPlotGamma(Eps)} (right). The plot of the
function $\rozC(\eps,\gopt)$ constructed in the proof is given on
Fig.\,\ref{Fig:aex+absRozPlotEps} (right).

Theorem~\ref{ThEsseenRozovskiiIneq} trivially yields

\begin{corollary}\label{CorEsseenRozClassicalConst<=}
The constants $C_1,$ $C_2,$ and $C_3$ in
inequalities~\eqref{EsseenIneqSupZ>0},
\eqref{EsseenIneqSupZbounded}, and \eqref{RozovskyIneq} satisfy
$$
C_1\le \essC(\infty,1)\le2.66,\quad C_2\le\essC(1,1)\le2.73,\quad C_3\le\rozC(1,1)\le2.73.
$$
\end{corollary}

The upper bounds for the constants $C_1,C_2,C_3$ presented in
Corollary~\ref{CorEsseenRozClassicalConst<=} are slightly greater
than the best known upper bound $1.87$, obtained
in~\cite{KorolevDorofeyeva2017}, for the absolute constant $C$
(which is no greater than $C_2$) in Osipov's
inequality~\eqref{OsipovIneqEps=1}. In
Theorem~\ref{ThEssRozOsipovLyapFracComparison}\,(iv) of
Section~\ref{SecComparisons} we provide various examples of
symmetric and non-symmetric distributions of the random summands
$X_1,\ldots,X_n$ for which the right-hand sides of both
inequalities~\eqref{MainEsseenIneq} and~\eqref{MainRozovskyIneq}
with $\eps=\gamma=1$, $\essC(1,1)=\rozC(1,1)=2.73$ are strictly less
than the right-hand side of the Osipov
inequality~\eqref{OsipovIneqEps=1} with $C=1.87$.

Similarly to Kolmogorov~\cite{Kolmogorov1953}, where the classical
Berry--Esseen inequality was discussed, we also introduce here the
so-called \textit{asymptotically exact constants}
in~\eqref{MainEsseenIneq}, \eqref{MainRozovskyIneq}
\begin{equation}\label{aexCessDef}
\aexess(\eps,\gamma)\coloneqq \limsup_{\ell \to 0} \sup_{n,F_1,\ldots,F_n} \left\{ \frac{\Delta_n(F_1,\ldots,F_n)}{\ell}\colon \es^3(\eps,\gamma)=\ell \right\},\quad \eps,\ \gamma>0,
\end{equation}
\begin{equation}\label{aexCrozDef}
\aexroz(\eps,\gamma)\coloneqq \limsup_{\ell \to 0} \sup_{n,F_1,\ldots,F_n} \left\{ \frac{\Delta_n(F_1,\ldots,F_n)}{\ell}\colon \roz^3(\eps,\gamma)=\ell \right\},\quad \eps,\ \gamma>0,
\end{equation}
and present their upper bounds for every $\eps>0$ and $\gamma>0.$

\begin{theorem}\label{ThAEXupperBounds}
For every $\eps>0$ and $\gamma>0$ we have
\begin{equation}\label{aexess(eps,gamma)<=}
\aexess(\eps,\gamma)\le \frac{4}{\sqrt{2\pi}}+ \frac{1}{\pi}  \bigg[\frac{\varkappa}{\eps}\gmmalow\Big(1,\frac{t_\gamma^2}{2\eps^2}\Big) +\frac{\eps}{12}\gmmalow\Big(2,\frac{t_\gamma^2}{2\eps^2}\Big) +\frac{\sqrt{2(6\varkappa\gamma^2+1)}} {6\gamma}\, \gmmahigh\Big(\frac{3}{2},\frac{t_\gamma^2}{2\eps^2}\Big) \bigg]
\eqqcolon \widehat\aexess(\eps,\gamma),
\end{equation}
\begin{multline}\label{aexroz(eps,all_gamma)<=}
\aexroz(\eps,\gamma)\le\frac{4}{\sqrt{2\pi}}+ \frac{1}{\pi}  \Big[
\frac{\varkappa}{\eps}
\gmmalow\Big(1,\frac{t_{1,\gamma}^{2}}{2\eps^2}\Big) +\frac{\eps}{12}
\gmmalow\Big(2,\frac{t_{1,\gamma}^{2}}{2\eps^2}\Big)+
\frac{\eps}{6}\gmmahigh\Big(2,\frac{t_{2,\gamma}^{2}}{2\eps^2}\Big)+
\\
+\frac{\sqrt2}{6\gamma}\Big(\frac{\sqrt\pi}{2}
-\gmmalow\Big(\frac32,\frac{t_{1,\gamma}^{2}}{2\eps^2}\Big)
-\gmmahigh\Big(\frac32,\frac{t_{2,\gamma}^{2}}{2\eps^2}\Big)\Big)
\Big]
\eqqcolon \widehat\aexroz(\eps,\gamma)
\end{multline}
where $\gmmahigh(r,x) \coloneqq \int_{x}^\infty t^{r-1} e^{-t} dt,$
$\gmmalow(r,x)  \coloneqq\int_0^{x} t^{r-1} e^{-t} dt =\gmmahigh(r)
- \gmmahigh(r,x),$ $r,x>0,$ are the upper and the lower incomplete
gamma functions,
$$
t_\gamma\coloneqq \tfrac2\gamma\big(\sqrt{(\gamma/\gopt)^2+1}-1\big),\quad
t_{2,\gamma}=2\max\big\{\gamma^{-1},\gopt^{-1}\big\}, \quad
t_{1,\gamma}\coloneqq t_{2,\gamma}\cdot\big(1-\sqrt{(1-(\gamma/\gopt)^2)_+}\,\big),
$$
$\varkappa=0.5315\ldots$ and $\gopt=0.5599\ldots$ are defined in
Theorem~$\ref{ThEsseenRozovskiiIneq}.$ Moreover, the function
$\widehat\aexess(\eps,\gamma)$ is monotonically decreasing with
respect to $\eps>0$ and $\gamma>0,$ the function
$\widehat\aexroz(\eps,\gamma)$ is monotonically decreasing with
respect $\gamma>0$ being constant for $\gamma\ge\gopt$ for every
fixed $\eps>0$.  In particular,
\begin{equation}\label{aexess(inf,gamma)<=}
\aexess(\infty,\gamma)\le 
\frac{1}{\sqrt{2\pi}}\Big(4 +\frac16{\sqrt{{\gopt}^{\!\!-2}+\gamma^{-2}}}\Big)=
\begin{cases}
\frac{1}{\sqrt{2\pi}}\left(4+\frac{\sqrt2}{6}\gopt^{-1}\right)=1.7636\ldots, &\gamma=\gopt,
\\[1mm]
\frac{1}{\sqrt{2\pi}}\left(4+\frac16{\sqrt{\gopt^{-2}+1}}\right)=1.7318\ldots, &\gamma=1,
\\[1mm]
\frac{1}{\sqrt{2\pi}}\big(4 +\frac1{6\gopt}\big)=1.7145\ldots,&\gamma\to\infty,
\end{cases}
\end{equation}
\begin{equation}\label{aexroz(eps,gamma)<=}
\sup_{\gamma\ge\gopt}\aexroz(\eps,\gamma)\le \widehat\aexroz(\eps,\gopt)
=\frac{4}{\sqrt{2\pi}}+ \frac{1}{\pi}  \left[\frac{\varkappa}{\eps}\gmmalow\left(1,\frac{2}{(\eps\gopt)^{2}}\right) +\frac{\eps}{12}\left(1+\gmmahigh\left(2,\frac{2}{(\eps\gopt)^{2}}\right)\right) \right],\quad \eps>0,
\end{equation}
$$
\inf_{\eps>0}\inf_{\gamma>0}\aexroz(\eps,\gamma)\le
\inf_{\eps>0}\sup_{\gamma\ge\gopt}\aexroz(\eps,\gamma)\le  \widehat \aexroz(1.89,\gopt)\le1.75.
$$
\end{theorem}

\begin{figure}[h!]
\center
\includegraphics[width=0.5\textwidth]{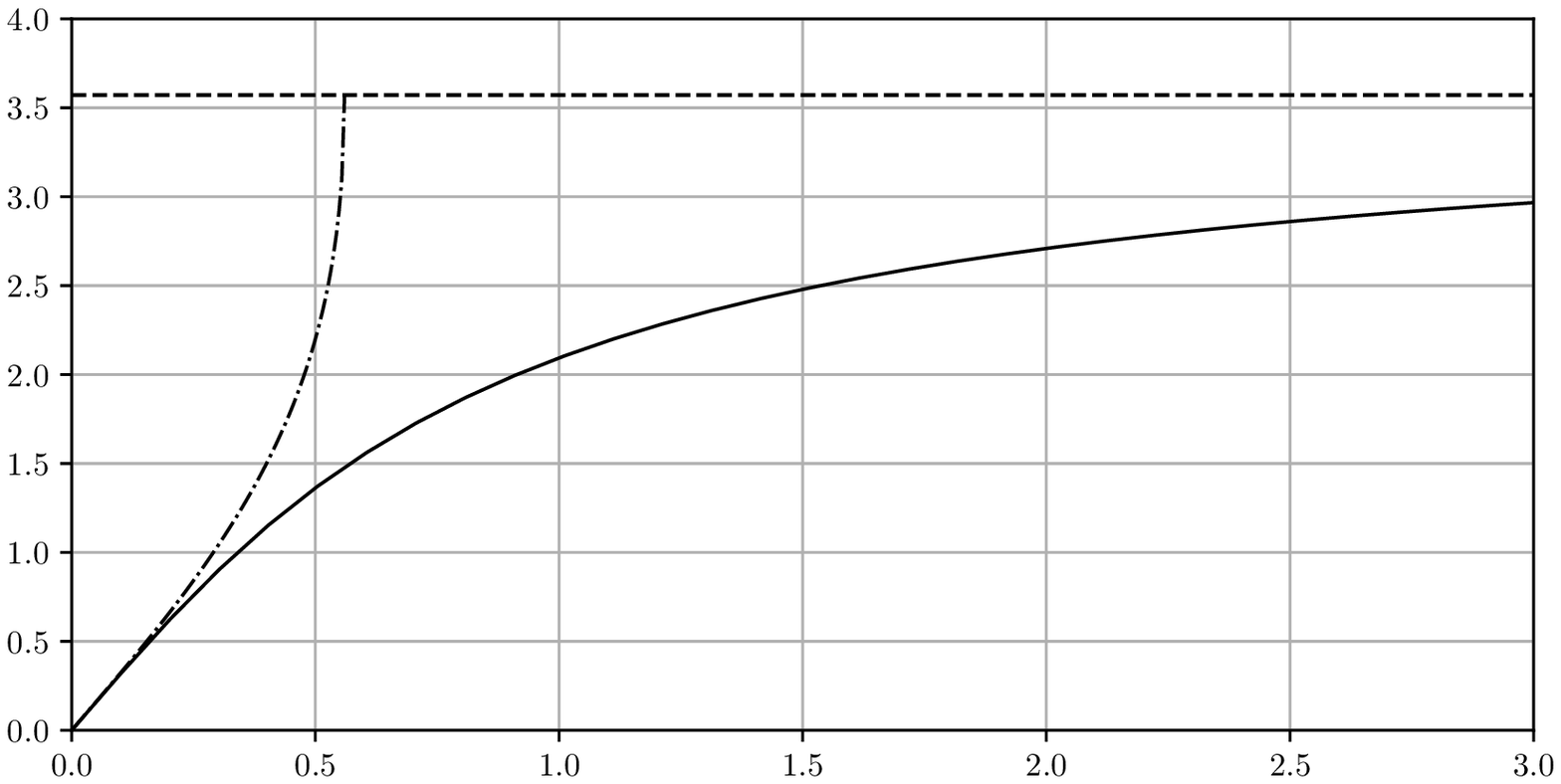}
\caption{Graphs of the functions $t_\gamma\coloneqq \tfrac2\gamma\big(\sqrt{(\gamma/\gopt)^2+1}-1\big)$ (solid line) and $t_{1,\gamma}\coloneqq\tfrac2\gamma\big(1-\sqrt{(1-(\gamma/\gopt)^2)_+}\,\big)$ (dashdot line)  for $\gamma>0$; dashed line represents the limiting value $t_\infty\coloneqq\lim\limits_{\gamma\to\infty}t_\gamma =t_{1,\gopt}=2/\gopt=3.5717\ldots\ .$}
\label{Fig:t_gamma+tilde(t_gamma)-plots}
\end{figure}

The values of the functions $\widehat\aexess(\eps,\gamma)$ and
$\widehat\aexroz(\eps,\gamma)$ for some $\eps>0$ and $\gamma>0$ are
given in the third columns of tables~\ref{Tab:EsseenC0(L)}
and~\ref{Tab:RozovskyC0(L)}, respectively. The plot of the  level curve
$\gamma=\gamma(\eps)$ with the constant value of $\widehat\aexess(\eps,\gamma)=1.72$,
closest to the minimal one $\widehat\aexess(\infty,\infty)$ up to
$10^{-2}$, is given on Fig.\,\ref{Fig:aex+absEssPlotGamma(Eps)}
(left).  The plot of the function $\widehat\aexroz\big(\eps,\gopt)$
is given on Fig.\,\ref{Fig:aex+absRozPlotEps} (left, solid line).
The graphs of the functions $t_\gamma,\ t_{1,\gamma},\ t_{2,\gamma}$
are given on Fig.\,\ref{Fig:t_gamma+tilde(t_gamma)-plots}.

\begin{remark}
Within the method used, the functions $\essC(\eps,\gamma)$ and
$\rozC(\eps,\gamma)$  constructed here in the proof of
Theorem~\ref{ThEsseenRozovskiiIneq} are bounded below by the
functions $\widehat\aexess(\eps,\gamma)$ and
$\widehat\aexroz(\eps,\gamma)$, respectively.

\end{remark}

\begin{remark}
Since
$$
\lim_{\eps\to0}\es(\eps,\gamma)=\lim_{\eps\to0}\roz(\eps,\gamma)=0\quad \text{for every }\ \gamma>0,
$$
the functions $\essC(\eps,\gamma),$ $\rozC(\eps,\gamma),$
$\aexess(\eps,\gamma),$ and $\aexroz(\eps,\gamma)$  must be
unbounded as $\eps\to0.$
\end{remark}

\begin{remark}
Though the first terms $\mu_k(\cdot),$ $k=1,2,\ldots,n,$ in the
definitions of $\es(\eps,\gamma),$ $\roz(\eps,\gamma)$ vanish for
symmetric distributions of random summands,
one cannot, in general, get rid of the third truncated moments
$\mu_k(\cdot)$ in~\eqref{MainEsseenIneq}, \eqref{MainRozovskyIneq},
that is, the ``constants'' $\essC(\eps,0+),$ $\rozC(\eps,0+),$
$\aexess(\eps,0+),$ and $\aexroz(\eps,0+)$ are no more bounded. This
fact follows from the observation that each of the above constants
is bounded below by the so-called \textit{asymptotically best
constant} (we follow here the terminology introduced
in~\cite{Shevtsova2010TVP})
$$
\sup_{F_1=\ldots=F_n\colon B_n>0}\ \limsup_{n\to\infty} \frac{\Delta_n(F_1,\ldots,F_n)}{\sup\limits_{0<z\le\eps}zL_n(z)}
$$
for which, in Theorem~\ref{ThAbsChBE(gamma=0,eps)=inf} of
Section~\ref{SecComparisons}, an infinite lower bound is constructed
for every $\eps>0$.
\end{remark}


The proofs of Theorems~\ref{ThEsseenRozovskiiIneq} and
\ref{ThAEXupperBounds} are based on the method of characteristic
functions (ch.f.'s) realized in the Prawitz smoothing
inequality~\cite{Prawitz1972}  (see Lemma~\ref{LemPrawSmoothIneq}
below) and new estimates for ch.f.'s presented in
Section~\ref{SecCh.F.Estims}. As an auxiliary result used in
bounding the absolute value of a ch.f., we prove a sharp inequality
$$
\E(X-X')^2\I(|X-X'|\ge z)\le 4\E X^2\I(|X|\ge z/2),\quad z\ge0,
$$
for arbitrary centered i.i.d. r.v.'s $X$ and $X'$ with finite second
moment, where the constant factor $4$ on the R.H.S. cannot be made
less (see Theorem~\ref{ThQuadrTailsSymmIneq} in
Section~\ref{SecQuadraticTails}). Similar inequality was proved and
used in the preceding works with the constant factor $40$ by
Esseen~\cite{Esseen1969}, $12$ and $8$ by Rozovskii~\cite{Rozovskii1974,Rozovskii1978c}, respectively.

The proof of the main Theorem~\ref{ThEsseenRozovskiiIneq} is given
in Section~\ref{SecMainThmProof} and consists of two steps. First,
we consider the case of small values of the fractions $\es\le L$ and
$\roz\le L$ with $L$ being small enough including the limit $L\to0$,
where the upper bounds for the corresponding ``constants'' (in fact,
depending on $L$) are obtained in the analytical form (see
Subsection~\ref{SubsecCaseSmallL}) yielding, as a particular case,
Theorem~\ref{ThAEXupperBounds}. Then we consider the remaining case,
where the upper bounds for the  appearing ``constants''  also depend
on the concrete value of the fraction $L\coloneqq\es$ or
$L\coloneqq\roz$ and have a more complicated form assuming numerical
evaluation by the use of a computer (see
Subsection~\ref{SubsecCaseLargeL}).

\renewcommand{\lyap}{L_{3,n}}

In Section~\ref{SecComparisons} we  compare  the fractions
$\es^3(\eps,\gamma)$, $\roz^3(\eps,\gamma)$, $\lyap$ and also the
right-hand sides of Osipov's inequality~\eqref{OsipovIneqEps=1} with
our new inequalities~\eqref{MainEsseenIneq},
\eqref{MainRozovskyIneq}.

All the figures and tables are presented in the concluding Appendix.

\section{An inequality for quadratic tails}
\label{SecQuadraticTails}

\begin{theorem}\label{ThQuadrTailsSymmIneq}
Let $X$ and $X'$ be i.i.d. with the d.f. $F$ and $\E X=0$. Then with
$$
\qt(z)\coloneqq\E X^2\I(|X|\ge z)=\int_{|x|\ge z}x^2dF(x),
$$
$$
\qt_s(z)\coloneqq\E (X-X')^2\I(|X-X'|\ge z)=\int\limits_{|x-y|\ge z}(x-y)^2dF(x)dF(y),\quad z\ge0,
$$
we have
$$
\qt_s(z)\le 2\qt(\alpha z)+2\qt((1-\alpha)z),\quad z\ge0,\quad \alpha\in[0,1].
$$
In particular, with $\alpha=1/2$ we have
\begin{equation}\label{QuadraticTailsSymmIneq}
\qt_s(z)\le 4\qt(z/2),\quad z\ge0,
\end{equation}
where the constant factor $4$ on the right-hand side is the best
possible in the sense of
\begin{equation}\label{QuadraticTailsSymmExtrema}
\sup_{X\colon \E X=0}\frac{\qt_s(z)}{\qt(z/2)}=4\quad \text{for every }\quad z>0,
\end{equation}
where the least upper bound is taken over all distributions of the
r.v. $X$ with $\E X=0$ and delivered by the sequence of two-point
distributions of the form $\Prob(X\hm=qz)\hm=p,$
$\Prob(X=-pz)=q\coloneqq1-p$ with $p\to\frac12+.$
\end{theorem}

\begin{remark}
In the original work of Esseen~\cite{Esseen1969}
inequality~\eqref{QuadraticTailsSymmIneq} was proved with the
constant factor $40$ on the right-hand side instead of $4$, and in
the works of Rozovskii~\cite{Rozovskii1974,Rozovskii1978c} this factor was successively lowered to $12$ and~$8$.
\end{remark}

\begin{proof}
By $F(x,y) = F(x) \cdot F(y)$ denote the d.f. of the random vector
$(X,X')$. Then for every $\alpha\in[0,1]$ with
$\beta\coloneqq1-\alpha$ we have
$$
\qt_s(z)=\int_{|x-y| \ge z} (x-y)^2 dF(x,y)\le
\int_{|x| \ge \alpha z}+ \int_{|y| \ge\beta z} (x-y)^2 dF(x,y)=
$$
$$
=\int_{|x| \ge \alpha z}+ \int_{|y| \ge\beta z} (x^2 - 2xy + y^2) dF(x,y)= \qt(\alpha z)+\qt(\beta z)+ \sigma^2\int_{|x| \ge\alpha z} dF(x)+\sigma^2\int_{|y| \ge\beta z} dF(y),
$$
since $\E X = \E X' = 0.$ Observing that
$$
\sigma^2\int_{|x| \ge\alpha z} dF(x)=\int_{|x| \ge\alpha z} dF(x) \int_{\R} y^2 dF(y)=
\int_{|y| < \alpha z \le |x|}y^2 dF(x) dF(y)+ \qt(\alpha z)\int_{|x| \ge\alpha z}dF(x)\le
$$
$$
\le \int_{|y| <\alpha z} dF(y)\int_{|x| \ge\alpha z} x^2 dF(x) + \qt(\alpha z)\int_{|x| \ge\alpha z}dF(x)=\qt(\alpha z),
$$
we finally obtain
$$
\qt_s(z)\le 2\qt(\alpha z)+2\qt(\beta z),\quad z\ge0.
$$

To prove~\eqref{QuadraticTailsSymmExtrema}, fix arbitrary $z>0$ and
consider a centred two-point distribution
$$
X = \begin{cases}
a, & \frac{b}{a+b},\\
-b, & \frac{a}{a+b},
\end{cases}
$$
where the choice of $b>a>0$ will depend on $z$. Then, with $X'$
being an independent copy of $X$, we have
\[
X-X' = \begin{cases}
a+b, & \frac{ab}{(a+b)^2},
\\[2mm]
0, & \frac{a^2+b^2}{(a+b)^2},
\\
-a-b, & \frac{ab}{(a+b)^2},
\end{cases}
\quad
\qt_s(z) =
\begin{cases}
2ab, & 0 \le z \le a+b,\\
0, &z > a+b,
\end{cases}
\quad
\qt\left(\frac z2\right) =
\begin{cases}
ab, & 0 \le \frac z2 \le a,
\\[2mm]
\frac{ab^2}{a+b}, & a < \frac z2 \le b,
\\
0, & z/2 > b.
\end{cases}
\]
Now let $a$ and $b$ satisfy $a+b=z$. Then $\qt_s(z)=2ab$ and
$\qt(z/2)=ab^2/(a+b)=ab^2/z$, and hence for every $z>0$
\[
\sup_{0<a<b\colon a+b=z} \frac{\qt_s(z)}{\qt\left(\tfrac z2\right)} =\sup_{z/2<b<z}\frac{2z}{b}= \lim_{b\to z/2+}\frac{2z}{b} = 4.
\]
\end{proof}

\section{Estimates for characteristic functions}
\label{SecCh.F.Estims}

In what follows we shall omit the arguments $\eps,\gamma$ of the
fractions $\es(\eps,\gamma)$ and $\roz(\eps,\gamma)$ assuming
$\gamma\in(0,\infty)$ and $\eps\in(0,\infty]$ for $\es(\eps,\gamma)$
or $\eps\in(0,\infty)$ for $\roz(\eps,\gamma)$ being fixed and even
simply using the notation $L$ for each of the fractions $\es,\roz$.
Furthermore, we shall also assume that the random summands
$X_1,\ldots,X_n$ are normalized so that
$$
B_n^2\coloneqq\sum_{k=1}^n\sigma_k^2=1.
$$

For $k=1,\ldots,n$ denote
$$
f_k(t)\coloneqq \E e^{itX_k},\quad \overline f_n(t)=\E \exp\left\{\frac{itS_n}{\sqrt{\D S_n}}\right\} =\prod_{k=1}^nf_k\left(t\right),\quad t\in\R,
$$
the ch.f.'s of $X_k$ and of the (already normalized) sum $S_n.$
Recall that $e^{-t^2/2},$ $t\in\R,$ is the ch.f. of the standard
normal distribution on~$\R$.

\subsection{Estimation of $\big|\overline f_n(t)\big|$}

This section is aimed at bounding the absolute value $\abs{\overline
f_n(t)}$ of the ch.f. of the normalized sum above. We provide two
kinds of such estimates. The first estimate which is the most exact,
but has a rather complicated form, is used below in the numerical
method described in Subsection~\ref{SubsecCaseLargeL} for the case
of $\roz,\es$ separated from the origin. The second estimate has a
simpler form convenient for analytical work and is used below in
Subsection~\ref{SubsecCaseSmallL} for the case of small $\roz,\es$.

We also note here that the form of the estimates for $\abs{\overline
f_n(t)}$ presented below is the same for the both fractions $\es^3$
and $\roz^3$ which is explained by the independence of the
constructed bounds of the function~$\essM(\cdot)$.

\begin{lemma}[see~\cite{Prawitz1991}]\label{LemPrawCos}
For every $\theta \in [0,2\pi]$ and $x \in \R$ we have
\begin{equation}\label{cos(x)<=1-ax^2+bx^4}
\cos x \le 1 - a(\theta) x^2 + b(\theta) x^4,
\end{equation}
with equality if and only if $x\in\{-\theta,0,\theta\},$ where
\[
a(\theta) = 2\frac{1-\cos \theta}{\theta^2} - \frac{\sin\theta}{2\theta},\quad
b(\theta) = \frac{1-\cos \theta}{\theta^4} - \frac{\sin\theta}{2\theta^3}\quad\text{for }\ \theta \in (0,2\pi]
\]
and $a(0)\coloneqq1/2,$ $b(0)\coloneqq1/24$ are defined  by continuity.
\end{lemma}

\begin{lemma}\label{Lem:a(theta)/b(theta)increases}
{\rm (i)} The functions $a(\theta)$ and $b(\theta)$ on the open
interval $\theta\in(0,2\pi)$ are both strictly monotonically
decreasing and strictly positive varying  within the range
$$
\frac12=a(0)>a(\theta)>a(2\pi)=0,\quad \frac1{24}=b(0)>b(\theta)>b(2\pi)=0.
$$

{\rm (ii)} The function $a(\theta)/b(\theta)$ is strictly
monotonically increasing for $\theta \in[0,2\pi)$, in particular,
$$
\sup_{0<\theta<2\pi}\frac{a(\theta)}{b(\theta)}= \lim_{\theta\to2\pi-}\frac{a(\theta)}{b(\theta)}=4\pi^2.
$$
\end{lemma}

\begin{proof}
For $\theta\in(0,2\pi)$ denote $g(\theta)\coloneqq
a(\theta)/b(\theta).$ Then, omitting the argument $\theta$ of the
functions $a(\theta)$ and $b(\theta)$ for short, we have
$$
a' = -4\,\frac{1-\cos \theta}{\theta^3} + 5\,\frac{\sin \theta}{2\theta^2}-\frac{\cos \theta}{\theta},
\quad
b' = -4\,\frac{1-\cos \theta}{\theta^5} + 5\,\frac{\sin \theta}{2\theta^4}-\frac{\cos \theta}{\theta^3}= \frac{a'(\theta)}{\theta^2},
$$
\begin{equation}\label{PropLimit.ineq0}
g'(\theta)=\frac{a'b-b'a}{b^2} =-\frac{b'(a-\theta^2b)}{b^2} =-\frac{b'}{b^2}\cdot \frac{1-\cos \theta}{\theta^2},
\end{equation}
hence, $\sgn g'(\theta)=-\sgn b'(\theta)$ for all $\theta\in(0,2\pi)$.

Observe that
\begin{equation}\label{PropLimit.ineq1}
b_1(\theta)\coloneqq2 \theta^5 b'(\theta) = -8(1-\cos \theta) + 5\theta\sin\theta - \theta^2\cos\theta,\quad b_1(0+) = 0,\quad b_1(2\pi-) = -4\pi^2,
\end{equation}
\begin{equation}\label{PropLimit.ineq2}
b_1'(\theta) = (\theta^2-3) \sin \theta + 3 \theta \cos \theta,\quad
b_1'(0+) = 0,\quad  b_1'(2\pi-) = 6\pi,
\end{equation}
\begin{equation}\label{PropLimit.ineq3}
b_2(\theta)\coloneqq \frac{b_1''(\theta)}{\theta^2} = \cos\theta-\frac{\sin\theta}{\theta},\quad
b_2(0+) = 0,\quad b_2(\pi) = -1,\quad b_2(2\pi) =1.
\end{equation}
It is easy to see that there exists a unique $\theta_0 \in
(\pi,2\pi)$ such that
\begin{eqnarray*}
b_2(\theta) < 0,&\quad 0 < \theta < \theta_0,
\\
b_2(\theta) > 0,&\quad \theta_0 < \theta < 2\pi,
\end{eqnarray*}
With the account of $\sgn b_2=\sgn b_1''(\theta)$ and $b_1'(0+)=0$,
$ b_1'(2\pi-) >0,$ we conclude that $b_1'$ has a unique sign change
in a point  $\theta_1\in(\theta_0,2\pi)\subset(\pi,2\pi)$ and
\begin{eqnarray*}
b_1'(\theta) < 0,&\quad 0 < \theta < \theta_1,\\
b_1'(\theta) > 0,&\quad \theta_1<\theta<2\pi.
\end{eqnarray*}
Hence, $b_1$ has the unique point of minimum $\theta_1$ on
$(0,2\pi)$ and for all $\theta\in(0,2\pi)$
$$
b_1(\theta)<\max\{b_1(0+),b_1(2\pi-)\}=0,
$$
i. e., $b'(\theta)=\frac{b_1(\theta)}{2\theta^5}<0$ for
$\theta\in(0,2\pi)$. Since also $\sgn a'=\sgn b'<0$ on $(0,2\pi)$,
we immediately obtain part (i) of the lemma. Finally, recalling that
$\sgn g'=-\sgn b'>0$ on $(0,2\pi)$ we get part (ii) of the lemma.
\end{proof}

\begin{theorem}\label{ThAbsCh.F.Estim}
For $\eps>0,$ $\tau\ge0,$ $u\ge0,$ and $\theta \in (0,2\pi]$ let
$$
k(\tau,u,\theta)\coloneqq a(\theta) - 4{{\tau}u^{-1}}\left(a(\theta) + b(\theta)u^2\right) \text{ for } u>0,\ \ k(0,0,\cdot)\coloneqq a(\theta),\ \ k(\tau,0,\cdot)\coloneqq -\infty \text{ for } \tau>0,
$$
$$
k(\tau,u) \coloneqq
\sup\left\{ k(\tau,v,\theta) \colon0\le v\le u\wedge \sqrt{\tfrac{a(\theta)}{b(\theta)}},\ \theta\in(0,2\pi)\right\}=
\sup_{0 < \theta<2\pi}  k\left(\tau,u\wedge\sqrt{\tfrac{a(\theta)}{b(\theta)}},\theta\right),
$$
$$
k\left(\tau\right) = \sup\limits_{\theta\in(0,2\pi],\,u>0}k(\tau,u,\theta)=\sup\limits_{0 < \theta<2\pi} k\left(\tau,\sqrt{\tfrac{a(\theta)}{b(\theta)}},\theta\right) =\sup\limits_{0 < \theta\le2\pi}\left\{ a(\theta) - 8\tau\sqrt{a(\theta)b(\theta)}\right\},
$$
$$
$$
where the functions $a(\theta), \, b(\theta)$ are defined in the
formulation of Lemma~\ref{LemPrawCos}. Then for each of the
fractions  $L\in\{\roz,\es\}$ we have

{\rm (i)} for every $ t \in \R$
\begin{equation} \label{AbsCh.F.EstimComp}
\abs{\Sfn} \le \exp \left\{-k\left(L^3|t|,2\eps|t|\right)t^2\right\} ;
\end{equation}

{\rm (ii)} for every  $L_0\ge L,$ $\tau\ge \underline{\tau}_1(L_0,\eps)\coloneqq{\pi L_0^3}/{\eps}$ and $|t| \le \tau/L^3$
\begin{eqnarray}\label{AbsCh.F.EstimAnalytically}
\abs{\Sfn} \le \exp\left\{-k(\tau)t^2\right\},
\end{eqnarray}
with $k(\tau)>0$ if and only if $\tau< \overline{\tau}_1:=
\frac{\pi}{4} = 0.7853\ldots$ and the interval
$[\underline{\tau}_1(L_0,\eps),\overline{\tau}_1)$ being nonempty if
and only if $L_0<(\eps/4)^{1/3}.$

{\rm (iii)} Statements {\rm(i)} and {\rm (ii)} remain true with
$L^3=\sup\limits_{0<z\le\eps}z\lind(z)$.
\end{theorem}

\begin{remark}\label{RemAbsChFEstimInreases}
The right-hand side of~\eqref{AbsCh.F.EstimComp} is monotonically
increasing with respect to $L>0$, monotonically decreasing with
respect to $\eps>0$ and does not depend on $\gamma$ with $L$ being
fixed. The right-hand side of~\eqref{AbsCh.F.EstimAnalytically} does
not depend neither on $L$, nor on $\eps$ and $\gamma$.
\end{remark}

\begin{proof}
Due to the symmetry, one can assume that $t\ge0.$ Moreover,
inequalities~\eqref{AbsCh.F.EstimComp}
and~\eqref{AbsCh.F.EstimAnalytically} hold trivially true with
$t=0$, so that, in what follows, we shall assume that $t>0.$

Consider the absolute value of the ch.f. of every single random
summand. Namely, let $X$ and $X'$ be i.i.d. centred r.v.'s with the
d.f. $F$, the ch.f. $f$, and $\sigma^2\coloneqq\E X^2<\infty$.
Denote $F^s(x)\coloneqq\Prob(X-X'<x),$ $x\in\R,$ the symmetrization
of $F$,
$$
\qt(z)\coloneqq\int_{|x|\ge z}x^2dF(x),\quad\text{and}\quad
\qt_s(z)\coloneqq\int_{|x|\ge z}x^2dF^s(x),\quad z\ge0.
$$
Then for every $t$ under consideration and  $\theta \in (0,2\pi]$ we have
\[
|f(t)|^2 = \E \cos t\left(X - X'\right)= 1 - 2a(\theta)\sigma^2t^2 + J,
\]
with $a(\theta)$ defined in Lemma~\ref{LemPrawCos} and
$$
J=J(t,\theta)\coloneqq\int_\R \left(\cos tx - 1 + a(\theta)t^2x^2\right)dF^s(x)
$$
Separating the range of integration in $J$ into two domains: $|x|\ge
z$ and $|x|<z$ with arbitrary $z \ge 0$, we use the inequality $\cos
x - 1 \le 0,$ $ x \in \R,$ to bound the integrand in the first
integral and Lemma~\ref{LemPrawCos} to bound the second and obtain
\[
J\le at^2\int_{\abs{x}\ge z} x^2dF^s(x) + bt^4\int_{\abs{x}< z} x^4dF^s(x)= \left(a-bt^2z^2\right)t^2\qt_s(z) + 2bt^4\int_0^z x\qt_s(x) dx,
\]
where, for short, we omitted the argument $\theta$ of the functions
$a(\theta)$ and $b(\theta)$ defined in Lemma~\ref{LemPrawCos} and
used the representation
\begin{equation}\label{Destr4Moment}
\int_{\abs{x} < z}x^4 dF(x) =-\int_0^z x^2 d\qt(x) = -z^2\qt(z)+ 2 \int_0^z x\qt(x) dx,\quad z\ge0,
\end{equation}
with $F\coloneqq F^s$, which is valid for arbitrary d.f. $F$ with
finite second moment. Now assuming that $a\ge bt^2z^2$ we may apply
Theorem~\ref{ThQuadrTailsSymmIneq} to bound $\qt_s(\cdot)$ in the
majorant for $J$ above and obtain
\[
J\le 4t^2\left(a - bt^2z^2\right)\qt(z/2) + 8bt^4 \int_{0}^{z} x\qt\left(\frac{x}{2}\right)dx,\quad z\ge0.
\]

Now repeating the same procedure for every random summand
$X\coloneqq X_k,$ $k=1,\ldots,n,$ we construct an estimate for the
absolute value of the corresponding ch.f. $|f_k(t)|$ in the form
$$
\abs{f_k(t)}^2=1 -2a\sigma_k^2 t^2 + J_k,
$$
where
$$
J_k=J_k(t,\theta)\le 4t^2\left(a - bt^2z^2\right)\qt_k(z/2) + 8bt^4 \int_{0}^{z} x\qt_k\left(\frac{x}{2}\right)dx,\quad z\ge0,
$$
with
$$
\sum_{k=1}^{n} J_k \le 4 t^2\left( a - bt^2z^2\right) \lind\left(\frac{z}{2}\right)  + 16bt^4 \int_{0}^{z}\frac{x}2\lind\left(\frac{x}{2}\right)dx \le
$$
$$
\le\left(8t^2z^{-1} \left( a - bt^2z^2\right) +16bt^4z\right)\sup\limits_{0 < x\le z/2} x\lind\left(x\right)
\le8t^2  z^{-1} \left( a + bt^2z^2\right) \sup_{0 < x\le\eps} x\lind\left(x\right)
$$
for all $\eps\ge z/2\ge0$, provided  that $a\ge bt^2z^2$. Thus,
using the inequality $1 + x \le e^x$, $x\in\R$, and observing that
$$
\sup\limits_{0 < x\le\eps} x\lind\left(x\right)\le\min\{\es^3,\roz^3\}\le L^3
$$
(compare with part (iii) of the theorem), we obtain
\begin{multline}\label{tech.LemMod1.1}
\abs{\overline f_n(t)}=\prod\limits_{k=1}^{n} \abs{f_k(t)}\le
\exp\left\{-at^2 + \frac12\,{\textstyle\sum\limits_{k=1}^{n} J_k}\right\}\le \exp\left\{-at^2+4L^3t^2z^{-1}(a+bt^2z^2)\right\}=
\\
=\exp\left\{-t^2\left(a+4L^3z^{-1}(a+bt^2z^2)\right)\right\}\eqqcolon \exp\left\{-t^2\cdot k(L^3t,tz,\theta)\right\}
\end{multline}
for all $t\in(0,\infty)$, $0<z\le\frac{1}{t}
\sqrt{\frac{a(\theta)}{b(\theta)}} \wedge 2\eps\eqqcolon
z_*(t,\theta,\eps),$ and $\theta\in(0,2\pi]$. Note that the function
$$
k(\tau,u,\theta)= a - 4\tau u^{-1}\left(a + bu^2\right)
$$
attains its global maximum with respect to $u\in(0,\infty)$ for
every fixed $\tau>0$ and $\theta\in(0,2\pi)$  at the point
$u=\sqrt{\frac{a(\theta)}{b(\theta)}}\eqqcolon u_*(\theta)$ with
$$
\max_{u>0}k(\tau,u,\theta)=k\left(\tau,u_*(\theta),\theta\right)=a(\theta) - 8\tau\sqrt{a(\theta)b(\theta)},\quad \tau>0,\quad\theta\in(0,2\pi),
$$
being monotonically increasing for $0<u<u_*(\theta)$ and
monotonically decreasing for $u>u_*(\theta)$ (in fact,
$k(\tau,u,\theta)$ is concave with respect $u>0$ for every fixed
$\tau>0$ and $\theta\in(0,2\pi)$). Now choosing $z=z_*$ to minimize
the right-hand side of~\eqref{tech.LemMod1.1} with respect to
$z\in(0,z_*]$ (which is equivalent to the choice of
$u=u_*(\theta)\wedge 2\eps t$) and optimizing then with respect to
$\theta\in(0,2\pi]$ we arrive at~\eqref{AbsCh.F.EstimComp} with the
right-hand side being a monotonically increasing function of $L>0$
as the least upper bound to a family of monotonically decreasing
functions.

Let us prove that~\eqref{AbsCh.F.EstimComp}
(or~\eqref{tech.LemMod1.1}) yields~\eqref{AbsCh.F.EstimAnalytically}
under the assumptions stated in the formulation of  part~(ii) of the
theorem. For this purpose, observe that for every
$\theta\in(0,2\pi)$
$$
-k(L^3t,tz_*,\theta)+a=
\begin{cases}
2aL^3\eps^{-1}+8bL^3\eps t^2\le 2aL^3\eps^{-1}+4L^3 t\sqrt{ab}, &  2\eps t\le\sqrt{a/b},
\\
8L^3t\sqrt{ab}, & 2\eps t>\sqrt{a/b},
\end{cases}
$$
so that for all $t\in\R$ we have
$$
-k(L^3t,tz_*,\theta)+a\le \max\left\{2aL^3\eps^{-1}+4L^3t\sqrt{ab}, 8L^3t\sqrt{ab}\right\}.
$$
Now let $\tau$ be an arbitrary positive number. Then for
$0<t\le\tau/L^3$ we trivially obtain
$$
-k(L^3t,tz_*,\theta)+a\le \max\left\{2aL^3\eps^{-1}+4\tau\sqrt{ab}, 8\tau\sqrt{ab}\right\}=8\tau\sqrt{ab}, \quad\text{if }\ \tau\ge\frac{L^3}{2\eps}\sqrt{\frac{a}{b}},
$$
and if for $L_0\ge L$ the uniform condition
$$
\tau\ge \frac{L_0^3}{2\eps}\sup_{0<\theta<2\pi}\sqrt{\frac{a(\theta)}{b(\theta)}} =\frac{L_0^3}{2\eps}\lim_{\theta\to2\pi-}\sqrt{\frac{a(\theta)}{b(\theta)}}=\frac{\pi L_0^3}{\eps}\eqqcolon\underline\tau_1(L_0,\eps)
$$
holds, where the equality sign is due to
Lemma~\ref{Lem:a(theta)/b(theta)increases}, then also
$$
-k(L^3t,2\eps t)\coloneqq -\sup_{0<\theta<2\pi}k(L^3t,tz_*(t,\theta,\eps),\theta)\le -\sup_{0<\theta<2\pi}\left\{a(\theta)-8\tau\sqrt{a(\theta)b(\theta)}\right\}\eqqcolon -k(\tau)
$$
for every $t\le\tau/L^3$ and $L\ge L_0.$ Note that for fixed  $\eps$
the quantity $\underline{\tau}_1(L_0,\eps)$ can be made arbitrarily
small by the choice of $L_0=L_0(\eps)>0$ small enough.

Finally, observe that the function $h(\tau,\theta)\coloneqq
a(\theta) - 8\tau\sqrt{a(\theta)b(\theta)},$ $\tau>0,$
$\theta\in(0,2\pi)$ is strictly positive if and only if
$\tau<\frac18\sqrt{\frac{a(\theta)}{b(\theta)}}\eqqcolon\tau_1(\theta)$,
and hence,
$k(\tau)\coloneqq\sup\limits_{0<\theta<2\pi}h(\tau,\theta)>0$ if and
only if there exists $\theta\in(0,2\pi)$ such that
$\tau<\tau_1(\theta).$ The latest condition is equivalent to $\tau<
\sup\limits_{0<\theta<2\pi}
\tau_1(\theta)=\pi/4\eqqcolon\overline\tau_1$. The proof is
completed by the remark that the set of admissible values of
$\tau\in\left[\underline{\tau}_1(L_0,\eps), \overline\tau_1\right)$
is not empty if and only if
$$
\underline{\tau}_1(L_0,\eps)\coloneqq\pi L_0^3/\eps<\pi/4\eqqcolon\overline\tau_1,
$$
i. e. $L_0<(\eps/4)^{1/3}.$
\end{proof}

\subsection{Estimation of $\big|\overline f_n(t)-e^{-t^2/2}\big|$.}

The present section is aimed at bounding the absolute value of the
difference $\big|\overline f_n(t)-e^{-t^2/2}\big|$ between the ch.f.
of the normalized sum and that of the standard normal law. Similarly
to the preceding section, we construct two estimates, the first
being the most exact and used in Subsection~\ref{SubsecCaseLargeL}
and the second being the most convenient for analytical work in
Subsection~\ref{SubsecCaseSmallL}.

However, unlike in the previous section, here significant
distinctions arise in the dependence of which of the fractions
$\es^3$ or $\roz^3$ is used. We pay special attention to the
appearing distinctions and explain their reasons.

Before passing to the main theorem of the present subsection we
proof four auxiliary statements.

\begin{lemma}\label{LemKappa2}
For every $x \in \R$ we have
$$
\abs{e^{ix} - 1 - ix - (ix)^2 / 2} \le \varkappa\cdot x^2,
$$
where $\varkappa := x^{-2}\sqrt{\vphantom{\frac12}(\cos x-1+x^2/2)^2+(\sin x-x)^2}\,\Big\vert_{x=x_0} = 0.531551\ldots$ and $x_0=5.487414\ldots$ is the unique root of the equation
$$
8(\cos x - 1) + 8x\sin x  - 4x^2\cos x-x^3\sin x = 0
$$
on the interval $(0,2\pi]$ lying in $(\pi,2\pi)$.
\end{lemma}

\begin{proof}
Due to the symmetry and triviality of the stated inequality for
$x=0$ we may assume that $x>0.$ Observe that the inequality of
interest is equivalent to the relation
$$
\max_{x>0}g(x)=g(x_0)=0.2825\ldots=\varkappa^2,\quad \text{where }\ g(x) := \frac{\left(\cos x - 1 + x^2/2\right)^2+\left(\sin x - x \right)^2}{x^4},\quad x>0.
$$
We have
$$
h(x) := x^5 g'(x) = 8(\cos x - 1) + 8x\sin x  - 4x^2\cos x-x^3\sin x,\quad h'(x) = x^2\left(\sin x - x\cos x\right).
$$
Note that the function $h'(x)$ has a unique root $x_1\in(\pi,
3\pi/2)$ on  $(0,2\pi]$, which is the point of maximum of $h(x).$
Since, $h(0)=0$ and $h$ is strictly increasing on $(0,x_1)$, we
conclude that $h(x_1)>0.$ Moreover, $h(2\pi)=-16\pi^2<0,$ hence, $h$
has a unique sign change on $(0,2\pi]$ and this sign change occurs
in the point $x_0\in(x_1,2\pi)\subset(\pi,2\pi)$ from $+$ to $-$.
Thus, the function $g$ has a unique stationary point
$x_0\in(\pi,2\pi)$ on $(0,2\pi],$ which is the point of maximum. So,
it remains to prove that $g(x)\le g(x_0)$ for $x>2\pi.$

If $x=2\pi+y\in(2\pi,9\pi/4]$, then, with the account of the
monotone increase of $\sin y$ and monotone decrease of $\cos y$ for
$y\in(0,\pi/4]$, we have
$$
\frac{h'(x)}{x^2}=\sin y-(y+2\pi)\cos y\le\sin\frac\pi4-2\pi\cos\frac{\pi}4=\frac{\sqrt2}{2}(1-2\pi)<0.
$$
Hence, $h(x)$ is strictly monotonically decreasing for $x\in(2\pi,9\pi/4]$ with
$$
\sup_{x\in(2\pi,9\pi/4]}h(x)=h(2\pi+)=-16\pi^2<0,
$$
so that $g'(x)\coloneqq x^{-5}h(x)<0$ for $x\in(2\pi,9\pi/4]$  and
$$
\sup_{x\in(2\pi,9\pi/4]}g(x)=g(2\pi+)=\tfrac14\left(1+\pi^{-2}\right)=0.2753\ldots<0.2825\ldots=\varkappa^2.
$$

Finally, for $x>9\pi/4$ we use the trivial bound
$$
g(x)=\frac{2(1-\cos x)+x^4/4+x^2\cos x-2x\sin x}{x^4}<\left(\frac{1}{4} + \frac{1}{x^2} + \frac{2}{x^3}\right)\bigg\vert_{x=9\pi/4} < 0.2757 < \varkappa^2,
$$
which completes the proof of the lemma.
\end{proof}

\begin{lemma}\label{LemAlpha2}
For $\eps > 0$ denote
\[
\alpha(\eps):=\inf_{0<x<\eps\wedge1}\left(x-x^3\right)^{-2/3} =
\begin{cases}
\left(\eps-\eps^3\right)^{-2/3},&\text{if } \eps\le 3^{-1/2},\\
3\cdot2^{-2/3}=1.8898\ldots,&\text{otherwise}.
\end{cases}
\]
Then for every $\eps>0$ we have
\begin{equation*}\label{LemAlpha2_Ineq1}
\max_{1\le k\le n}\sigma_k^2\le \alpha(\eps)\sup_{0<x\le\eps}\big(x\lind(x)\big)^{2/3}.
\end{equation*}
\end{lemma}

\begin{proof}
From the definition of $\lind(z)\coloneqq\sum_{k=1}^n\qt_k(z)$ it follows that
\[
\sup_{0<x\le\eps}x\lind(x)\ge z\qt_k(z)
\]
for every ${z\in(0,\eps]}$ and $k=1,\ldots,n$. On the other hand, by the definition of $\sigma_k^2\coloneqq\E X_k^2$, we have
\[
\sigma_k^2 = \E X_k^2\I(|X_k| < z) + \E X_k^2\I(|X_k|\ge z)\le z^2 + \qt_k(z),\quad z>0,\quad k=1,\ldots,n,
\]
whence it follows that
$$
\sup_{0<x\le\eps}x\lind(x)\ge z\left(\sigma_k^2 - z^2\right)\quad\text{for }\ 0 < z \le \eps.
$$
Now the choice of $z\coloneqq\sigma_k/\sqrt{3} \wedge \eps
\sigma_k\le\eps$ as a maximizer to the right-hand side of the latest
inequality completes the proof of the lemma.
\end{proof}

\begin{lemma}\label{LemDist2}
For all $t\in\R$ and $0<z\le\eps\le\infty$ we have
$$
\sum_{k=1}^n \abs{f_k(t) - 1}^2
\le t^4\bigg[\frac{z}{2}\sup_{0 < x\le\eps}x\lind(x) +\frac{1}{4}\sum_{k=1}^n\sigma_k^2\qt_k(z)\bigg]\le t^4\bigg[\frac{z}{2}\sup_{0<x\le\eps}x\lind(x) +\frac{\alpha(\eps)}{4z}\sup_{0<x\le\eps}(x\lind(x))^{5/3}\bigg].
$$
\end{lemma}

\begin{proof}
For every random summand $X$ with the d.f. $F$, the ch.f. $f$ and $\E X=0,$
$\qt(z)\coloneqq\E X^2\I(|X|\ge z),$ $z\ge0,$ $\sigma^2\coloneqq\qt(0)=\E X^2$ we have
\[
\abs{f(t)-1}^2 \le \Big(\frac{t^2\sigma^2}2\Big)^2
=\frac{t^4}{4} \bigg[\int_{|x| < z} x^2 dF(x) + \qt(z) \bigg]^2=
\]
\[
=\frac{t^4}{4} \bigg[\bigg(\int_{|x| < z} x^2 dF(x)\bigg)^2 +2\qt(z)\int_{|x|<z}x^2dF(x)+(\qt(z))^2\bigg],\quad z\ge0.
\]
Using the  Jensen inequality for the first term, both the
representation $\int_{|x|<z}x^2dF(x)=\sigma^2-\qt(z)$ and the bound
$\int_{|x|<z}x^2dF(x)\le z^2$ for the second term in square brackets
in the first step below, and then equality~\eqref{Destr4Moment} we
obtain
\[
\abs{f(t)-1}^2\le\frac{t^4}{4} \bigg[\int_{|x| < z} x^4 dF(x) +z^2\qt(z) +\sigma^2\qt(z)\bigg]
=\frac{t^4}{4}\bigg[2\int_0^z x\qt(x)dx + \sigma^2\qt(z)\bigg].
\]

Summing up the constructed bounds for every random summand
$X\coloneqq X_k$ over all $k=1,\ldots,n,$ we get
$$
\sum_{k=1}^n\abs{f_k(t)-1}^2\le \frac{t^4}{4}\bigg[2\int_0^z x\lind(x)dx +\sum_{k=1}^n\sigma_k^2\qt_k(z)\bigg]\le t^4\bigg[\frac{z}{2}\sup_{0<x\le z} x\lind(x) +\frac14\sum_{k=1}^n\sigma_k^2\qt_k(z)\bigg].
$$
Now the observation that the least upper bound with respect to
$x\in(0,z]$ is a non-decreasing function of $z$ yields the first
claim of the lemma. The second claim follows from
Lemma~\ref{LemAlpha2} yielding the chain of inequalities
$$
\sum_{k=1}^n\sigma_k^2\qt_k(z)\le \alpha(\eps)\sup_{0<x\le\eps}\big(x\lind(x)\big)^{2/3}L_n(z)
\le \frac{\alpha(\eps)}{z}\sup_{0<x\le\eps}\big(x\lind(x)\big)^{5/3},\quad z\in(0,\eps].
$$
\end{proof}

\begin{lemma}\label{LemDist1}
For all $t\in\R,$ $\eps>0,$ and $\gamma > 0$ we have
\begin{equation}\label{LemDist1_Ineq1}
\bigg|\sum\limits_{k=1}^{n} \left(f_k(t) - 1 +\sigma_k^2 t^2 / 2\right)\bigg|\le p_\scE\left(t,\eps,\gamma\right) \cdot\es^3(\eps,\gamma),
\end{equation}
\begin{equation}\label{LemDist1_Ineq2}
\bigg|\sum\limits_{k=1}^{n} \left(f_k(t) - 1 +\sigma_k^2 t^2 / 2\right)\bigg|\le p_\scR(t,\eps,\gamma)\cdot\roz^3(\eps,\gamma),
\end{equation}
where
\begin{equation}\label{pE(t,eps,gamma)def}
p_\scE(t,\eps,\gamma) \coloneqq t^2\min_{0<z\le\eps}\left[\frac{zt^2}{12}+\max\left\{
\frac{|t|}{6\gamma},\frac{\varkappa}{z} - \frac{zt^2}{24}\right\}\right]
=\begin{cases}
\frac{\varkappa t^2}{\eps}+\frac{\eps t^4}{24}, &\eps|t|\le t_\gamma, \\[2mm]
\frac{\sqrt{6\varkappa\gamma^2+1}}{6\gamma}\,|t|^3, &\eps|t|>t_\gamma,
\end{cases}
\end{equation}
\begin{multline}\label{pR(t,eps,gamma)def}
p_\scR(t,\eps,\gamma)
\coloneqq t^2\max\left\{\frac{\abs{t}}{6\gamma},\frac{\varkappa}{\eps} + \frac{\eps t^2}{24},\frac{\eps t^2}{12}\right\}=
\\
=\begin{cases}
 \begin{cases}
    \frac{\varkappa t^2}{\eps}+\frac{\eps t^4}{24}, &\eps|t|\le t_\infty,
    \\
    \frac{\eps t^4}{12}, &\eps|t|>t_\infty,
    \end{cases}
    & \gamma \ge \gopt,
\\[3mm]
 \begin{cases}
    \frac{\varkappa t^2}{\eps}+\frac{\eps t^4}{24}, &\eps|t|\le t_{1,\gamma},
    \\[1mm]
    \frac{|t|^3}{6\gamma},&t_{1,\gamma} < \eps|t| \le t_{2,\gamma},
    \\[1mm]
    \frac{\eps t^4}{12}, &\eps|t|>t_{2,\gamma},
 \end{cases}
&\gamma < \gopt,
\end{cases}
=p_\scR\big(t,\eps,\gamma\wedge\gopt\big),
\end{multline}
$$
t_\gamma\coloneqq \tfrac2\gamma\big(\sqrt{(\gamma/\gopt)^2+1}-1\big), \quad
t_\infty\coloneqq\lim\limits_{\gamma \to \infty} t_\gamma
=2/\gopt=3.5717\ldots,
$$
$$
t_{2,\gamma}\coloneqq2\max\big\{\gamma^{-1},\gopt^{-1}\big\},\quad
t_{1,\gamma}\coloneqq t_{2,\gamma} \big(1-\sqrt{(1-(\gamma/\gopt)^2)_+}\,\big),
$$
$\gopt={1}/{\sqrt{6\varkappa}}=0.5599\ldots$ and $\varkappa$ is defined in Lemma~\ref{LemKappa2}. Moreover, the functions $t_\gamma$ and  $t_{1,\gamma}$  are monotonically increasing with respect to $\gamma,$  and $t_\gamma\le  t_{1,\gamma}\le t_\infty$ for all $\gamma>0$, the functions  $p_\scE(t,\eps,\gamma),$ $p_\scR(t,\eps,\gamma)$ are monotonically decreasing with respect to $\gamma>0$ with $p_\scR(t,\eps,\gamma)$ do not depending on $\gamma$
for $\gamma \ge\gopt,$ and $p_\scE(t,\eps,\gamma)$ is also monotonically decreasing with respect to $\eps.$
\end{lemma}

The values of the function $t_\gamma
={2}{\gamma}^{-1}\big(\sqrt{(\gamma/\gopt)^2+1}-1\big)$ for some
$\gamma$ are given in Table~\ref{Tab:t0(eps,gamma)}. The plots of
the functions $t_\gamma$ and $t_{1,\gamma}$ are given on
Fig.\,\ref{Fig:t_gamma+tilde(t_gamma)-plots} in the Appendix.

\begin{table}[h]
\begin{center}
\begin{tabular}{|c|c|c|c|c|c|c|c|c|c|}
\hline
$\gamma\vphantom{\displaystyle\frac12}$
& $0+$
& $0.25$
& $0.41$
&$\gopt$
& $0.73$
& $1$
& $1.25$
&$1.5$
&$\infty$
\\ \hline
$t_\gamma\vphantom{\displaystyle\frac12}$
&$0$
&$0.7611$
&$1.1678$
&$1.4794$
&$1.7617$
&$2.0935$
&$2.3137$
&$2.4791$
&$3.5717$
\\ \hline
\end{tabular}
\caption{Values of the function $t_\gamma={2}{\gamma}^{-1}\big(\sqrt{(\gamma/\gopt)^2+1}-1\big)$ for some $\gamma$  (rounded down).}
\label{Tab:t0(eps,gamma)}
\end{center}
\end{table}

\begin{remark}
For $\gamma\ge\gopt$ we have $t_{1,\gamma}=t_{2,\gamma}=2/\gopt=t_\infty$.
\end{remark}

\begin{remark}\label{Rem:p_E=p_R,|t|<=t_0}
The functions $p_\scE(t,\eps,\gamma)$ and $p_\scR(t,\eps,\gamma)$  coincide on the interval $|t|\le t_\gamma/\eps$ for every ${\gamma>0}$.
\end{remark}

\begin{proof}
Since $\E X_k=0$ for all $k=1,\ldots,n$, for every $t\in\R$ and $z>0$ we have
\begin{multline*}
I\coloneqq\abs{\sum_{k=1}^{n} \left(f_k(t) - 1 + \sigma_k^2 t^2 / 2 \right)} = \abs{\sum_{k=1}^{n}\int_{\R} \left(
e^{itx}-1-itx- \tfrac{(itx)^2}{2}\right)dF_k(x)}\le
\\
\le\sum_{k=1}^{n}\int_{\abs{x} < z} \abs{e^{itx}-1-itx-\frac{(itx)^2}{2}-\frac{(itx)^3}{6}}dF_k(x)+
\\
+\abs{\sum_{k=1}^{n}\int_{\abs{x} < z} \frac{(itx)^3}{6}dF_k(x)} +\sum_{k=1}^{n}\int_{\abs{x} \ge z} \abs{
e^{itx}-1-itx-\frac{(itx)^2}{2}}dF_k(x).
\end{multline*}
Using the inequalities $|e^{iy} - 1 - iy - \frac{1}{2}(iy)^2 -\frac{1}{6}(iy)^3| \le \frac{1}{24} y^4$ and $|e^{iy} - 1 - iy - \frac{1}{2}(iy)^2| \le \varkappa y^2$ valid for all $y\in\R$ and with the account of~\eqref{Destr4Moment} we obtain
\begin{multline}\label{tech.LemDist1.1}
I\le\frac{t^4}{24}\sum_{k=1}^{n}\int_{\abs{x} < z}x^4dF_k(x) +\frac{\abs{t}^3}{6}\abs{M_n(z)} +\varkappa t^2\sum_{k=1}^{n}\qt_k(z)=
\\
=\frac{t^4}{12} \int_{0}^{z} x\lind(x) dx -\frac{z^2t^4}{24}\lind(z)
+\frac{\abs{t}^3}{6}\abs{M_n(z)} +\varkappa t^2\lind(z)\le
\\
\le\frac{zt^4}{12} \sup_{0<x\le z} x\lind(x)
+\left(\frac{\varkappa t^2}z-\frac{zt^4}{24}\right)z\lind(z)
+\frac{\abs{t}^3}{6}\abs{M_n(z)},\quad z>0,\ t\in\R.
\end{multline}

To prove~\eqref{LemDist1_Ineq1}, observe that~\eqref{tech.LemDist1.1} yields
\[
I\le\frac{zt^4}{12}\sup_{0 < x\le z} x\lind(x)
+\max\left\{\frac{\varkappa t^2}{z}-\frac{zt^4}{24},
\frac{|t|^3}{6\gamma}\right\}
\left(z\lind(z)+\gamma\abs{M_n(z)}\right)
\]
$$
\le  t^2\left(\frac{zt^2}{12}+\max\left\{
\frac{\varkappa}{z} - \frac{zt^2}{24},\frac{|t|}{6\gamma}\right\}\right) \sup_{0<x\le z} \left\{\gamma\abs{M_n(x)}+x\lind(x)\right\}
$$
for all $t\in\R$ and $z>0.$ Now choosing $z=\frac{2}{\gamma|t|}\left(\sqrt{6\varkappa \gamma^2 + 1} - 1\right)\wedge\eps$ to minimize the R.H.S. of the expression in the large brackets with respect to $z\in(0,\eps]$ and observing that
$$
\sup_{0<x\le z} \left\{\gamma\abs{M_n(x)}+x\lind(x)\right\}=\es(z,\gamma)\le\es(\eps,\gamma)
$$
for every $0<z\le\eps$ and $\gamma>0$, we arrive at~\eqref{LemDist1_Ineq1}.

To prove~\eqref{LemDist1_Ineq2}, observe that~\eqref{tech.LemDist1.1} yields
$$
I\le\left(\frac{zt^4}{12} +\left(\frac{\varkappa t^2}z-\frac{zt^4}{24}\right)_+\right)
\sup_{0 < x\le z} x\lind(x)+\frac{\abs{t}^3}{6}\abs{M_n(z)} \le
$$
$$
\le \max\left\{\frac{zt^4}{12} +\Big(\frac{\varkappa t^2}z-\frac{zt^4}{24}\Big)_+,
\frac{\abs{t}^3}{6\gamma}\right\} \left(\sup_{0 < x\le z} x\lind(x)+\gamma\abs{M_n(z)} \right)=p_\scR(t,z,\gamma)\roz^3(z,\gamma)
$$
for every $z>0$ and $t\in\R$. The choice of $z=\eps$ completes the
proof of inequality~\eqref{LemDist1_Ineq2}. The stated properties of
the functions $t_\gamma$,  $\widetilde t_\gamma$, $p_{\scE}$, and
$p_{\scR}$, as well as representations for $p_{\scE}$ and $p_{\scR}$
in the right-hand sides of~\eqref{pE(t,eps,gamma)def}
and~\eqref{pR(t,eps,gamma)def} are proved by use of elementary
analysis.
\end{proof}

Now we are ready to prove the main result of the present subsection.

\begin{theorem}\label{ThAbsChFDiffEstim}
For $\eps > 0,\ \gamma>0,\ L >0,\ \tau>0$ denote
\[
\alpha(\eps):=\inf_{0<x<\eps\wedge1}\left(x-x^3\right)^{-2/3} =
\begin{cases}
\left(\eps-\eps^3\right)^{-2/3},&\text{if } \eps\le 1/\sqrt3=0.5773\ldots,\\
3\cdot2^{-2/3}=1.8898\ldots,&\text{otherwise},
\end{cases}
\]
$$
\textstyle
\overline\tau_0(\eps)\coloneqq\sqrt{\frac{2}{\alpha(\eps)}},\qquad
\overline L_0(\eps)\coloneqq\eps\sqrt{\frac2{\alpha(\eps)}}=\eps\overline\tau_0(\eps),
$$
$$
\alpha_1(\eps,L) = \min_{0<u\le\eps L^{-1}}\Big(\frac{u}{2}+\frac{\alpha(\eps)}{4u}\Big)=
\begin{cases}
\sqrt{\alpha(\eps)/2},& L\le \overline L_0(\eps),
\\[2mm]
\frac{\eps}{2L}+\frac{\alpha(\eps)}{4\eps}L,& L>\overline L_0(\eps),
\end{cases}
$$
$$
\textstyle
\varkappa= \sup\limits_{x>0}x^{-2}\sqrt{\vphantom{\frac12}(\cos x-1+x^2/2)^2+(\sin x-x)^2} = 0.5315\ldots,\quad \gopt=1/\sqrt{6\varkappa}=0.5599\ldots\ .
$$
For $\eps  > 0,\ \gamma>0,\ L>0,\ 0<\tau< \overline\tau_0(\eps)$  also introduce the functions
$$
B(\tau,\eps,L)= -\frac{4\alpha_1(\eps,L)}{\alpha^2(\eps)\tau^4} \Big[ \ln\left(1-\tfrac12\alpha(\eps)\tau^2\right)+\tfrac12\alpha(\eps)\tau^2\Big],\quad B(\tau,\eps)= \frac{\sqrt{\alpha(\eps)/2}}{2-\alpha(\eps)\tau^2},
$$
$$
A(\tau,\eps)=\exp\left\{\tau^4B(\tau,\eps)\right\},\quad A_{\scE}(\tau,\eps,\gamma,L)=
A(\tau,\eps) \exp\left\{\frac{\tau^3}{6\gamma} +\sqrt{\frac{\varkappa}3}\tau^3 +\frac{\varkappa}\eps L\tau^2 \right\},
$$
$$
A_{\scR}(\tau,\eps,\gamma,L) =A(\tau,\eps) \exp\left\{\frac{\tau^3}{6\gamma}\I(\gamma<\gopt) +\frac{\varkappa}\eps L\tau^2 \right\},
$$
\noindent{\rm\bf(i)} For every $L\in\{\es(\eps,\gamma),\roz(\eps,\gamma)\},$ $\eps>0,$ $\gamma>0,$ and $L|t|<\overline{\tau}_0(\eps)$  we have
\begin{equation}\label{AbsChFDiffEstimComp}
r_n(t)\coloneqq\big|\Sfn - e^{-t^2/2}\big| \le \left( \exp\left\{ L^3p(t,\eps,\gamma) +L^4t^4B(L|t|,\eps,L) \right\} -1\right) e^{-t^2/2},
\end{equation}
where $p=p_\scE$ for $L=\es(\eps,\gamma)$ and $p=p_\scR$ for $L=\roz(\eps,\gamma)$ with
$p_\scE,$ $p_\scR$ defined in Lemma~\ref{LemDist1}. In particular, for $L=\es(\eps,\gamma)$ with $\eps=\infty$ we have
$$
r_n(t)  \le \left( \exp\left\{\tfrac{\sqrt{6\varkappa\gamma^2+1}}{6\gamma} L^3 -4\sqrt2\cdot3^{-3/2}\Big[ \ln\left(1-3\cdot2^{-5/3}\tau^2\right)+3\cdot2^{-5/3}\tau^2\Big]\tau^{-4} t^4L^4 \right\} -1\right) e^{-t^2/2}
$$
for $L|t|<2^{5/6}/\sqrt3=1.0287\ldots\,.$ Moreover, the right-hand side of~\eqref{AbsChFDiffEstimComp} is
monotonically increasing with respect to $L>0$ and monotonically decreasing with respect to $\gamma>0$ and,
for $L=\es$, also with respect to $\eps>0.$  The right-hand side of~\eqref{AbsChFDiffEstimComp} with
$L=\roz(\eps,\gamma)$ does not depend on $\gamma$ for $\gamma\ge\gopt$.

\noindent{\rm\bf(ii)} For $L=\es(\eps,\gamma) \le \overline{L}_0(\eps)$  and $L|t|<\tau\le\overline{\tau}_0(\eps),$ we have
\begin{equation}\label{AbsChFDiffEstimAnalytEsseen}
r_n(t) \le
A_{\scE}(\tau,\eps,\gamma,L)\left( p_{\scE}(t,\eps,\gamma)+B(\tau,\eps)Lt^4\right) \cdot L^3e^{-t^2/2},
\end{equation}
in particular, with $\eps=\infty,$
$$
r_n(t)\le  \exp\Big\{\tfrac{\tau^3}{6\gamma} +\sqrt{\tfrac{\varkappa}6}\tau^3 +\tau^4B(\tau)\Big\} \Big(\tfrac{\sqrt{6\varkappa\gamma^2+1}}{6\gamma}+  B(\tau)L|t|\Big) \cdot  L^3|t|^3e^{-t^2/2},\quad  t\in\R,
$$
where $B(\tau)\coloneqq B(\tau,\infty) =\left(
\sqrt3\cdot2^{-1/6}\right) /\left(2^{5/3}-3\tau^2\right)$. Moreover,
the right-hand side of~\eqref{AbsChFDiffEstimAnalytEsseen} is
monotonically increasing with respect to $L>0$ and monotonically
decreasing with respect to $\eps>0$ and $\gamma>0$.

\smallskip

For $L=\roz(\eps,\gamma) \le \overline{L}_0(\eps)$ and $L|t|<\tau\le\overline{\tau}_0(\eps)$ with $B\coloneqq B(\tau,\eps),$  we have
\begin{equation}\label{AbsChFDiffEstimAnalytRozovsky}
r_n(t) \le
A_{\scR}(\tau,\eps,\gamma,L)\left( p_{\scR}(t,\eps,\gamma)+B(\tau,\eps)Lt^4\right) \cdot L^3e^{-h(\tau,\eps,L)t^2/2},
\end{equation}
where $h(\tau,\eps,L)\coloneqq1-{\eps\tau^2L}/{6}>0$ for all
$\tau\in\big(0,\overline\tau_0(\eps)\big)$ if and only if
$L\le6/(\eps\overline\tau_0^2(\eps))$. Moreover, the right-hand side
of~\eqref{AbsChFDiffEstimAnalytRozovsky} is monotonically increasing
with respect to $L>0$ and monotonically decreasing with respect to
$\gamma>0$.
\end{theorem}

\begin{remark}\label{ThDistance_increasing}
The right-hand sides of~\eqref{AbsChFDiffEstimComp} with
$L=\roz(\eps,\gamma)$ and~\eqref{AbsChFDiffEstimAnalytRozovsky}  are
unbounded as $\eps\to\infty$.
\end{remark}

\begin{proof}
Fix arbitrary $\eps,\gamma>0$ and let $L=\es(\eps,\gamma)$ or $L=\roz(\eps,\gamma)$. With $u_k\coloneqq f_k\left(t\right) - 1$ Taylor's formula and Lemma~\ref{LemAlpha2} yield
\begin{equation}\label{|f_k(t)-1|<=alpha2t^2L^2/2}
\abs{u_k} =\abs{f_k\left(t\right) - 1}\le \frac{t^2 \sigma_k^2}{2} \le \frac{\alpha(\eps)}{2}L^2 t^2
\end{equation}
for every $k=1,\ldots,n$ and $t \in \R,$ whence it follows that $|u_k|<1$ for $L|t|<\sqrt{2/\alpha(\eps)}\eqqcolon\overline{\tau}_0(\eps)$ and that $\Sfn = \prod_{k=1}^n (1+u_k)$ does not vanish for $L|t|<\overline{\tau}_0(\eps)$. Hence, the logarithm $\ln{\Sfn}$ is well-defined for $L|t|<\overline{\tau}_0(\eps)$. Considering the main branch of the logarithm and using the inequality $\abs{e^z-1}\le e^{\abs{z}} - 1$, $z \in \C$, we get
\begin{equation}\label{|f_n-e|<=(e-1)e}
r_n(t)=\abs{\exp\left\{\frac{t^2}{2} + \ln\Sfn\right\}-1}e^{-t^2/2}
\le\big(e^{\delta_n(t)} - 1\big)e^{-t^2/2}\le\delta_n(t)e^{\delta_n(t)} e^{-t^2/2},\quad L|t|<\overline{\tau}_0(\eps),
\end{equation}
where
\begin{multline*}
\delta_n(t) \coloneqq\abs{\frac{t^2}{2} + \ln \Sfn}=
\abs{\sum\limits_{k=1}^n\left(\ln\left(1 + u_k\right) + \frac{\sigma_k^2t^2}{2} \right)}
\le \sum_{k=1}^n\abs{\ln\left(1 + u_k\right) - u_k}+ \abs{\sum_{k=1}^n \left(u_k + \frac{\sigma_k^2t^2}{2}\right)}\le
\\
\le\sum_{k=1}^n\sum_{j=2}^\infty\frac{|u_k|^j}{j} + \abs{\sum_{k=1}^n \left(f_k(t)-1+ \frac{\sigma_k^2t^2}{2}\right)}\eqqcolon I_2+I_1.
\end{multline*}
By Lemma~\ref{LemDist1}  we have
$$
I_1\le p(t,\eps,\gamma)\cdot L^3\quad \text{for}\quad t\in\R,\ \eps>0,\ \gamma>0,
$$
where $p=p_\scE$ if $L=\es(\eps,\gamma)$ and $p=p_\scR$ if $L=\roz(\eps,\gamma)$.

To bound $I_2$ above, observe that
inequality~\eqref{|f_k(t)-1|<=alpha2t^2L^2/2} yields
$$
I_2
\le \sum_{k=1}^n|u_k|^2\sum_{j=2}^\infty\frac1j\left(\frac{\alpha(\eps)}2L^2t^2\right)^{j-2} =
-\frac{4}{\alpha^2(\eps)L^4t^4} \left[\ln\left(1-\frac{\alpha(\eps)}2L^2t^2\right)+\frac{\alpha(\eps)}2L^2t^2\right]\sum_{k=1}^n|u_k|^2
$$
for $L|t|<\overline{\tau}_0(\eps),$ while Lemma~\ref{LemDist2}, with
the account of the inequality $\sup_{0 < x\le\eps}x\lind(x)\le L^3,$
implies that
$$
\sum_{k=1}^n|u_k|^2 = \sum_{k=1}^n \abs{f_k(t) - 1}^2
\le t^4 \bigg[\frac{z}{2}L^3 +\frac{\alpha(\eps)}{4z}L^5\bigg] \eqqcolon t^4L^4\cdot\widetilde \alpha_1\left(\eps,zL^{-1}\right),\quad t\in\R,\ z\in[0,\eps],
$$
where
$$
\widetilde \alpha_1(\eps,u)=\frac{u}2 + \frac{\alpha(\eps)}{4u},\quad \eps>0,\ u>0.
$$
Choosing $u=zL^{-1}=\sqrt{{\alpha(\eps)}/{2}}\wedge\eps
L^{-1}\coloneqq u_*(L,\eps)$ to  minimize $\widetilde
\alpha_1(\eps,u)$ with respect to $u\in\left(0,\eps L^{-1}\right]$
we obtain $\sum_{k=1}^n\abs{u_k}^2\le\alpha_1(\eps,L)\cdot L^4t^4$
with $\alpha_1(\eps,L)\coloneqq\min\limits_{0<u\le\eps
L^{-1}}\widetilde \alpha_1(\eps,u)=\alpha_1(\eps, u_*(L,\eps))$ as
in the formulation of the lemma. Thus, for
$L|t|<\overline{\tau}_0(\eps)$ we have
$$
I_2\le -\frac{4\alpha_1(\eps,L)}{\alpha^2(\eps)} \Big[ \ln\left(1-\tfrac12\alpha(\eps)L^2t^2\right)+\tfrac12\alpha(\eps)L^2t^2
\Big] \eqqcolon B(L|t|,\eps,L)\cdot L^4t^4,
$$
and
\begin{equation}\label{delta_n(t)<=pL^3+BL^4t^4}
\delta_n(t)\le I_1+I_2\le p(t,\eps,\gamma)\cdot L^3 +B(L|t|,\eps,L)\cdot L^4t^4,
\end{equation}
whence, with the account of the penultimate inequality
in~\eqref{|f_n-e|<=(e-1)e}, we obtain~\eqref{AbsChFDiffEstimComp}.
The observation that $p(t,\eps,\gamma)$ is monotonically decreasing
with respect to $\gamma$, the functions $p_{\scE}(t,\eps,\gamma)$,
$\alpha(\eps),$ $\alpha_1(\eps,L)$ are monotonically decreasing with
respect to $\eps$ and
$$
B(L|t|,\eps,L)=\alpha_1(\eps,L) \sum_{j=0}^\infty\frac{(\alpha(\eps)L^2t^2/2)^j}{j+2}
$$
is monotonically increasing with respect to $L>0$ and monotonically
decreasing with respect to $\eps>0$ yields the stated properties of
the right-hand side of~\eqref{AbsChFDiffEstimComp}.

Now let us assume that $L$ is small and prove slightly rougher
bounds than~\eqref{AbsChFDiffEstimComp} that are more convenient for
analytical work. By virtue of the inequality
$-(\ln\left(1-x\right)+x) =\sum_{j=2}^\infty{x^j}/{j}\le
\frac{1}{2}\cdot \frac{x^2}{1-x},$ $0\le x<1,$ applied to
$x\coloneqq\alpha(\eps)L^2t^2/2$ we have
$$
B(L|t|,\eps,L)\le \frac{\alpha_1(\eps,L)}{2-\alpha(\eps)L^2t^2}\le \frac{\alpha_1(\eps,L)}{2-\alpha(\eps)\tau^2}\quad\text{for}\quad L|t|<\tau\le\sqrt{2/\alpha(\eps)}\eqqcolon\overline{\tau}_0(\eps),
$$
and  if, in addition, $L\le\eps\sqrt{2/\alpha(\eps)}\eqqcolon\overline{L}_0(\eps)$, then  $\alpha_1(\eps,L)=\sqrt{\alpha(\eps)/2}$ and we also get
$$
B(L|t|,\eps,L)\le \frac{\sqrt{\alpha(\eps)/2}}{2-\alpha(\eps)\tau^2} \eqqcolon B(\tau,\eps),
$$
whence, with the account of~\eqref{|f_n-e|<=(e-1)e}, we obtain the factors $(p(t,\eps,\gamma)+B(\tau,\eps)Lt^4)\cdot L^3e^{-t^2/2}$ in~\eqref{AbsChFDiffEstimAnalytEsseen} and~\eqref{AbsChFDiffEstimAnalytRozovsky}, which are monotonically decreasing with respect to $\gamma>0$ and $\eps>0$ and monotonically increasing with respect to $L>0.$
Finally, to bound $e^{\delta_n(t)}$ in~\eqref{|f_n-e|<=(e-1)e} we observe that  for all $L|t|\le\tau$
$$
p_{\scE}(t,\eps,\gamma)=
t^2\min_{0<z\le\eps}\left[\frac{zt^2}{24}+\max\left\{
\frac{|t|}{6\gamma}+\frac{zt^2}{24},\frac{\varkappa}{z} \right\}\right]\le
L^{-2}\tau^2\min_{0<z\le\eps}\left[\frac{zt^2}{24}+\max\left\{
\frac{\tau}{6\gamma L}+\frac{zt^2}{24},\frac{\varkappa}{z} \right\}\right].
$$
Majorizing now the maximum of two non-negative numbers by their sum we obtain
$$
L^3p_{\scE}(t,\eps,\gamma)\le
\frac{\tau^3}{6\gamma} + L\tau^2\min_{0<z\le\eps}\left\{ \frac{\varkappa}{z}+\frac{zt^2}{12} \right\} =\frac{\tau^3}{6\gamma} + L\tau^2\times
\begin{cases}
\frac{\eps t^2}{12}+\frac{\varkappa}{\eps}\le \sqrt{\frac\varkappa{12}}\,|t|+\frac{\varkappa}{\eps}, &\eps|t|\le\sqrt{12\varkappa},
\\[3mm]
\sqrt{\frac\varkappa3}\,|t|, &\eps|t|>\sqrt{12\varkappa},
\end{cases}\le
$$
$$
\le  \frac{\tau^3}{6\gamma} + L\tau^2\left(\sqrt{\frac\varkappa{3}}\,|t|+\frac{\varkappa}{\eps} \right)
\le  \frac{\tau^3}{6\gamma} + \sqrt{\frac\varkappa{3}}\, \tau^3+ \frac{\varkappa}{\eps}L\tau^2,
$$
with  the latest expression being monotonically decreasing with
respect to $\gamma>0$ and $\eps>0$ and monotonically increasing with
respect to $L>0.$ Thus the desired estimate takes the form
$$
e^{\delta_n(t)}\le \exp\left\{\tfrac{\tau^3}{6\gamma} +\sqrt{\tfrac{\varkappa}3}\tau^3 +\tfrac{\varkappa}\eps L\tau^2 +B(\tau,\eps)\tau^4\right\}\eqqcolon A_{\scE}(\tau,\eps,\gamma,L)\quad \text{for }\ L|t|\le\tau,
$$
whence, with the account of~\eqref{|f_n-e|<=(e-1)e} and the above
constructed estimate for $\delta_n(t)$, we
obtain~\eqref{AbsChFDiffEstimAnalytEsseen}. Similarly, for $L=\roz,$
using the explicit representation for $p_{\scR}$ in the right-hand
side of~\eqref{pR(t,eps,gamma)def}, we have
$$
L^3p_\scR(t,\eps,\gamma)
\le L^3\left(\frac{|t|^3}{6\gamma}\,\I(\gamma<\gopt)+  \frac{\varkappa t^2}{\eps} + \frac{\eps t^4}{12}\right)
\le \frac{\tau^3}{6\gamma}\,\I(\gamma<\gopt)+  \frac{\varkappa}{\eps}L\tau^2+ \frac{\eps L\tau^2}{12}t^2,
$$
$$
e^{\delta_n(t)}\le \exp\big\{\tfrac{\tau^3}{6\gamma}\I(\gamma<\gopt) +\tfrac{\varkappa}\eps L\tau^2 +\tfrac{\eps L\tau^2}{12}t^2  +B(\tau,\eps)\tau^4\big\}\eqqcolon A_{\scR}(\tau,\eps,\gamma,L)e^{\eps L\tau^2t^2/12}\quad \text{for }\ L|t|\le\tau,
$$
with the right-hand side monotonically increasing with respect to
$L>0$ and monotonically decreasing with respect to $\gamma>0.$
Hence, with the account of~\eqref{|f_n-e|<=(e-1)e}, we
obtain~\eqref{AbsChFDiffEstimAnalytRozovsky}.
\end{proof}

\section{Proofs of the main results}
\label{SecMainThmProof}

The proof of Theorem~\ref{ThEsseenRozovskiiIneq} will be given
simultaneously for both fractions $L=\es(\eps,\gamma)$ and
$L=\roz(\eps,\gamma)$ with $\eps>0$ and $\gamma>0$ being fixed. In
what follows by $C=C(\eps,\gamma)$ we mean $\essC(\eps,\gamma)$ or
$\rozC(\eps,\gamma)$, sometimes omitting the arguments
$\eps,\gamma$. Following the outline of Zolotarev's reasoning
employed in~\cite{Zolotarev1966} we construct an upper bound for
$C(\eps,\gamma)$ in the form $\sup_{L>0}C(\eps,\gamma,L)$, where
$C(\eps,\gamma,L)$ is such a function of~$L$ that the inequality
\begin{equation}\label{Delta_n<=L^3C(L)}
\Delta_n \le C(\eps,\gamma,L) \cdot L^3,
\end{equation}
holds for all $n\in\N$ and all distributions of independent centered
r.v.'s $X_1,\ldots,X_n$ with fixed value of the fraction of interest
$\es(\eps,\gamma)=L$ or $\roz(\eps,\gamma)=L$ for every $L>0$. Due
to the boundedness of the Kolmogorov distance $\Delta_n\le1$ for
arbitrary d.f.'s, one can exclude from consideration large values of
$L$. Moreover, the region of values of $L$ to be considered can be
restricted even more by use of the following sharpened upper bound
for $\Delta_n$ for standardized distributions.

\begin{lemma}[see~\cite{BhattacharyaRangaRao1982,Agnew1957}]\label{LemMaxDist}
For arbitrary r.v. $X$ with $0<\D X<\infty$ we have
$$
\sup_{x\in\R}\bigg|\Prob\bigg(\frac{X-\E X}{\sqrt{\D X}}<x\bigg)-\Phi(x)\bigg|\le  \sup_{x>0}\bigg|\frac1{1+x^2}-\Phi(-x)\bigg|=0.5409\ldots\ .
$$
\end{lemma}

Lemma~\ref{LemMaxDist} implies that, in order to prove the
inequality $\Delta_n \le C \cdot L^3$ with the constant $C \ge
C_{\min}\coloneqq2$, say, it suffices to consider $L<
(0.541/C_{\min})^{1/3}<0.65 =: L_1$. As this is so, we separate then
the interval $L\in(0,L_1]$ into the two regions:

(i) $L\in(0,L_0]$ with some $L_0=L_0(\eps)$ small enough, where mostly \textit{analytical} techniques is used;

(ii) $L\in(L_0,L_1]$, where \textit{numerical} computations with the help of a computer are needed;
\\
and construct an upper bound for $C$ as the maximum of the two corresponding bounds and the absolute constant $C_{\min}$:
$$
C(\eps,\gamma)\le\max\Big\{C_{\min},\max_{0<L\le L_1}C(\eps,\gamma,L)\Big\}=\max\Big\{C_{\min}, \max_{0<L\le L_0}C(\eps,\gamma,L), \max_{L_0< L\le L_1}C(\eps,\gamma,L)\Big\},
$$
where, in fact, the third term turns to be extremal. However, upper
bounds for the asymptotically exact constants
$\aexess(\eps,\gamma)$, $\aexroz(\eps,\gamma)$ are obtained as
limiting values of $C(\eps,\gamma,L)$ as $L\to0.$

To bound $\Delta_n$ above on each of the intervals $(0,L_0]$ and
$(L_0,L_1]$, we use the method of characteristic functions realized
by the Prawitz smoothing inequality.

\begin{lemma}[see \cite{Prawitz1972}]\label{LemPrawSmoothIneq}
For arbitrary d.f. $F$ with the ch.f. $f$ and for all $0 < T_0 \le T_1$ we have
$$
\sup\limits_{x \in \R}\abs{F(x) - \Phi(x)}\le
\frac{2}{T_1} \int_0^{T_0} \abs{K\left(\frac{t}{T_1}\right)} \cdot \abs{f(t) - e^{-t^2/2}}dt
+\frac{2}{T_1} \int_{T_0}^{T_1} \abs{K\left(\frac{t}{T_1}\right)} \cdot \abs{f(t)}dt+
$$
$$
+\frac{2}{T_1} \int_0^{T_0} \abs{K\left(\frac{t}{T_1}\right) - \frac{iT_1}{2\pi t}}e^{-t^2/2}dt
+\frac{1}{\pi} \int_{T_0}^\infty e^{-t^2/2} \frac{dt}{t},
$$
where
\[
K(t) = \frac{1}{2} \left(1 - \abs{t}\right) + \frac{i}{2}\left[ \left(1 - \abs{t}\right)\cot\pi t +  \frac{\sgn t}{\pi} \right], \quad t\in(-1,1)\setminus\{0\},
\]
and $K(\pm1)\coloneqq0$ defined by continuity. Moreover, the function $K(t)$ for all $t\in[-1,1]\setminus\{0\}$ satisfies
\begin{equation} \label{LemPrawSmoothIneq.Estims}
\abs{K(t)} \le \frac{1.0253}{2\pi \abs{t}},\quad
\abs{K(t) - \frac{i}{2\pi t}} \le
\frac{1}{2} \left(
1 - \abs{t} + \frac{\pi^2t^2}{18}
\right) \le \frac{1}{2}.
\end{equation}
\end{lemma}

By Lemma~\ref{LemPrawSmoothIneq} we have
\begin{equation}\label{proof.MainCompIneq}
\Delta_n \le I_1 + I_2 + I_3 + I_4,
\end{equation}
where
\[
I_1=\frac{2}{T_1} \int_0^{T_0} \abs{K\left(\frac{t}{T_1}\right)} \cdot \abs{\Sfn - e^{-t^2/2}}dt,
\quad
I_2=\frac{2}{T_1} \int_{T_0}^{T_1} \abs{K\left(\frac{t}{T_1}\right)} \cdot \abs{\Sfn}dt,
\]
\[
I_3=\frac{2}{T_1} \int_0^{T_0} \abs{K\left(\frac{t}{T_1}\right) - \frac{iT_1}{2\pi t}}e^{-t^2/2}dt,
\quad
I_4=\frac{1}{\pi} \int_{T_0}^\infty e^{-t^2/2} \frac{dt}{t},\quad 0<T_0\le T_1.
\]
In what follows we use the notation
\[
T_0(\tau_0,L)\coloneqq \frac{\tau_0}{L}, \quad T_1(\tau_1,L) \coloneqq \frac{\tau_1}{L^3},\quad\text{where}\quad \tau_0 \in (0,\overline\tau_0(\eps)),\quad \tau_1 \in (0,\overline{\tau}_1)
\]
are free parameters to be chosen below with $\overline\tau_1\coloneqq\pi/4$, $\overline\tau_0(\eps)\coloneqq\sqrt{2/\alpha(\eps)}$  defined in Theorems~\ref{ThAbsCh.F.Estim},~\ref{ThAbsChFDiffEstim}, respectively.

\subsection{The case (i) $L\in(0,L_0]$ and the proof of Theorem~\ref{ThAEXupperBounds}}
\label{SubsecCaseSmallL}

The purpose of the present subsection is to bound $C(\eps,\gamma,L)$
above for  $L\in(0,L_0]$ by an increasing function of $L$ with the
best possible (within the method used) limiting value as $L\to0$.

Let $L_0 > 0$ satisfy the conditions of
Theorem~\ref{ThAbsCh.F.Estim}\,(ii) and
Theorem~\ref{ThAbsChFDiffEstim}\,(ii):
$$
L_0<(\eps/4)^{1/ 3}\quad \text{and}\quad L_0\le\overline L_0(\eps)\coloneqq\eps\sqrt{{2}/{\alpha(\eps)}}=\eps\overline\tau_0(\eps),
$$
so that for $\eps\ge0.0557$ arbitrary $L_0\le0.03$ surely fits. For
$L=\roz$ we additionally assume that
$L_0\le6/(\eps\overline\tau_0^2(\eps))=3\alpha(\eps)/\eps$, which
provides the positiveness of the exponent in the right-hand side
of~\eqref{AbsChFDiffEstimAnalytRozovsky}. Since
$\alpha(\eps)\ge3\cdot2^{-2/3}$, this assumption is surely satisfied
if
$$
\eps L_0\le9\cdot 2^{-2/3}=5.6696\ldots,
$$
in particular, for $\eps\le188$ we may consider arbitrary $L_0\le0.03.$

Let us describe the process of estimation of each term
in~\eqref{proof.MainCompIneq} assuming that $T_0\coloneqq
T_0(\tau_0,L)$, $T_1\coloneqq  T_1(\tau_1,L)$ with $\tau_1$ also
satisfying the condition of Theorem~\ref{ThAbsCh.F.Estim}\,(ii):
$\tau_1\ge{\pi L_0^3}/{\eps}\eqqcolon\underline{\tau}_1(L_0,\eps)$.

Adding and substracting $(iT_1)/(2\pi t)$ from
$K\big(\frac{t}{T_1}\big)$ under the modulus sign in the integrand
in $I_1$ and applying then the inequality
$\abs{K\big(\frac{t}{T_1}\big) - \frac{iT_1}{2\pi t}}\le\frac12,$
$|t|\le T_1,$  from Lemma~\ref{LemPrawSmoothIneq}, we obtain
$$
I_1 \le \frac{1}{\pi} \int_0^{T_0} \frac{r_n(t)}{t}\,dt +  \frac{2}{T_1} \int_0^{T_0} \abs{K\left(\frac{t}{T_1}\right) - \frac{iT_1}{2\pi t}} r_n(t)\,dt \le I_{11} + I_{12},
$$
where
\[
r_n(t)\coloneqq\abs{\Sfn - e^{-t^2/2}},\quad I_{11}\coloneqq \frac{1}{\pi} \int_0^{T_0} \frac{r_n(t)}{t}\,dt,\quad
I_{12}\coloneqq \frac{1}{T_1} \int_0^{T_0} r_n(t)\,dt.
\]
Further we use inequalities~\eqref{AbsChFDiffEstimAnalytEsseen} and \eqref{AbsChFDiffEstimAnalytRozovsky} from Theorem~\ref{ThAbsChFDiffEstim} to estimate the integrands $r_n(t)$ in $I_{11}$ and $I_{12}$ and enlarge then the region of integration from $(0,T_0)$ to $(0,\infty)$. With the definitions of the upper and the lower incomplete gamma-functions yielding
\begin{equation}\label{proof.GammaEquations}
\int_{x}^{\infty} t^s e^{-kt^2} dt = \frac{k^{-\tfrac{s+1}{2}}}{2} \gmmahigh\left(\tfrac{s+1}{2},kx^2\right), \quad
\int_{0}^{x} t^se^{-kt^2} dt = \frac{k^{-\tfrac{s+1}{2}}}{2} \gmmalow\left(\tfrac{s+1}{2},kx^2\right), \quad  s,k,x > 0,
\end{equation}

\noindent $\bullet$ for the Esseen-type fraction $L=\es(\eps,\gamma)$  with arbitrary $\eps,\gamma>0$ we have
\begin{multline*}
\frac{I_{11}}{L^3} \le \frac{A_\scE}{\pi}  \bigg[ \int_{0}^{t_\gamma/\eps} \left(\tfrac{\varkappa t}{\eps} + \tfrac{\eps t^3}{24}\right) e^{-t^2/2}dt + \tfrac{\sqrt{6\varkappa \gamma^2 + 1}}{6\gamma}\int_{t_\gamma/\eps}^{\infty}t^2e^{-t^2/2} dt +BL\int_{0}^\infty t^3e^{-t^2/2}dt\bigg]=
\\
=\frac{A_\scE}{\pi}  \bigg[\frac{\varkappa}{\eps}\, \gmmalow\Big(1,\frac{t_\gamma^2}{2\eps^2}\Big)
+\frac{\eps}{12}\gmmalow\Big(2,\frac{t_\gamma^2}{2\eps^2}\Big) +\frac{\sqrt{2(6\varkappa\gamma^2+1)}}{6\gamma} \gmmahigh\Big(\frac{3}{2},\frac{t_\gamma^2}{2\eps^2}\Big)+2BL\bigg],
\end{multline*}
\begin{multline*}
\frac{I_{12}}{L^6}\le \frac{A_\scE}{\tau_1} \bigg[ \int_{0}^{t_\gamma/\eps} \left(\tfrac{\varkappa t^2}{\eps} + \tfrac{\eps t^4}{24}\right) e^{-{t^2}/2}dt
+\tfrac{\sqrt{6\varkappa \gamma^2 + 1}}{6\gamma}\int_{t_\gamma/\eps}^{\infty}t^3e^{-t^2/2} dt + BL\int_0^{\infty} t^4e^{-t^2/2} dt\bigg]=
\\
=\frac{A_\scE}{\tau_1}  \bigg[\frac{\sqrt2\varkappa}{\eps}\, \gmmalow\Big(\frac32,\frac{t_\gamma^2}{2\eps^2}\Big)
+\frac{\sqrt2\eps}{6}\gmmalow\Big(\frac52,\frac{t_\gamma^2}{2\eps^2}\Big) +\frac{\sqrt{6\varkappa\gamma^2+1}}{3\gamma} \gmmahigh\Big(2,\frac{t_\gamma^2}{2\eps^2}\Big) +\frac32\sqrt{2\pi}BL\bigg],
\end{multline*}
where $A_\scE= A_\scE(\tau_0,\eps,\gamma,L),$ $B = B(\tau_0,\eps);$

\noindent $\bullet$ in particular, for the Esseen-type fraction $L=\es(\infty,\gamma)$ with $\eps=\infty$  we obtain
$$
\frac{I_{11}}{L^3} \le
\frac{A_\scE}{\pi} \bigg[\frac{\sqrt{2\pi(6\varkappa\gamma^2+1)}}{12\gamma} +2BL\bigg],
\quad
\frac{I_{12}}{L^6} \le \frac{A_\scE}{\tau_1} \bigg[\frac{\sqrt{6\varkappa\gamma^2+1}}{3\gamma}+\frac32\sqrt{2\pi}BL\bigg];
$$

\noindent $\bullet$ for the Rozovskii-type fraction $L=\roz(\eps,\gamma)$  with arbitrary
$\eps,\gamma>0$ we have
\begin{multline*}
\frac{I_{11}}{L^3} \le \frac{A_{\scR}}{\pi}\bigg[  \int_{0}^{t_{1,\gamma}/\eps} \left(\frac{\varkappa t}{\eps} + \frac{\eps t^3}{24}\right) e^{-h{t^2}/2}dt
+\frac1{6\gamma} \int_{t_{1,\gamma}/\eps}^{t_{2,\gamma}/\eps} t^2e^{-h{t^2}/2}dt +\frac{\eps}{12} \int_{t_{2,\gamma}/\eps}^{\infty} t^3e^{-h{t^2}/2}dt+
\\
+BL\int_0^\infty t^3e^{-h{t^2}/2}dt \bigg]=
\frac{ A_{\scR}}{\pi} \bigg[
\frac{\varkappa}{\eps h}\,
\gmmalow\bigg(1,\frac{ht_{1,\gamma}^{2}}{2\eps^2}\bigg) +\frac{\eps}{12h^2}
\gmmalow\bigg(2,\frac{ht_{1,\gamma}^{2}}{2\eps^2}\bigg)
+
\\
+\frac{\sqrt2}{6\gamma h^{3/2}}\bigg(\frac{\sqrt\pi}{2}
-\gmmalow\bigg(\frac32,\frac{ht_{1,\gamma}^{2}}{2\eps^2}\bigg)
-\gmmahigh\bigg(\frac32,\frac{ht_{2,\gamma}^{2}}{2\eps^2}\bigg)\bigg)+
\frac{\eps}{6h^2}\, \gmmahigh\bigg(2,\frac{ht_{2,\gamma}^{2}}{2\eps^2}\bigg)
+\frac{2BL}{h^2}
\bigg],
\end{multline*}
\begin{multline*}
\frac{I_{12}}{L^6} \le \frac{A_{\scR}}{\tau_1}\bigg[  \int_{0}^{t_{1,\gamma}/\eps} \left(\frac{\varkappa t^2}{\eps} + \frac{\eps t^4}{24}\right) e^{-h{t^2}/2}dt
+\frac1{6\gamma} \int_{t_{1,\gamma}/\eps}^{t_{2,\gamma}/\eps} t^3e^{-h{t^2}/2}dt +\frac{\eps}{12} \int_{t_{2,\gamma}/\eps}^{\infty} t^4e^{-h{t^2}/2}dt+
\\
+BL\int_0^\infty t^4e^{-h{t^2}/2}dt \bigg]
=\frac{ A_{\scR}}{\tau_1} \bigg[
\frac{\sqrt2\varkappa}{\eps h^{3/2}}\,
\gmmalow\bigg(\frac32,\frac{ht_{1,\gamma}^{2}}{2\eps^2}\bigg) +\frac{\sqrt2\eps}{12h^{5/2}}
\gmmalow\bigg(\frac52,\frac{ht_{1,\gamma}^{2}}{2\eps^2}\bigg)
+
\\
+\frac{1}{3\gamma h^2}\bigg(1
-\gmmalow\bigg(2,\frac{ht_{1,\gamma}^{2}}{2\eps^2}\bigg)
-\gmmahigh\bigg(2,\frac{ht_{2,\gamma}^{2}}{2\eps^2}\bigg)\bigg)+
\frac{\sqrt2\eps}{6h^{5/2}}\, \gmmahigh\bigg(\frac52,\frac{ht_{2,\gamma}^{2}}{2\eps^2}\bigg)
+\frac{3\sqrt{2\pi}BL}{2h^{5/2}} \bigg],
\end{multline*}
where $A_\scR = A_\scR(\tau_0,\eps,\gamma,L),$ $B = B(\tau_0,\eps),$ $h=h(\tau_0,\eps,L)$
provided that $\eps L_0<9\cdot 2^{-2/3}=5.6696\ldots;$

\noindent $\bullet$ in particular, for the Rozovskii-type fraction $L=\roz(\eps,\gamma)$ with $\gamma\ge\gopt$, taking into account that $t_{1,\gamma}=t_{2,\gamma}=t_\infty=
2/\gopt$, we obtain
$$
\frac{I_{11}}{L^3} \le
\frac{ A_{\scR}}{\pi} \bigg[
\frac{\varkappa}{\eps h}\,
\gmmalow\bigg(1,\frac{2h}{(\eps\gopt)^2}\bigg) +\frac{\eps}{12h^2}\bigg(1+
\gmmahigh\bigg(2,\frac{2h}{(\eps\gopt)^2}\bigg)\bigg)
+\frac{2BL}{h^2}
\bigg],
$$
$$
\frac{I_{12}}{L^6} \le
\frac{ A_{\scR}}{\tau_1} \bigg[
\frac{\sqrt2\varkappa}{\eps h^{3/2}}\,
\gmmalow\bigg(\frac32,\frac{2h}{(\eps\gopt)^2}\bigg) +\frac{\sqrt2\eps}{12h^{5/2}}\bigg( \frac34\sqrt\pi+\gmmahigh\bigg(\frac52,\frac{2h}{(\eps\gopt)^2}\bigg)\bigg)
+\frac{3\sqrt{2\pi}BL}{2h^{5/2}}
\bigg].
$$
Note that the constructed upper bounds for $I_{11}/L^3$ and $I_{12}/L^6$ are monotonically increasing with respect to $L,$ monotonically decreasing with respect to $\gamma>0$ and, for $L=\es$, also with respect to $\eps>0$ as integrals of the functions possessing the stated properties. Moreover, they do not depend on $\gamma$ for $L=\roz(\eps,\gamma)$ with $\gamma\ge\gopt$, so that the further increase of $\gamma$ has no effect on the value of the resulting constant $C(\eps,\gamma,L)$, hence, the value $\gamma =\gopt$ is an optimal one.

To estimate $I_2$, we use the first inequality in~\eqref{LemPrawSmoothIneq.Estims}  from Lemma~\ref{LemPrawSmoothIneq} and bound~\eqref{AbsCh.F.EstimAnalytically} from Theorem~\ref{ThAbsCh.F.Estim} to get
$$
\frac{I_2}{L^3} \le \frac{1.0253}{\pi L^3} \int_{T_0}^{T_1} \frac{\abs{\Sfn}}{t} dt \le
\frac{1.0253}{\pi L^3T_0^3} \int_{T_0}^{\infty} t^2e^{-k(\tau_1)t^2} dt=
\frac{1.0253}{2 \pi \tau_0^3\left(k(\tau_1)\right)^{{3}/{2}} } \, \gmmahigh\left(\frac32, \frac{k(\tau_1)\tau_0^2}{L^2} \right),
$$
with $k(\tau_1)$ defined in the formulation of Theorem~\ref{ThAbsCh.F.Estim} (recall that $k(\tau_1)>0$ iff $0<\tau_1 < \overline{\tau}_1$). We also observe that the constructed upper bound for $I_2/L^3$ holds true for both fractions $L=\es(\eps,\gamma)$ and $L=\roz(\eps,\gamma)$, is independent of $\eps$ and $\gamma$, and is monotonically increasing with respect to $L$, moreover, $I_2=\mathcal O(L^\nu)$ as $L\to0$ for every $\nu>0$.

Using the inequality  $\abs{K\big(\frac{t}{T_1}\big) - \frac{iT_1}{2\pi t}}\le\frac12$ from Lemma~\ref{LemPrawSmoothIneq} once again to bound $I_3$ and the condition $1\le t^3/T_0^3$ defining the region of integration in $I_4$, we estimate the sum of the integrals $I_3$ and $I_4$  as
\[
\frac{I_3 + I_4}{L^3}\le \frac{1}{\tau_1} \int_{0}^{\infty} e^{-t^2/2} dt + \frac{1}{\pi L^3T_0^3} \int_{T_0}^\infty t^2 e^{-t^2/2} dt
=
\frac{1}{\tau_1} \sqrt{\frac{\pi}{2}} + \frac{\sqrt{2}}{\pi \tau_0^3} \, \gmmahigh\left(\frac{3}{2},\frac{\tau_0^2}{2L^2}\right)
\]
with the latest expression being monotonically increasing with respect to $L$ and independent of $\eps$ and $\gamma$.

Summing up the obtained bounds for the integrals $I_{11},I_{12},I_2,I_3,I_4,$ we obtain an upper bound for $\Delta_n$ in the form
\begin{equation}\label{Delta_n<=L^3C_0(L,tau0,tau1,eps,gamma)}
\Delta_n \le L^3\cdot C_0(\eps,\gamma,L,\tau_0,\tau_1)
\end{equation}
with the function $C_0(\eps,\gamma,L,\tau_0,\tau_1)$ being monotonically increasing with respect to~$L$, monotonically decreasing with respect to $\gamma>0$ and, for $L=\es(\eps,\gamma),$ also with respect to $\eps>0$ and satisfying
\begin{multline*}
C_0(\eps,\gamma,\tau_0,\tau_1)\coloneqq \lim_{L\to0}C_0(\eps,\gamma,L,\tau_0,\tau_1) =\frac{1}{\tau_1} \sqrt{\frac{\pi}{2}}+
\\
+
\begin{cases}
\frac{A_\scE(\tau_0,\eps,\gamma,0+)}{\pi}  \Big[\frac{\varkappa}{\eps} \gmmalow\Big(1,\frac{t_\gamma^2}{2\eps^2}\Big)
+\frac{\eps}{12}\gmmalow\Big(2,\frac{t_\gamma^2}{2\eps^2}\Big) +\frac{\sqrt{2(6\varkappa\gamma^2+1)}}{6\gamma} \gmmahigh\Big(\frac{3}{2},\frac{t_\gamma^2}{2\eps^2}\Big)\Big],
& L=\es(\eps,\gamma),
\\[2mm]
\frac{ A_\scR(\tau_0,\eps,\gamma,0+)}{\pi} \Big[
\frac{\varkappa}{\eps h}
\gmmalow\Big(1,\frac{ht_{1,\gamma}^{2}}{2\eps^2}\Big) +\frac{\eps}{12h^2}
\gmmalow\Big(2,\frac{ht_{1,\gamma}^{2}}{2\eps^2}\Big)+
\frac{\eps}{6h^2}\gmmahigh\Big(2,\frac{ht_{2,\gamma}^{2}}{2\eps^2}\Big)+
&
\\
\qquad\qquad\quad\ +\frac{\sqrt2}{6\gamma h^{3/2}}\Big(\frac{\sqrt\pi}{2}
-\gmmalow\Big(\frac32,\frac{ht_{1,\gamma}^{2}}{2\eps^2}\Big)
-\gmmahigh\Big(\frac32,\frac{ht_{2,\gamma}^{2}}{2\eps^2}\Big)\Big)
\Big],
& L=\roz(\eps,\gamma),
\end{cases}
\end{multline*}
where $h=h(\tau_0,\eps,L)\coloneqq1-\eps\tau_0^2L/6$.

Inequality~\eqref{Delta_n<=L^3C_0(L,tau0,tau1,eps,gamma)} yields the following upper bound for the function $C(\eps,\gamma,L)$ from~\eqref{Delta_n<=L^3C(L)}:
\begin{equation}\label{C0(L)=infC0(L,tau0,tau1,eps,gamma)}
C(\eps,\gamma,L) \le C_{0}(\eps,\gamma,L)\coloneqq \inf\left\{
C_0(\eps,\gamma,L,\tau_0,\tau_1)\colon \tau_0 \in\big(0,\overline{\tau}_0(\eps)\big),\ \tau_1 \in \big(\underline{\tau}_1(L_0,\eps), \overline{\tau}_1\big),\ \tau_1 \ge L^2\tau_0
\right\}
\end{equation}
for every $\eps>0,\,\gamma>0,$ and $L\le L_0$ with
\begin{equation}\label{L_0<=restrictions}
L_0\le\overline L_0(\eps)\wedge(\eps/4)^{1/3},\ \text{ also $L_0\le 9\cdot 2^{-2/3}\eps^{-1}$ for $L=\roz$,}
\end{equation}
and $C_{0}(\eps,\gamma,L)$  being a monotonically increasing function of $L$ as the greatest lower bound to the increasing function $C_0(\eps,\gamma,L,\tau_0,\tau_1)$, where  the greatest lower bound  is taken over a  decreasing system of  sets. Hence,
\begin{equation*}
\max_{0<L\le L_0}C(\eps,\gamma,L) \le\max_{0<L\le L_0}C_0(\eps,\gamma,L)= C_0(\eps,\gamma,L_0).
\end{equation*}
Moreover, $C_0(\eps,\gamma,L)$ is monotonically decreasing with respect to $\gamma>0$ and, for $L=\es(\eps,\gamma),$ also with respect to $\eps>0$.

Inequality~\eqref{Delta_n<=L^3C_0(L,tau0,tau1,eps,gamma)} yields upper bounds for the asymptotically exact constants $C^*\in\left\{\aexess,\aexroz\right\}$ defined in~\eqref{aexCessDef}, \eqref{aexCrozDef}. Namely, with the observation that conditions~\eqref{L_0<=restrictions} are trivially satisfied for every $\eps>0$ if $L\to0$, we have
$$
C^*(\eps,\gamma)\le  \lim_{L\to0} C(\eps,\gamma,L)\le C_0(\eps,\gamma,0+)\le 
\inf\left\{C_0(\eps,\gamma,\tau_0,\tau_1)\colon \tau_0 \in\big(0,\overline{\tau}_0(\eps)\big),\ \tau_1 \in \big(0, \overline{\tau}_1\big) \right\} =
$$
$$
= C_0\big(\eps,\gamma,0+,\overline \tau_1-\big)\quad \text{for every }\ \eps,\gamma>0.
$$
Hence, recalling that $\overline \tau_1=\pi/4,$ $A_\scE(0+,\eps,\gamma,0+)=A_\scR(0+,\eps,\gamma,0+)=A(0+,\eps)=1$, $h(\tau_0,\eps,0+)=1$, for every $\eps>0$ and $\gamma>0$ we obtain:
\\
\noindent$\bullet$ in the Esseen-type inequality
$$
\textstyle \aexess(\eps,\gamma)\le C_0\big(\eps,\gamma,0+,\overline \tau_1-\big)=\frac{4}{\sqrt{2\pi}}+ \frac{1}{\pi}  \Big[\frac{\varkappa}{\eps}\gmmalow\big(1,\frac{t_\gamma^2}{2\eps^2}\big) +\frac{\eps}{12}\gmmalow\big(2,\frac{t_\gamma^2}{2\eps^2}\big) +\frac{\sqrt{2(6\varkappa\gamma^2+1)}} {6\gamma}\, \gmmahigh\big(\frac{3}{2},\frac{t_\gamma^2}{2\eps^2}\big) \Big]= \widehat\aexess(\eps,\gamma),
$$
\\
\noindent$\bullet$ in the Rozovskii-type inequality
\begin{multline*}
\aexess(\eps,\gamma)\le C_0\big(\eps,\gamma,0+,\overline \tau_1-\big)=\frac{4}{\sqrt{2\pi}}+ \frac{1}{\pi}  \Big[
\frac{\varkappa}{\eps}
\gmmalow\Big(1,\frac{t_{1,\gamma}^{2}}{2\eps^2}\Big) +\frac{\eps}{12}
\gmmalow\Big(2,\frac{t_{1,\gamma}^{2}}{2\eps^2}\Big)+
\frac{\eps}{6}\gmmahigh\Big(2,\frac{t_{2,\gamma}^{2}}{2\eps^2}\Big)+
\\
+\frac{\sqrt2}{6\gamma}\Big(\frac{\sqrt\pi}{2}
-\gmmalow\Big(\frac32,\frac{t_{1,\gamma}^{2}}{2\eps^2}\Big)
-\gmmahigh\Big(\frac32,\frac{t_{2,\gamma}^{2}}{2\eps^2}\Big)\Big)
\Big]= \widehat\aexroz(\eps,\gamma),
\end{multline*}
with $\widehat\aexess(\eps,\gamma)$, $\widehat\aexroz(\eps,\gamma)$ defined in the formulation of Theorem~\ref{ThAEXupperBounds} (see~\eqref{aexess(eps,gamma)<=}, \eqref{aexroz(eps,all_gamma)<=}). Moreover, the function $\widehat\aexess(\eps,\gamma)$ decreases with respect to $\eps>0$ and $\gamma>0$ as an integral and then a limit of a function with the similar properties and is unbounded as $\eps\to0+$ or $\gamma\to0+$. Similarly, the function $\widehat\aexroz(\eps,\gamma)$ decreases with respect to $\gamma>0$ being constant for $\gamma\ge\gopt$ for every fixed $\eps>0$ and infinitely grows as $\eps\to0+$, $\eps\to\infty,$ or $\gamma\to0+$.

\subsection{The case (ii) $L\in[L_0,L_1]$}
\label{SubsecCaseLargeL}

Let $0 < T_0 \le T_0(\overline{\tau}_0(\eps),L),$ $T_0\le T_1 \le T_1(\overline{\tau}_1,L)$, and $0<L_0\le L\le L_1<\infty$.

Though the function $K(t)$ has a singularity of order  $\mathcal O\big(|t|^{-1}\big)$ as $t\to0$, the integrands in $I_1$ and $I_3$ have no singularities due to the presence of the factor $r_n(t)=\mathcal O(t^2)$, $t\to0,$ in the integrand in $I_1$ and to the  boundedness of the function $\big|K\big(\frac{t}{T_1}\big) - \frac{iT_1}{2\pi t}\big|$ by~\eqref{LemPrawSmoothIneq.Estims} in $I_3$.

Using estimates~\eqref{AbsChFDiffEstimComp} and~\eqref{AbsCh.F.EstimComp} from Theorems~\ref{ThAbsChFDiffEstim} and~\ref{ThAbsCh.F.Estim} to bound integrands in $I_1$ and $I_2$,  we obtain, by~\eqref{proof.MainCompIneq}, an upper bound for $\Delta_n$ in the form
$$
\Delta_n\le D(\eps,\gamma,L,T_0,T_1),
$$
which is uniform in the class of all distributions of random summands with fixed value of the fraction under consideration  $L\in[L_0,L_1]$, where $T_0\le T_1$  are free parameters. Moreover, the function $D(\eps,\gamma,L,T_0,T_1)$ here is monotonically increasing with respect to $L$, monotonically decreasing with respect to $\gamma>0$ and, for $L=\es,$ also with respect to $\eps>0$.  Hence, we may construct an upper bound for  $C(\eps,\gamma,L)$  on the interval $L_0 \le L \le L_1$ in the form
\begin{equation}\label{C1(L)def}
C(\eps,\gamma,L) \le C_1(\eps,\gamma,L) := \inf \left \{
\frac{D(\eps,\gamma,L,T_0,T_1)}{L^3} \colon 0 < T_0 < \frac{\overline{\tau}_0(\eps)}L,\ T_0 < T_1 < \frac{\overline{\tau}_1}{L^3} \right \},
\end{equation}
where $\overline{\tau}_1\coloneqq\pi/4$ and $\overline\tau_0(\eps)\coloneqq  \sqrt{2/{\alpha(\eps)}}$ are defined in Theorems~\ref{ThAbsCh.F.Estim} and~\ref{ThAbsChFDiffEstim}, respectively. The maximum value $\max_{L_0\le L\le L_1}C_1(\eps,\gamma,L)$  can be estimated similarly to~\cite{Zolotarev1966} by computation of $C_1(\eps,\gamma,L)$ in a finite number of points using the inequality
\[
\max\limits_{L'\le L \le L''} C_1(\eps,\gamma,L) \le C_1(\eps,\gamma,L'')\cdot\big({L''}/{L'}\big)^3,
\]
which is valid due to the monotone growth of the function $C_1(\eps,\gamma,L)\cdot L^3$ as the greatest lower bound to the monotonically increasing function $D(\eps,\gamma,L,T_0,T_1)$
over a decreasing system of sets. Furthermore, since the function $D(\eps,\gamma,L,T_0,T_1)$ is monotonically decreasing with respect to $\gamma$ and, for $L=\es$, with respect to $\eps$, so is $\max\limits_{L_0\le L\le L_1}C_1(\eps,\gamma,L)$.

\subsection{Numerical results}

Summarizing what was said above, as the constants $C(\eps,\gamma)\in\{\essC(\eps,\gamma),\rozC(\eps,\gamma)\}$ in inequalities~\eqref{MainEsseenIneq} and~\eqref{MainRozovskyIneq} we can take
$$
C(\eps,\gamma)\coloneqq \max\big\{C_{\min}, \sup_{0 < L \le L_1} C(\eps,\gamma,L)\big\}
\ \text{ with }\ C(\eps,\gamma,L)\coloneqq
\begin{cases}
C_0(\eps,\gamma,L),& 0<L\le L_0(\eps),
\\
C_1(\eps,\gamma,L),&L_0(\eps)<L\le L_1,
\end{cases}
$$
$C_{\min}\coloneqq 2$, $L_1\coloneqq0.65>(0.541/C_{\min})^{1/3}$,
$L_0(\eps)\coloneqq (\eps/4)^{1/ 3}\wedge \eps\overline\tau_0(\eps)$ if $L=\es$ and $L_0(\eps)\coloneqq (\eps/4)^{1/ 3}\wedge \eps\overline\tau_0(\eps)\wedge6/(\eps\overline\tau_0^2(\eps))$ if $L=\roz,$ for every $\eps,\gamma>0$. Since $C_0(\eps,\gamma,L)$ and $C_1(\eps,\gamma,L)$ are both monotonically decreasing with respect to $\gamma>0$ and, for $L=\es$, also with respect to $\eps$, so is $C(\eps,\gamma)$.

The concrete numerical values of $C_0(\eps,\gamma,L)$ and
$C_1(\eps,\gamma,L)$ are processed with the help of a computer. Our
computations were carried out in Python~3.6 using the library
Scipy~1.0.0. The values of $\max_{0<L\le
L_0}C_0(\eps,\gamma,L)=C_0(\eps,\gamma,L_0)$
for some $\eps$ and $\gamma$ with $L_0=0.001$ and $L_0=0.03$ are given
for the Esseen-type fraction $L=\es(\eps,\gamma)$ in table~\ref{Tab:EsseenC0(L)}
in the fourth and seventh columns, respectively, accompanied by the optimal values
of the parameters $\tau_0$ and $\tau_1$ in~\eqref{C0(L)=infC0(L,tau0,tau1,eps,gamma)}.
The values of $\max_{0.03\le L\le L_1}C_1(\eps,\gamma,L)=C_1(\eps,\gamma,L_*)$
for some $\eps$ and $\gamma$ are given in table~\ref{Tab:EsseenC1(L)} for the
Esseen fraction $L=\es(\eps,\gamma)$  accompanied by the optimal values of the
parameters~$T_0$ and~$T_1$ in~\eqref{C1(L)def} specified in the form
$\tau_0\coloneqq T_0L$ and $\tau_1\coloneqq T_1L^3$, 
in the extremal point $L=L_*$ which is also given in the fourth
column of the same table. Columns 7\,--\,10 contain the normalized
contributions of the integrals $I_k/L_*^3$, $k=1,2,3,4,$ into the
extreme value~$C_1(\eps,\gamma,L_*)$, so that
$C_1(\eps,\gamma,L_*)=(I_1+I_2+I_3+I_4)/L_*^3$.
Tables~\ref{Tab:RozovskyC0(L)} and~\ref{Tab:RozovskyC1(L)} contain
similar results for the Rozovskii-type fraction
$L=\roz(\eps,\gamma)$.

\begin{figure}[h!]
\includegraphics[width=0.5\textwidth]{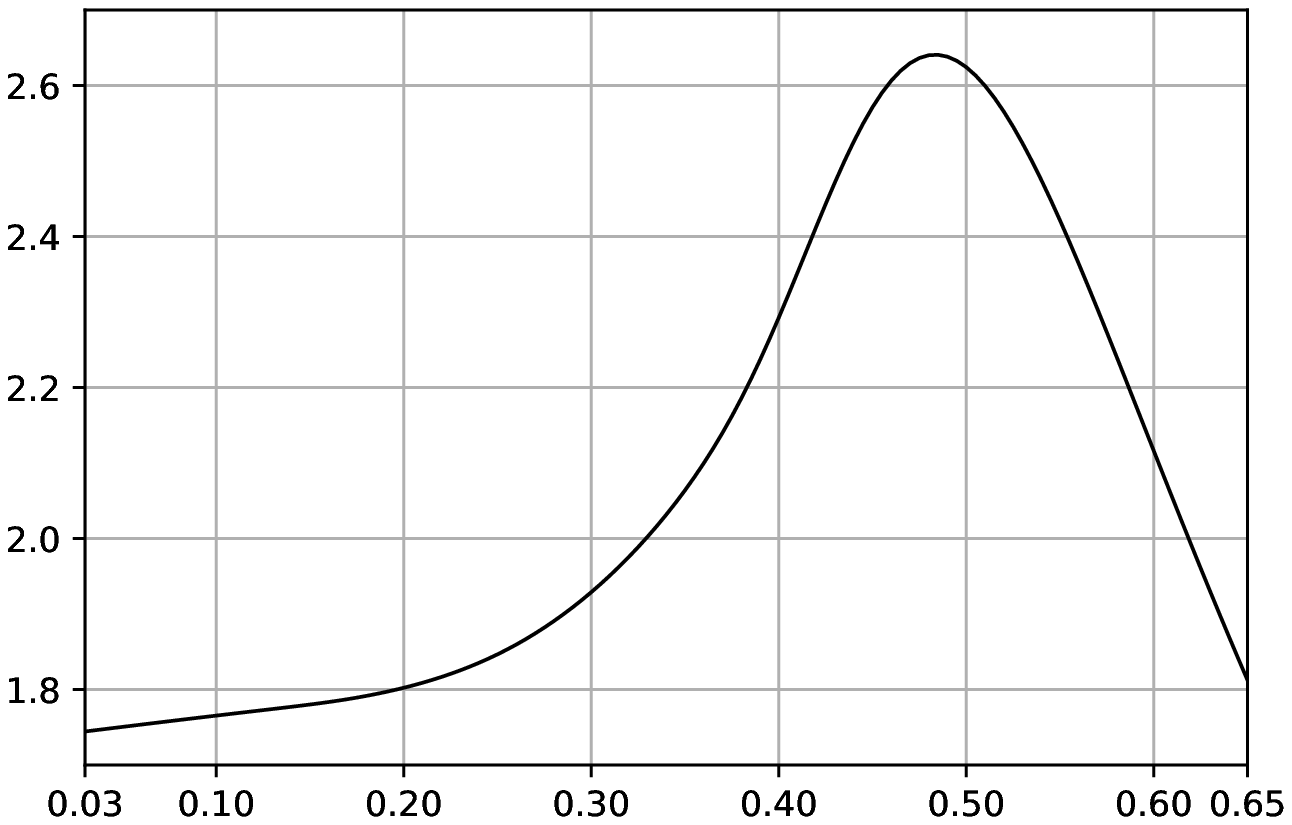}
\includegraphics[width=0.5\textwidth]{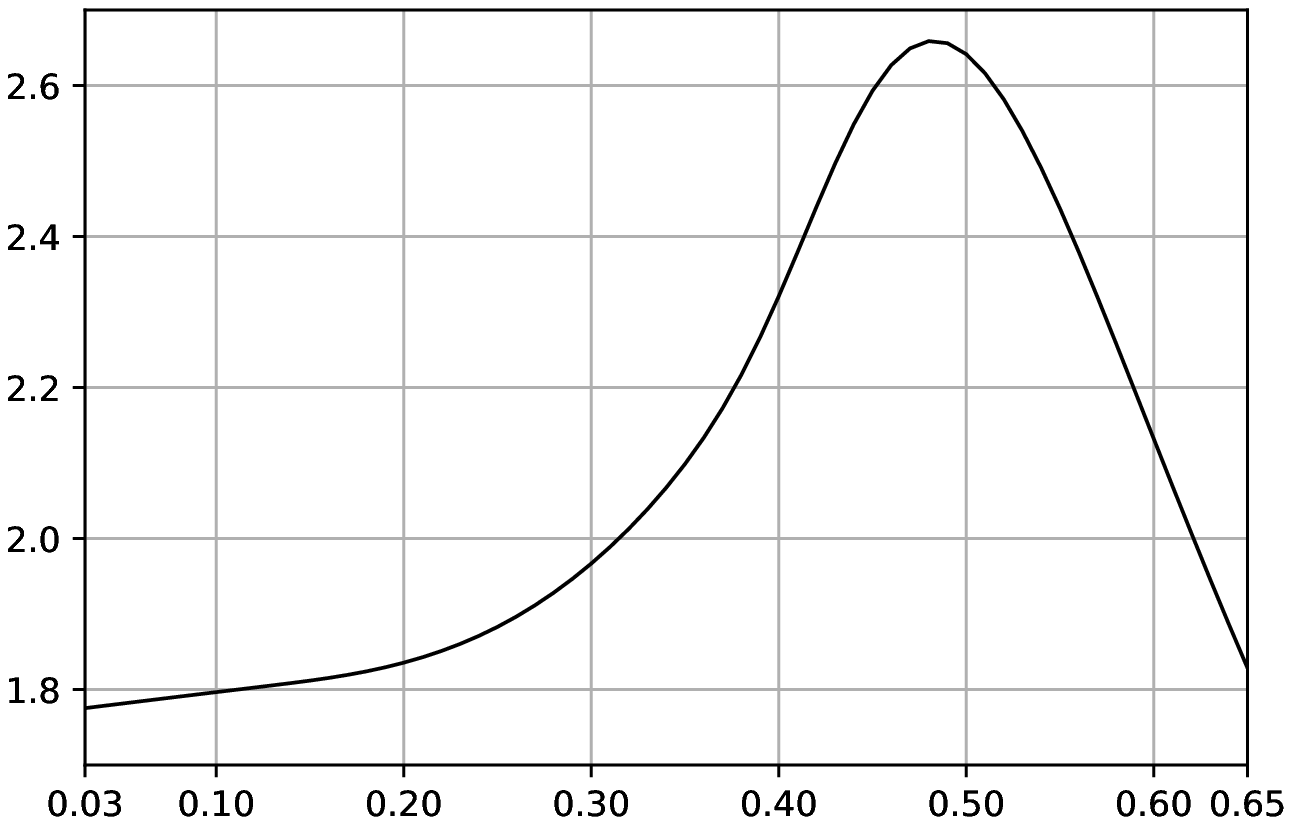}
\caption{Demonstration to the evaluation of the absolute constants $\essC(\infty,\infty)$ and $\rozC(2.12,\gopt)$: Plots of the functions $C_1(\infty,\infty,L)$, $L\in[L_0,L_1]$, defined in~\eqref{C1(L)def} with $L\coloneqq\es(\infty,\infty)$ (left) and $L\coloneqq\roz(2.12,\gopt)$ (right).}
\label{Fig:C1(L)plot}
\end{figure}

For example, for the minimum value $\essC(\infty,\infty)$ of the
constant $\essC(\eps,\gamma)$ in the Esseen-type
inequality~\eqref{MainEsseenIneq}  we have
$$
\sup_{0<L\le0.03}C(\infty,\infty,L)\le C_0(\infty,\infty,0.03)\le2.28,
\quad  \sup_{0.03<L\le L_1}C(\infty,\infty,L)\le C_1(\infty,\infty,0.4833\ldots)\le2.65
$$
(the plot of the function $C_1(\infty,\infty,L)$ for
$L\coloneqq\es(\infty,\infty)\in[0.03,L_1]$ is given on
Fig.\,\ref{Fig:C1(L)plot} (left)), hence,
$$
\inf_{\eps>0,\gamma>0}\essC(\eps,\gamma)=
\essC(\infty,\infty)\le\max\{C_{\min},\,2.28,\,2.65\}=2.65.
$$
However, the same upper bound $\essC(\eps,\gamma)\le2.65$ can be
reached, due to the rounding gap, already for finite values of
$\eps$ and $\gamma=\gamma(\eps)$ plotted on
Fig.\,\ref{Fig:aex+absEssPlotGamma(Eps)} (right) and also at the
points with one infinite component, for example,
$$
(\eps,\gamma)\in\big\{(2.56,\infty),\ (2.74,3),\ (4,1.62),\ (\infty,1.43)\big\},
$$
for which, due to the monotonicity and according to tables~\ref{Tab:EsseenC0(L)}, \ref{Tab:EsseenC1(L)}, we have
$$
\max\big\{ C(2.56,\infty,L),\,C(2.74,3,L),\,C(4,1.62,L),\,C(\infty,1.43,L) \big\}=
$$
$$
=\begin{cases}
 C_0(0.6,0.3,0.03)\le2.64,&L\le0.03\quad\text{(see table~\ref{Tab:EsseenC0(L)})},
\\
\max\big\{C_1(2.56,\infty,0.4833\ldots),\ldots\big\}\le2.65,&0.03<L\le L_1\quad\text{(see table~\ref{Tab:EsseenC1(L)})}.
\end{cases}
$$
Similar level curve
$\{(\eps,\gamma)\colon \widehat\aexess(\eps,\gamma)=1.72\}$ of the upper
bound $\widehat\aexess(\eps,\gamma)\coloneqq C_0(\eps,\gamma,0+)$ for the
asymptotically exact constant $\aexess(\eps,\gamma)$ is plotted on Fig.\,\ref{Fig:aex+absEssPlotGamma(Eps)} (left). 

\begin{figure}[h!]
\includegraphics[width=0.5\textwidth]{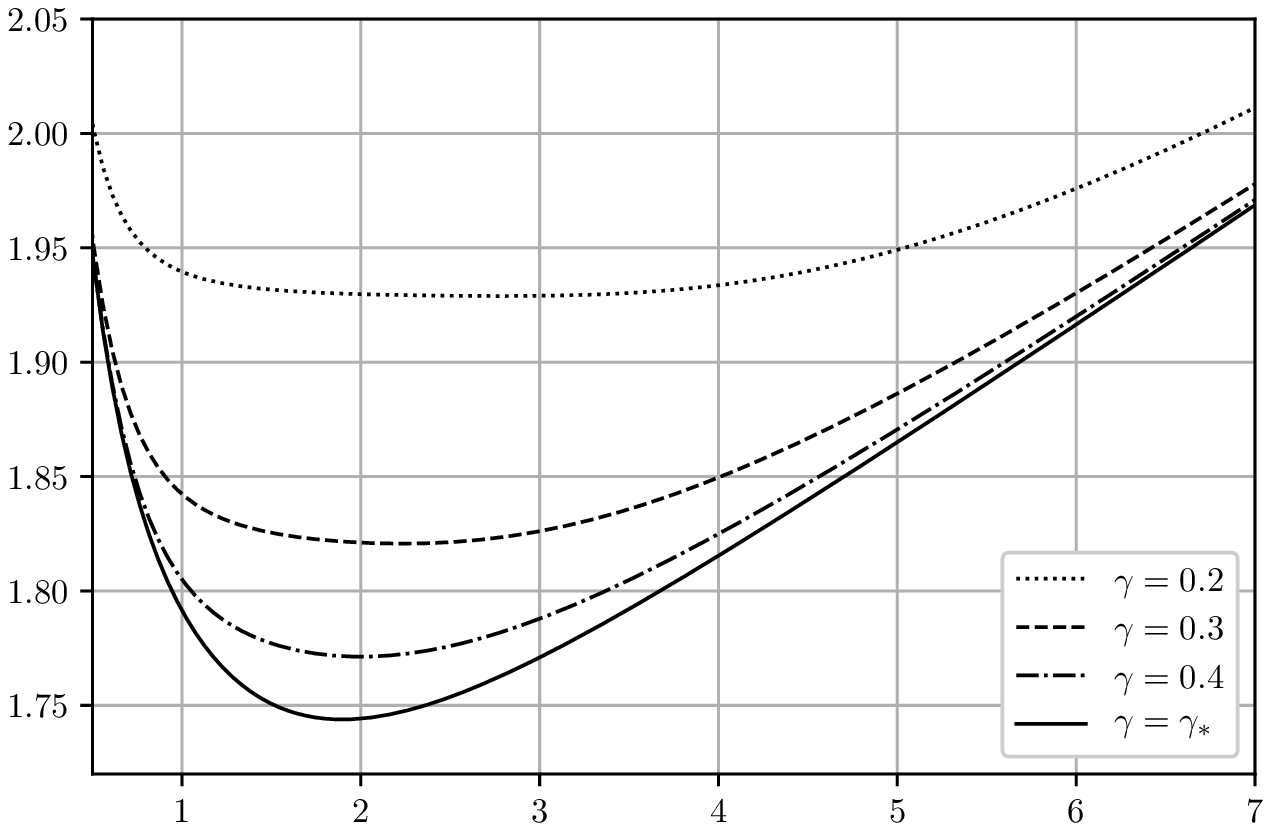}
\includegraphics[width=0.5\textwidth]{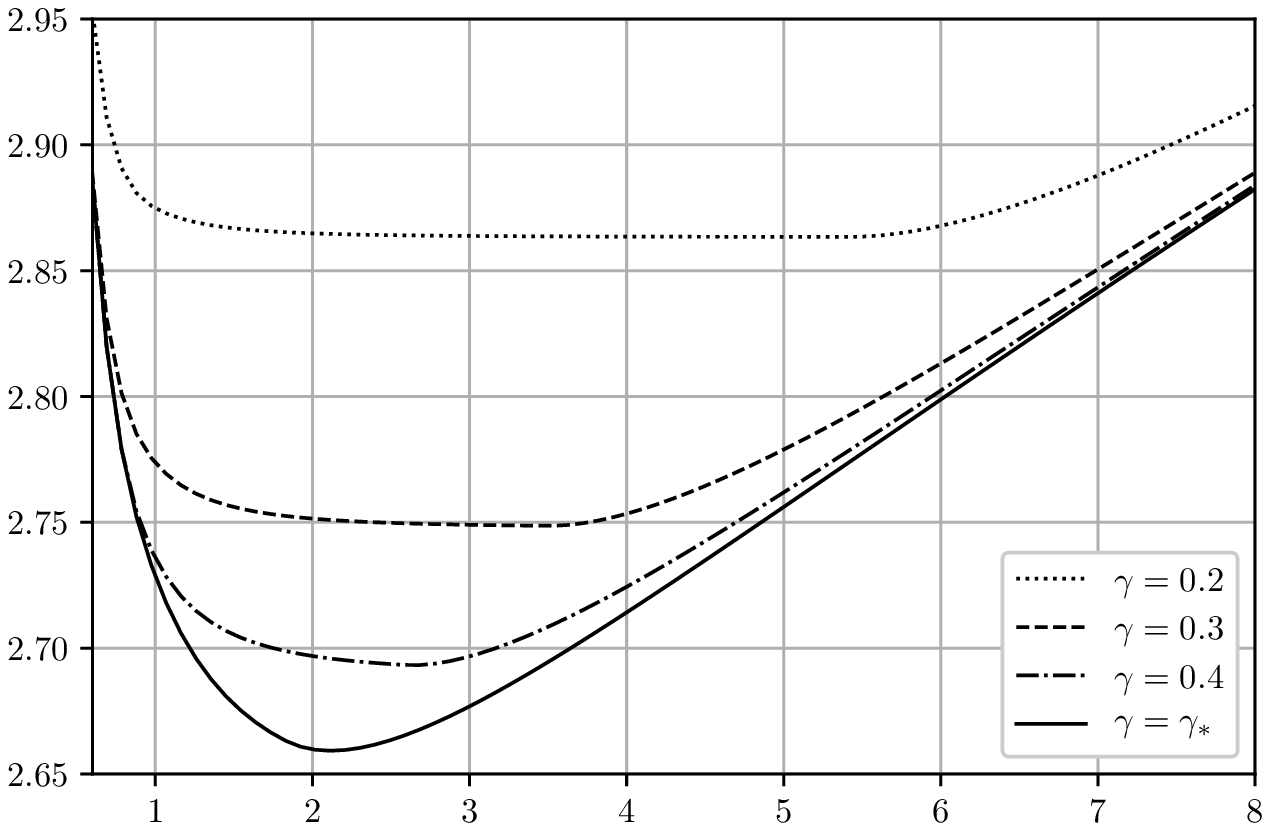}
\caption{Left:
Graphs of the function $\widehat\aexroz\big(\eps,\gamma\big)$ (see~\eqref{aexroz(eps,gamma)<=}) which bounds from above the asymptotically exact constant $\aexroz(\eps,\gamma)$ (see~\eqref{aexCrozDef}), with respect to $\eps$, for $\gamma\ge\gopt=0.5599\ldots$ (solid line), $\gamma=0.4$ (dashdot line), $\gamma=0.3$ (dashed line), and $\gamma=0.2$ (dotted line).
\\
Right: Graphs of the upper bounds $\max_{L_0\le L\le L_1}C_1(\eps,\gamma,L)$ (see~\eqref{C1(L)def}) for the ``constants'' $\rozC(\eps,\gamma)$ in the Rozovskii-type inequality, with respect to $\eps$, for $\gamma\ge\gopt=0.5599\ldots$ (solid line), $\gamma=0.4$ (dashdot line), $\gamma=0.3$ (dashed line), and $\gamma=0.2$ (dotted line).}
\label{Fig:aex+absRozPlotEps}
\end{figure}

Since the constructed upper bound for $\rozC(\eps,\gamma)$ is
monotonically decreasing with respect to $\gamma$ and is constant
for $\gamma\ge\gopt$, its global minimum with respect to $\gamma$ is
attained at $\gamma=\gopt=0.5599\ldots\ .$ Plot of the function
$\max_{0.03\le L\le L_1}C_1(\eps,\gopt,L)$ with respect to $\eps$ is
given on Fig.~\ref{Fig:aex+absRozPlotEps} (right, solid line), and
numerical computations show that its minimum is attained around the
point $\eps=2.12$, for which, according to
tables~\ref{Tab:RozovskyC1(L)} and~\ref{Tab:RozovskyC0(L)}, we have
$$
C_0(2.12,\gopt,0.03)\le2.29,\quad
\max_{0.03\le L\le L_1}C_1(2.12,\gopt,L)=
C_1(2.12,\gopt,0.4827\ldots)\le2.66
$$
(the plot of the function $C_1(2.12,\gopt,L)$ for
$L\coloneqq\roz(2.12,\gopt)\in[0.03,L_1]$ is given on
Fig.\,\ref{Fig:C1(L)plot} (right)).
Hence,
$$
\inf_{\eps>0,\,\gamma>0}\rozC(\eps,\gamma)\le \rozC(2.12,\gopt) \le\max\{C_{\min},\,2.29,\,2.66\}=2.66.
$$
We also provide the least value of $\eps$ (up to the second decimal
digit), for which, within the numerical method used, the same upper
bound $2.66$ holds,
$$
\inf\big\{\eps>0\colon \rozC(\eps,\gopt)\le2.66\big\}\approx1.99.
$$
The interest to exactly the least value of $\eps$ is stipulated by
that the second term ($\sup_{0<z\le\eps}\{\ldots\}$) in the
definition of $\roz(\eps,\gamma)$ is monotonically increasing with
respect $\eps>0.$

Similarly, the constructed upper bound
$\widehat\aexroz(\eps,\gamma)$ (defined
in~\eqref{aexroz(eps,gamma)<=}) for the asymptotically exact
constant  $\aexroz(\eps,\gamma)$ (defined in~\eqref{aexCrozDef})
attains its minimum value at the point $\gamma=\gopt$,
$\eps\approx1.89$ with
$$
\inf_{\eps>0,\gamma>0}\widehat\aexroz(\eps,\gamma)\le \widehat\aexroz(1.89,\gopt)\le1.75;
$$
the plot of $\widehat\aexroz(\eps,\gopt)$ is given on
Fig.~\ref{Fig:aex+absRozPlotEps} (left, solid line), while the least
value of $\eps$, for which the inequality
$\widehat\aexroz(\eps,\gopt)\le1.75$ still holds, is around $1.52.$
Furthermore, as Fig.~\ref{Fig:aex+absRozPlotEps} demonstrates, with
the decrease of $\gamma$, the graphs of
$\widehat\aexroz(\,\cdot\,,\gamma)$ and $\rozC(\,\cdot\,,\gamma)$
become more and more flat around the points of minimum, so that the
intervals, where the values of the functions
$\widehat\aexroz(\,\cdot\,,\gamma)$ and $\rozC(\,\cdot\,,\gamma)$
differ slightly from their minimal values,
$\min_{\eps>0}\widehat\aexroz(\eps,\gamma)$ and
$\min_{\eps>0}\rozC(\eps,\gamma)$, enlarges. For example, for
$\gamma=0.4$ (see Fig.~\ref{Fig:aex+absRozPlotEps} (left, dashdot
line) and table~\ref{Tab:RozovskyC0(L)}) the point of minimum of
$\widehat\aexroz(\,\cdot\,,\gamma)$ is located around
$\eps\approx1.99$ with $\widehat\aexroz(1.99,0.4)\le1.78$, while the
inequality $\widehat\aexroz(\eps,0.4)\le1.78$ still holds for
$\eps=1.41$; similarly, for $\gamma=0.2$
$$
\inf_{\eps>0}\widehat\aexroz(\eps,0.2)\le\widehat\aexroz(2.77,0.2)\le1.93,\quad
\inf\{\eps>0\colon \widehat\aexroz(\eps,0.2)\le1.93\}\approx1.89,
$$
see Fig.~\ref{Fig:aex+absRozPlotEps} (left, dotted line) and
table~\ref{Tab:RozovskyC0(L)}. For the ``absolute'' constant
$\rozC(\eps,\gamma)$ with $\gamma\in\{0.2,0.4\}$ we have
$$
\inf_{\eps>0,\,\gamma\ge0.4}\rozC(\eps,\gamma)\le \rozC(2.63,0.4) \le\max\{C_{\min},\,C_0(2.63,0.4,0.03),C_1(2.63,0.4,0.4822\ldots\} \le2.70,
$$
$$
\inf\{\eps>0\colon \rozC(\eps,0.4)\le2.70\}\approx1.76,\quad \rozC(1.76,0.4)\le\max\{2,\,2.37,\,2.70\}=2.70;
$$
$$
\inf_{\eps>0,\,\gamma\ge0.2}\rozC(\eps,\gamma)\le \rozC(5.39,0.2) \le\max\{2,\,2.68,\,2.87\}=2.87,\qquad\qquad\qquad\qquad\qquad\qquad\qquad\quad
$$
$$
\inf\{\eps>0\colon \rozC(\eps,0.2)\le2.87\}\approx1.21,\quad \rozC(1.21,0.2)\le\max\{2,\,2.64,\,2.87\}=2.87;
$$
see Fig.~\ref{Fig:aex+absRozPlotEps} (right, dashdot and dotted lines) and tables~\ref{Tab:RozovskyC0(L)} and~\ref{Tab:RozovskyC1(L)}.

Another particular case is concerned with the historical values
$\gamma=1$ and $\eps\in\{1,\infty\}$, for which, according to
tables~\ref{Tab:EsseenC0(L)}, \ref{Tab:RozovskyC0(L)},
\ref{Tab:EsseenC1(L)}, and~\ref{Tab:RozovskyC1(L)},  we have
$$
\rozC(1,1)\le\rozC(1,\gopt)\le 2\vee C_0(1,\gopt,0.03)\vee  C_1(1,\gopt,0.4834\ldots)
\le2\vee 2.35\vee2.73=2.73,
$$
$$
\essC(1,1)\le\essC(1,0.72)\le2.73\ge \essC(1,\infty),
$$
$$
\essC(\infty,1)\le2.66,\ \text{ which also follows from }\
\essC(\infty,0.97)\le2.66\ \text{ or }\ \essC(4.35,1)\le 2.66.
$$

It is interesting to note that, as it follows from
Tables~\ref{Tab:EsseenC1(L)} and~\ref{Tab:RozovskyC1(L)}, the
largest contribution into extreme values of
$C(\eps,\gamma)=C_1(\eps,\gamma,L_*)$ in both
inequalities~\eqref{MainEsseenIneq}, \eqref{MainRozovskyIneq} for
all the presented values of $\eps$ and $\gamma$ is provided by the
integral $I_3$ which depends on the constructed estimates for
characteristic functions through the maximal length of the interval,
where the absolute value of a characteristic function can be
estimated by a majorant strictly less than $1$. Hence, to get
further improvements of the constructed upper bounds for
$C(\eps,\gamma)$, one should improve, in first turn, upper bounds
for absolute values of characteristic functions presented in
Theorem~\ref{ThAbsCh.F.Estim}.

\section{The comparison of Osipov's, Lyapunov's, and modified Esseen's and Rozovskii's fractions. Lower bound for $\gamma\to0$.}
\label{SecComparisons}

In the present section we compare the fractions $\es^3(\eps,1)$,
$\roz^3(\eps,1)$ with~$\lyap$ and $\osiM(1)+\lind(1)$, and
demonstrate that our new inequalities~\eqref{MainEsseenIneq},
\eqref{MainRozovskyIneq} with $\essC(1,1)=2.73=\rozC(1,1)$ may be
sharper than Osipov's inequality~\eqref{OsipovIneqEps=1} with the
best known constant $C=1.87$~\cite{KorolevDorofeyeva2017}.
In what follows we emphasize the dependence of the above fractions
on the distributions of random summands $X_1,\ldots,X_n$ with the
d.f.'s $F_1,\ldots,F_n$ by writing
$\es^3(\eps,\gamma,F_1,\ldots,F_n)$,
$\roz^3(\eps,\gamma,F_1,\ldots,F_n)$, and $\lyap(F_1,\ldots,F_n)$
using  the three-argument notation $\es^3(\eps,\gamma,F)$,
$\roz^3(\eps,\gamma,F)$, $\lyap(F)$, $n\in\N$, in the i.i.d. case,
that is, for $F_1=\ldots=F_n=F$. Let~$\mathcal F$ denote the set of
all d.f.'s on $\R$ with zero mean and finite second-order moment and
$F_p\in\mathcal F$ be the d.f. of the two-point distribution
prescribing the masses $p\in\big[\frac12,1\big)$ and $q=1-p$ to the
points $\sqrt{q/p}$ and $-\sqrt{p/q}$. It easy to see that $F_p$ has
zero mean and unit variance.

\begin{theorem}\label{ThEssRozOsipovLyapFracComparison}
{\rm (i) }
For all $n\in\N,$ $\eps>0,$ and $F_1,\ldots,F_n\in\mathcal F$ such that $B_n>0$ we have
$$
\es^3(\eps,1)\le \lyap,
$$
where the equality takes place for every $n\in\N$ and $\eps>0$ such that $n\eps^{2}\ge1$. As for the extremal, one can take a common d.f. $F_1=\ldots=F_n=F_p$ with arbitrary $p\in\big[\frac12,1\big)$ satisfying $p/(1-p)\le n\eps^2$.

{\rm (ii) } For all $n\in\N,$ $\eps>0,$ and $p\in\big(\frac12,\frac{\sqrt5-1}{2}\big)$ such that $n\eps^2 > (1-p)/p$ we have
$$
\roz^3(\eps,1,F_p)>\lyap(F_p),
$$
in particular,  $\roz^3(\eps,1,F_p)>\es^3(\infty,1,F_p)$,  $\roz^3(\eps,1,F_p)>\esfr^3(\infty)$.

{\rm (iii) }
There exists a d.f. $F\in\mathcal F$ such that for all $n\in\N$ and $\eps>0$ satisfying the condition $n\eps^2\ge\frac{49}{25}=1.96$ we have
$$
\roz^3(\eps,1,F)<\es^3(\eps,1,F),
$$
in particular,  $\sup_{\eps>0}\roz^3(\eps,1,F)<\lyap(F),$  where one can consider $F$ to be a discrete tree-point d.f.

{\rm (iv) } There  exists a common distribution of random summands with the d.f. from $\mathcal F$ such that Esseen--Rozovskii-type inequalities~\eqref{MainEsseenIneq}, \eqref{MainRozovskyIneq} with $\essC(1,1)=\rozC(1,1)=2.73$ are sharper than Osipov's inequality~\eqref{OsipovIneqEps=1} with $C=1.87,$ namely,
$$
2.73\cdot\max\left\{\es^3(1,1),\roz^3(1,1)\right\}\ <\ 1.87 (\osiM(1)+\lind(1)) \quad\text{for every }\ n\ge9.
$$
\end{theorem}

\begin{proof}
The inequality in (i) follows from~\eqref{es^3(eps,1)<=lyap} with $\delta=1.$ To prove that equality also occurs, consider the sequence of i.i.d. r.v.'s $X_1,\ldots,X_n$ with the d.f. $F_p\in\mathcal F$ for $p\in\left[\frac12,1\right)$:
$$
X_1 = \begin{cases}
\sqrt{\frac{q}{p}}, &p,
\\
-\sqrt{\frac{p}{q}}, &q\coloneqq1-p.
\end{cases}
$$
We have  $B_n = \sqrt{n}$,
$$
\sqrt{n}\lyap=\int_{\R}|x|^3dF_p(x)=\frac{p^2+q^2}{\sqrt{pq}} 
\quad\text{for all }n\in\N,
$$
$$
\mu(z)\coloneqq\int_{|x|<z}x^3dF_p(x)=
\begin{cases}
0, & 0<z\le \sqrt{\frac{q}{p}},
\\
\frac{q^2}{\sqrt{pq}}, &\sqrt{\frac{q}{p}} < z \le \sqrt{\frac{p}{q}},
\\
\frac{p-q}{\sqrt{pq}}, &z > \sqrt{\frac{p}{q}},
\end{cases}
\quad
\qt(z)\coloneqq \int_{|x|\ge z}x^2dF_p(x)=
\begin{cases}
1, & 0<z\le \sqrt{\frac{q}{p}},
\\
p, &\sqrt{\frac{q}{p}} < z \le \sqrt{\frac{p}{q}},
\\
0, &z > \sqrt{\frac{p}{q}}.
\end{cases}
$$
With the account of $\mu(z)\ge0$ for all $z>0$ and of  the left-continuity of the functions $\mu(z)$ and $\qt(z)$ for $z>0$, for all $n\in\N$ and $p\in\left[\frac12,1\right)$ we obtain  for every $\eps>0$
\begin{multline*}
\sqrt{n}\es^3(\eps,1,F_p) = \sup_{0<z\le\eps\sqrt{n}}\left\{\mu(z)+z\qt(z)\right\}=
\\
=\begin{cases}
\eps\sqrt{n}, &n\eps^2\le\frac{q}{p},
\\[2mm]
\max\left\{ \sqrt{\frac qp}, \frac{q^2}{\sqrt{pq}}+p\eps\sqrt{n} \right\}=  \frac{q^2}{\sqrt{pq}}+p\eps\sqrt{n}, &\frac{q}{p}<n\eps^2\le\frac{p}{q},
\\[2mm]
\max\left\{ \sqrt{\frac qp}, \frac{q^2}{\sqrt{pq}}+p\sqrt{\frac pq}, \frac{p-q}{\sqrt{pq}}\right\}=
\frac{q^2+p^2}{\sqrt{pq}},
&n\eps^2>\frac{p}{q},
\end{cases}
\end{multline*}
\begin{multline*}
\sqrt{n}\roz^3(\eps,1,F_p) =  \mu\left(\eps\sqrt{n}\right)+\sup_{0<z\le\eps\sqrt{n}}z\qt(z) =
\\
=\begin{cases}
\eps\sqrt{n}, &n\eps^2\le\frac{q}{p},
\\[2mm]
\frac{q^2}{\sqrt{pq}}+\max\left\{ \sqrt{\frac qp}, p\eps\sqrt{n} \right\}
=\frac{q^2+\max\left\{ q,\,\eps p\sqrt{npq} \right\}}{\sqrt{pq}}, & \frac{q}{p}<n\eps^2\le\frac{p}{q},
\\
\frac{p-q}{\sqrt{pq}}+\max\left\{ \sqrt{\frac qp}, p\sqrt{\frac pq}\right\} =
\frac{p-q+\max\left\{ q,p^2 \right\}}{\sqrt{pq}}=
\begin{cases}
\frac{p}{\sqrt{pq}},&\frac12\le p \le \frac{\sqrt{5}-1}{2},\\
\frac{p-q+p^2}{\sqrt{pq}},&\frac{\sqrt{5}-1}{2} < p < 1,
\end{cases}
&n\eps^2 > {\frac{p}{q}}.
\end{cases}
\end{multline*}
Now it is easy to see that for all $p\in\left[\tfrac12,1\right)$ and $n\eps^2\ge p/q$,
$$
\sqrt{n}\es^3(\eps,1,F_{p})= \frac{q^2}{\sqrt{pq}}+p\sqrt{\frac pq} = \frac{q^2+p^2}{\sqrt{pq}} =\sqrt{n}\lyap(F_p),
$$
where, for $n\eps^2=1$, the unique  value $p=\frac12$ is admissible. Thus, (i) is proved.

Now let $\frac12<p<\frac{\sqrt5-1}{2}$. If  $n\eps^2>p/q$, then
$$
\sqrt{n}\roz^3(\eps,1,F_p)=\frac{p}{\sqrt{pq}}\ \ >\ \ \frac{p+q(q-p)}{\sqrt{pq}}= \frac{p^2+q^2}{\sqrt{pq}} =\sqrt{n}\lyap(F_p),
$$
while for $q/p<n\eps^2\le p/q$ we have $q>p^2\ge p\eps\sqrt{npq}$ and, hence,
$$
\sqrt{n}\roz^3(\eps,1,F_p) =\frac{q^2+q}{\sqrt{pq}}\ \  >\ \ \frac{q^2+p^2}{\sqrt{pq}}=\sqrt{n}\lyap(F_p), 
$$
so that $\roz^3(\eps,1,F_p) >\lyap(F_p)\ge\sup_{\eps'>0}\es^3(\eps',1,F_p)$ for all $n\in\N$, $\eps>0,$ and $p\in\big(\frac12,\frac{\sqrt5-1}{2}\big)$ such that $n\eps^2>(1-p)/p$, which proves (ii).

To prove (iii), consider the common three-point distribution of the random summands $X_1,\ldots,X_n$ concentrated in the points
$$
x_1 = \tfrac45,\quad x_2 = -1, \quad x_3 = \tfrac75,
$$
with masses $p,q,r\ge0$ such that
\[
\begin{cases}
p + q + r = 1,\\
\E X_1= px_1+qx_2+rx_3=0,\\
\E X_1^2= px_1^2+qx_2^2+rx_3^2=1,
\end{cases}
\]
that is,
$$
p = \frac{10}{27}= 0.3703\ldots, \quad q = \frac{53}{108}= 0.4907\ldots,\quad r = \frac{5}{36}=0.1388\ldots\,.
$$
Now for $n\eps^2\ge(|x_1|\vee|x_2|\vee|x_3|)^2= x_3^2=\frac{49}{25}=1.96$ with $B_n=\sqrt{n}$ we have
\[
\sqrt{n}\es^3(\eps,1)= \max\left\{
\abs{x_1},
\abs{x_1^3p} + (1-x_1^2p)\abs{x_2},
\abs{x_1^3p+x_2^3q} + x_3^2r \cdot \abs{x_3},
\abs{x_1^3p+x_2^3q + x_3^3r}\right\}=
\]
\[
 =\frac{643}{675}=0.9525\ldots,
\]
\[
\sqrt{n}\roz^3(\eps,1) = \abs{x_1^3p+x_2^3q + x_3^3r} + \max\left\{
\abs{x_1},(1-x_1^2p)\abs{x_2},x_3^2r \cdot \abs{x_3}
\right\} = \frac{22}{25}=0.88<\sqrt{n}\es^3(\eps,1),
\]
which proves (iii) for  $n\eps^2\ge\frac{49}{25}.$

To prove (iv), let us consider the common symmetric four-point distribution of the random summands $X_1,\ldots,X_n$ of the form
\[
X_1 = \begin{cases}
\pm x_1, &p/2\ge0,\\
\pm x_2, &q/2\ge0,
\end{cases}
\quad p+q = 1,
\]
and such that $\E X_1=0,$ $\E X_1^2=1.$ Setting $x_1 = 0.9$ and $x_2=3,$ we find $p = \frac{800}{819}=0.9768\ldots$ and $q=\frac{19}{819}=0.0231\ldots,$ so that for $n\eps^2\ge x_1^2\vee x_2^2= x_2^2=9$ with $B_n=\sqrt{n}$ we have
$$
\sqrt{n}\es^3(\eps,1) = \sqrt{n}\roz^3(\eps,1) = \max \left\{x_1,x_2\left(1-px_1^2\right)\right\}=0.9,
$$
while for $n\ge x_1^2\vee x_2^2=9$ we obtain
$$
\sqrt{n}(\osiM(1)+\lind(1)) = x_1^3p+x_2^3q = \frac{87}{65}=1.3384\ldots\ .
$$
In particular, for $\eps=1$ and $n\ge9$
$$
2.73 \cdot \sqrt{n}\max\{\es^3(1,1),\roz^3(1,1)\} = 2.457\ \ < \ \ 2.5029\ldots = 1.87 \cdot \sqrt{n}(\osiM(1)+\lind(1)),
$$
which proves (iv). We can also propose another ``non-symmetric'' example for even $n$. Let $X_1,\ldots,X_n$ be independent r.v.'s with the distributions $X_k\stackrel{d}{=}X$ if $k$ is odd and $X_k\stackrel{d}{=}-X$ if $k$ is even, where $X$ is a three-point r.v. with
$$
X =
\begin{cases}
    x_1=1/2, &p=\frac{4}{9},\\
    x_2=-1, &q=\frac{4}{9},\\
    x_3=2, &r=\frac{1}{9},
\end{cases}
\qquad p+q+r = 1,\quad\text{so that }\ \E X = 0,\quad\E X^2 = 1.
$$
Taking into account that $\sum_{k=1}^{n}\E X_k^3\I\left(\abs{X_k}<z\right)=0$ for even $n$ and all $z>0$, for even $n\ge\max\{x_1^2,x_2^2,x_3^2\}=4$ we have
$$
\sqrt{n}\es^3(1,1) = \sqrt{n}\roz^3(1,1) = \sup_{0<z\le\sqrt{n}}\left\{z\E X^2\I\left(\abs{X}\ge z\right)\right\}=\max\left\{x_1,(1-px_1^2)\abs{x_2},rx_3^3\right\} = \frac{8}{9},
$$
$$
\sqrt{n}(\osiM(1)+\lind(1)) = px_1^3+q\abs{x_2}^3+rx_3^3 = \frac{25}{18}=1.3888\ldots,
$$
$$
2.73 \cdot \sqrt{n}\max\left\{\es^3(1,1),\roz^3(1,1)\right\} = 2.4266\ldots\ \ < \ \ 2.5972\ldots = 1.87 \cdot \sqrt{n}(\osiM(1)+\lind(1)).
$$
\end{proof}

\begin{theorem}\label{ThAbsChBE(gamma=0,eps)=inf}
For the asymptotically best constant we have
$$
\sup_{F_1=\ldots=F_n\in\mathcal F\colon B_n>0}\ \limsup_{n\to\infty} \frac{\Delta_n(F_1,\ldots,F_n)}{\sup\limits_{0<z\le\eps}zL_n(z)}=\infty \quad\text{for every}\quad \eps>0.
$$
\end{theorem}

\begin{proof}
Using his asymptotic expansion, Esseen deduced~\cite{Esseen1945}
that in the i.i.d. case with $\E X_1=0,$ $\E X_1^2=\sigma^2>0,$
$\E|X_1|^3<\infty$ the limit below exists and
$$
\lim_{n\to\infty}\sqrt{n}\Delta_n=\frac{|\alpha_3|+3h\sigma^2}{6\sqrt{2\pi}\sigma^3},
$$
where $\alpha_3\coloneqq\E X_1^3$, $h>0$ is the span in case of a lattice distribution of $X_1$, and $h\coloneqq0$ otherwise. Now let us consider an absolute continuous distribution of $X_1$ whose d.f. $F_\theta$ is defined by  the density
$$
p_\theta(x)=\begin{cases}
ax^{-4-\theta},&x>1
\\
0, &|x|\le1,
\\
b|x|^{-5},&x<-1,
\end{cases}
\quad\text{with }\ a=a(\theta)\coloneqq \frac{4(2+\theta)(3+\theta)}{17+7\theta},\ \ b=b(\theta)\coloneqq \frac{12(3+\theta)}{17+7\theta}, \ \ \theta\in(0,1).
$$
Then  $h=0,$ $\E X_1=0,$
$$
\sigma^2=\sigma^2(\theta) =\frac{a}{1+\theta} +\frac{b}2 =\frac{2(3+\theta)(7+5\theta)}{(1+\theta)(17+7\theta)},\quad
\alpha_3=\alpha_3(\theta) =\frac{a}{\theta}-b =\frac{4(3+\theta)(2-8\theta-3\theta^2)}{\theta(17+7\theta)},
$$
$$
\E X_1^2\I(|X_1|\ge z)=
\begin{cases}
\sigma^2,&z\in(0,1],
\\
\frac{a}{1+\theta}z^{-1-\theta}+\frac{b}{2}z^{-2},&z>1,
\end{cases}
$$
$$
\lim_{n\to\infty}\sqrt{n}\sup_{0<z\le\eps}zL_n(z) =\sigma^{-3}\lim_{n\to\infty}\sup_{0<z\le\eps\sigma\sqrt{n}}z\E X_1^2\I(|X_1|\ge z) =\sigma^{-3}\sup_{z>0}z\E X_1^2\I(|X_1|\ge z)
=\sigma^{-1},
$$
so that
$$
\sup_{F_1=\ldots=F_n\in\mathcal F\colon B_n>0}\ \limsup_{n\to\infty} \frac{\Delta_n(F_1,\ldots,F_n)}{\sup\limits_{0<z\le\eps}zL_n(z)} \ge \sup_{\theta\in(0,1]}\frac{\lim\limits_{n\to\infty}\sqrt{n}\Delta_n(F_\theta,\ldots,F_\theta)} {\lim\limits_{n\to\infty}\sqrt{n}\sup\limits_{0<z\le\eps}zL_n(z)} = \frac1{6\sqrt{2\pi}}\lim_{\theta\to0}\frac{|\alpha_3(\theta)|}{\sigma^2(\theta)}=\infty.
$$
\end{proof}

\clearpage

\section{Appendix: Figures and Tables}

\begin{table}[h!]
\begin{center}
\begin{tabular}{||c|c||c||c|c|c||c|c|c||}
\hline
\multirow{2}{*}{$\eps$}
&\multirow{2}{*}{$\gamma$}
&\multirow{2}{*}{$C_{0}(\eps,\gamma,0+)$}
&$L=0.001$
&\multicolumn{2}{c||}{Optimal}
&$L=0.03$
&\multicolumn{2}{c||}{Optimal}
\\ \cline{5-6}\cline{8-9}
&&&$C_{0}(\eps,\gamma,L)$&$\tau_0$&$\tau_1$&$C_{0}(\eps,\gamma,L)$&$\tau_0$&$\tau_1$
\\ \hline
$0.6$&$0.3$&$1.92245$&$1.98369$&$0.33196$&$0.76348$&$2.63565$&$0.61269$&$0.58888$\\ \hline 
$1.21$&$0.2$&$1.95457$&$2.02040$&$0.30850$&$0.76192$&$2.70866$&$0.58750$&$0.58212$\\ \hline 
$2.06$&$0.2$&$1.94999$&$2.01563$&$0.30911$&$0.76197$&$2.70121$&$0.58906$&$0.58256$\\ \hline 
$\infty$&$0.2$&$1.94879$&$2.01437$&$0.30943$&$0.76199$&$2.69812$&$0.58998$&$0.58279$\\ \hline 

$1.48$&$0.4$&$1.80997$&$1.86401$&$0.37487$&$0.76589$&$2.45637$&$0.65732$&$0.59975$\\ \hline 
$\infty$&$0.4$&$1.80005$&$1.85350$&$0.37849$&$0.76607$&$2.44001$&$0.66177$&$0.60081$\\ \hline 
$1.89$&$\gopt$&$1.77136$&$1.82167$&$0.39960$&$0.76705$&$2.38689$&$0.68022$&$0.60488$\\ \hline 
$2.03$&$\gopt$&$1.76995$&$1.82017$&$0.40115$&$0.76711$&$2.38467$&$0.68083$&$0.60502$\\ \hline 
$\infty$&$\gopt$&$1.76370$&$1.81351$&$0.40416$&$0.76725$&$2.37413$&$0.68386$&$0.60563$\\ \hline 
$1$&$\gopt$&$1.80596$&$1.85831$&$0.38771$&$0.76651$&$2.43955$&$0.66679$&$0.60193$\\ \hline 
$1$&$0.67$&$1.79961$&$1.85099$&$0.39328$&$0.76673$&$2.42541$&$0.67230$&$0.60312$\\ \hline 
$1$&$\infty$&$1.79149$&$1.83892$&$0.42035$&$0.76791$&$2.38889$&$0.69303$&$0.60759$\\ \hline 
$2.24$&$1$&$1.73996$&$1.78661$&$0.43002$&$0.76828$&$2.32719$&$0.70218$&$0.60952$\\ \hline 
$\infty$&$1$&$1.73186$&$1.77796$&$0.44157$&$0.76870$&$2.31385$&$0.70646$&$0.61040$\\ \hline 
$3.07$&$\infty$&$1.71998$&$1.76313$&$0.45717$&$0.76925$&$2.28233$&$0.72259$&$0.61358$\\ \hline 
$3.2$&$5$&$1.71997$&$1.76354$&$0.45694$&$0.76926$&$2.28502$&$0.72045$&$0.61316$\\ \hline 
$3.28$&$4$&$1.71999$&$1.76368$&$0.45267$&$0.76914$&$2.28573$&$0.71991$&$0.61306$\\ \hline 
$4$&$2.4$&$1.71998$&$1.76401$&$0.45217$&$0.76907$&$2.28788$&$0.71820$&$0.61272$\\ \hline 
$5$&$2.06$&$1.71997$&$1.76413$&$0.45034$&$0.76902$&$2.28870$&$0.71753$&$0.61259$\\ \hline 
$5.37$&$2$&$1.72000$&$1.76420$&$0.45155$&$0.76907$&$2.28892$&$0.71737$&$0.61256$\\ \hline 
$\infty$&$1.83$&$1.71995$&$1.76423$&$0.45135$&$0.76907$&$2.28913$&$0.71712$&$0.61251$\\ \hline 
$\infty$&$\infty$&$1.71451$&$1.75725$&$0.46485$&$0.76952$&$2.27337$&$0.72554$&$0.61416$\\ \hline 
\end{tabular}
\end{center}
\caption{Demonstration to the evaluation of  the upper bound $C_0(\eps,\gamma,L_0)$
defined in~\eqref{C0(L)=infC0(L,tau0,tau1,eps,gamma)}  for the constant $\essC(\eps,\gamma)$ in the Esseen-type inequality~\eqref{MainEsseenIneq} for small values of $L\coloneqq\es(\eps,\gamma)\le L_0$ and some $\eps,\gamma$: Values of  $C_0(\eps,\gamma,L_0)$, rounded up,  for $L_0=0.001$ (the fourth column) and $L_0=0.03$ (the seventh column) accompanied by the corresponding optimal values of the parameters $\tau_0$, $\tau_1$ in~\eqref{C0(L)=infC0(L,tau0,tau1,eps,gamma)}.
The third column provides values of the function $\widehat \aexess(\eps,\gamma)\coloneqq C_0(\eps,\gamma,0+)$ defined in~\eqref{aexess(eps,gamma)<=} which bounds from above the asymptotically exact constant  $\aexess(\eps,\gamma)$ defined in~\eqref{aexCessDef}. Recall that $\gopt=0.5599\ldots\ .$}
\label{Tab:EsseenC0(L)}
\end{table}

\begin{table}[h!]
\begin{center}

\begin{tabular}{||c|c||c||c|c|c||c|c|c||}
\hline
\multirow{2}{*}{$\eps$}
&\multirow{2}{*}{$\gamma$}
&\multirow{2}{*}{$C_{0}(\eps,\gamma,0)$}
&$L=0.001$
&\multicolumn{2}{c||}{Optimal}
&$L=0.03$
&\multicolumn{2}{c||}{Optimal}
\\ \cline{5-6}\cline{8-9}
&&&$C_{0}(\eps,\gamma,L)$&$\tau_0$&$\tau_1$&$C_{0}(\eps,\gamma,L)$&$\tau_0$&$\tau_1$
\\ \hline
$1.21$&$0.2$&$1.93474$&$1.99463$&$0.33997$&$0.76401$&$2.63576$&$0.62018$&$0.59077$\\ \hline 
$1.89$&$0.2$&$1.92998$&$1.98967$&$0.33862$&$0.76391$&$2.62929$&$0.62130$&$0.59105$\\ \hline 
$2.77$&$0.2$&$1.92890$&$1.98855$&$0.33849$&$0.76389$&$2.62857$&$0.62119$&$0.59103$\\ \hline 
$5.39$&$0.2$&$1.95832$&$2.01917$&$0.33339$&$0.76358$&$2.67331$&$0.61213$&$0.58874$\\ \hline 
$1.41$&$0.4$&$1.77974$&$1.82650$&$0.42470$&$0.76809$&$2.37242$&$0.69715$&$0.60847$\\ \hline 
$1.76$&$0.4$&$1.77249$&$1.81886$&$0.43083$&$0.76830$&$2.36208$&$0.69946$&$0.60894$\\ \hline 
$1.99$&$0.4$&$1.77128$&$1.81759$&$0.43528$&$0.76848$&$2.36048$&$0.69982$&$0.60902$\\ \hline 
$2.63$&$0.4$&$1.77841$&$1.82511$&$0.42595$&$0.76814$&$2.37148$&$0.69707$&$0.60846$\\ \hline 
$0.5$&$\gopt$&$1.94743$&$1.99139$&$0.43510$&$0.76847$&$2.54844$&$0.68154$&$0.60513$\\ \hline 
$1$&$\gopt$&$1.79154$&$1.83112$&$0.48527$&$0.77017$&$2.34604$&$0.72358$&$0.61378$\\ \hline 
$1.52$&$\gopt$&$1.74995$&$1.78796$&$0.49835$&$0.77051$&$2.29031$&$0.73728$&$0.61641$\\ \hline 
$1.89$&$\gopt$&$1.74383$&$1.78159$&$0.50808$&$0.77078$&$2.28222$&$0.73934$&$0.61679$\\ \hline 
$1.99$&$\gopt$&$1.74412$&$1.78189$&$0.50748$&$0.77083$&$2.28271$&$0.73918$&$0.61676$\\ \hline 
$2.12$&$\gopt$&$1.74542$&$1.78324$&$0.50216$&$0.77063$&$2.28462$&$0.73863$&$0.61666$\\ \hline 
$3$&$\gopt$&$1.77092$&$1.80977$&$0.48990$&$0.77029$&$2.32025$&$0.72921$&$0.61487$\\ \hline 
$5$&$\gopt$&$1.86500$&$1.90688$&$0.45666$&$0.76924$&$2.44632$&$0.70070$&$0.60920$\\ \hline 
\end{tabular}
\end{center}
\caption{Demonstration to the evaluation of  the upper bound $C_0(\eps,\gamma,L_0)$ 
defined in~\eqref{C0(L)=infC0(L,tau0,tau1,eps,gamma)}  for the constant $\rozC(\eps,\gamma)$ in the Rozovskii-type inequality~\eqref{MainRozovskyIneq} for small values of $L\coloneqq\roz(\eps,\gamma)\le L_0$ and some $\eps,\gamma$: Values of  $C_0(\eps,\gamma,L_0)$, rounded up,  for $L_0=0.001$ (the fourth column) and $L_0=0.03$ (the seventh column) accompanied by the corresponding optimal values of the parameters $\tau_0$, $\tau_1$ in~\eqref{C0(L)=infC0(L,tau0,tau1,eps,gamma)}.
The third column provides values of the function $\widehat \aexroz(\eps,\gamma)\coloneqq C_0(\eps,\gamma,0+)$ defined in~\eqref{aexroz(eps,gamma)<=} which bounds from above the asymptotically exact constant  $\aexroz(\eps,\gamma)$ defined in~\eqref{aexCrozDef}. Recall that $\gopt=0.5599\ldots\ .$}
\label{Tab:RozovskyC0(L)}
\end{table}

\begin{table}[h!]
\begin{center}
\begin{tabular}{||c|c|c|c||c|c||c|c|c|c||}
\hline
\multirow{2}{*}{$\eps$}
&\multirow{2}{*}{$\gamma$}
&\multirow{2}{*}{$C_1\le$}
&\multirow{2}{*}{$L_*$}
&\multicolumn{2}{c||}{Optimal $T_0,T_1$}
&\multicolumn{4}{c||}{Contributions of $I_k/L_*^3$\vphantom{$\frac{A^A}{}$} }
\\ \cline{5-10}
&&&&${T_0}L_*$& $T_1L_*^3\vphantom{\displaystyle\frac12}$& $k=1$&$k=2$&$k=3$&$k=4$
\\ \hline
$1.21$&$0.2$&$2.89038$&$0.47850$&$0.87755$&$0.69402$&$0.44928$&$0.81565$&$1.51412$&$0.11135$\\ \hline 
$1.24$&$0.2$&$2.88998$&$0.47856$&$0.87744$&$0.69400$&$0.44866$&$0.81589$&$1.51400$&$0.11144$\\ \hline 
$\infty$&$0.2$&$2.88457$&$0.47857$&$0.87749$&$0.69398$&$0.44331$&$0.81581$&$1.51405$&$0.11142$\\ \hline 

$1.76$&$0.4$&$2.73593$&$0.48170$&$0.92023$&$0.69306$&$0.35252$&$0.76761$&$1.52685$&$0.08897$\\ \hline 
$5.94$&$0.4$&$2.73000$&$0.48173$&$0.92020$&$0.69311$&$0.34652$&$0.76781$&$1.52668$&$0.08900$\\ \hline 
$\infty$&$0.4$&$2.72984$&$0.48173$&$0.92020$&$0.69311$&$0.34635$&$0.76781$&$1.52669$&$0.08900$\\ \hline 
$1$&$\gopt$&$2.73662$&$0.48267$&$0.92963$&$0.69278$&$0.36492$&$0.75795$&$1.52899$&$0.08477$\\ \hline 
$1.87$&$\gopt$&$2.69989$&$0.48245$&$0.92958$&$0.69285$&$0.32825$&$0.75774$&$1.52921$&$0.08471$\\ \hline 
$\infty$&$\gopt$&$2.69190$&$0.48244$&$0.92960$&$0.69285$&$0.32028$&$0.75770$&$1.52924$&$0.08470$\\ \hline 
$1$&$0.72$&$2.72979$&$0.48304$&$0.93393$&$0.69261$&$0.36342$&$0.75335$&$1.53015$&$0.08289$\\ \hline 
$1$&$\infty$&$2.72857$&$0.48346$&$0.93595$&$0.69239$&$0.36458$&$0.75136$&$1.53053$&$0.08211$\\ \hline 
$4.35$&$1$&$2.66000$&$0.48305$&$0.93735$&$0.69266$&$0.29795$&$0.74958$&$1.53117$&$0.08131$\\ \hline 
$\infty$&$1$&$2.65879$&$0.48305$&$0.93734$&$0.69266$&$0.29673$&$0.74959$&$1.53117$&$0.08132$\\ \hline 
$\infty$&$0.97$&$2.65985$&$0.48303$&$0.93710$&$0.69266$&$0.29749$&$0.74984$&$1.53111$&$0.08142$\\ \hline 
$2.56$&$\infty$&$2.64999$&$0.48339$&$0.94137$&$0.69256$&$0.29282$&$0.74546$&$1.53211$&$0.07962$\\ \hline 
$2.62$&$5$&$2.64992$&$0.48338$&$0.94121$&$0.69256$&$0.29257$&$0.74561$&$1.53207$&$0.07968$\\ \hline 
$2.65$&$4$&$2.64996$&$0.48342$&$0.94113$&$0.69253$&$0.29249$&$0.74571$&$1.53205$&$0.07974$\\ \hline 
$2.74$&$3$&$2.64995$&$0.48335$&$0.94092$&$0.69257$&$0.29224$&$0.74592$&$1.53200$&$0.07980$\\ \hline 
$3.13$&$2$&$2.64997$&$0.48329$&$0.94031$&$0.69259$&$0.29153$&$0.74652$&$1.53188$&$0.08005$\\ \hline 
$4$&$1.62$&$2.64999$&$0.48325$&$0.93979$&$0.69259$&$0.29091$&$0.74705$&$1.53176$&$0.08027$\\ \hline 
$5.37$&$1.5$&$2.64993$&$0.48323$&$0.93954$&$0.69260$&$0.29055$&$0.74732$&$1.53170$&$0.08038$\\ \hline 
$\infty$&$1.43$&$2.64991$&$0.48321$&$0.93935$&$0.69260$&$0.29030$&$0.74750$&$1.53166$&$0.08046$\\ \hline 
$\infty$&$\infty$&$2.64082$&$0.48338$&$0.94138$&$0.69255$&$0.28367$&$0.74541$&$1.53214$&$0.07961$\\ \hline 
\end{tabular}
\end{center}
\caption{Demonstration to the evaluation of the upper bound $\max\limits_{L_0\le L\le L_1}C_1(\eps,\gamma,L)$ (see~\eqref{C1(L)def}) for the constant $\essC(\eps,\gamma)$ in the Esseen-type inequality~\eqref{MainEsseenIneq} for moderate values of $L\coloneqq\es(\eps,\gamma)\in[L_0,L_1]=[0.03,0.65]$ and some $\eps,\gamma$: Extreme values of $C_1(\eps,\gamma,L)\le C_1(\eps,\gamma,L_*)$ on the interval $L\in[L_0,L_1]$ (column~3); maximizer $L_*$, rounded down (column~4); optimal values of the parameters  $T_0$ and $T_1$ in~\eqref{C1(L)def} multiplied  by $L$ and $L^3$ in the extremal point $L=L_*$, rounded down (columns 5, 6);  values of the normalized integrals $I_k/L_*^3,$ $k=1,2,3,4,$ rounded down (columns 7--10), so that $C_1(\eps,\gamma,L_*)=(I_1+I_2+I_3+I_4)/L_*^3$. Recall that $\gopt=0.5599\ldots\ .$}
\label{Tab:EsseenC1(L)}
\end{table}

\begin{table}[h!]
\begin{center}
\begin{tabular}{||c|c|c|c||c|c||c|c|c|c||}
\hline
\multirow{2}{*}{$\eps$}
&\multirow{2}{*}{$\gamma$}
&\multirow{2}{*}{$C_1\le$}
&\multirow{2}{*}{$L_*$}
&\multicolumn{2}{c||}{Optimal $T_0,T_1$}
&\multicolumn{4}{c||}{Contributions of $I_k/L_*^3$\vphantom{$\frac{A^A}{}$} }
\\ \cline{5-10}
&&&&${T_0}{L_*}$& ${T_1}{L_*^3}\vphantom{\displaystyle\frac12}$ &$k=1$&$k=2$&$k=3$&$k=4$
\\ \hline
$1.21$&$0.2$&$2.86991$&$0.47901$&$0.88360$&$0.69388$&$0.43734$&$0.80872$&$1.51596$&$0.10790$\\ \hline 
$5.39$&$0.2$&$2.86343$&$0.47902$&$0.88371$&$0.69384$&$0.43101$&$0.80851$&$1.51608$&$0.10784$\\ \hline 
$1.76$&$0.4$&$2.69985$&$0.48249$&$0.93002$&$0.69288$&$0.32875$&$0.75737$&$1.52923$&$0.08452$\\ \hline 
$2.63$&$0.4$&$2.69323$&$0.48229$&$0.92948$&$0.69296$&$0.32154$&$0.75779$&$1.52922$&$0.08470$\\ \hline 
$0.5$&$\gopt$&$3.03953$&$0.50011$&$0.92200$&$0.66074$&$0.52576$&$0.85947$&$1.55924$&$0.09507$\\ \hline 
$1$&$\gopt$&$2.72857$&$0.48346$&$0.93595$&$0.69239$&$0.36458$&$0.75136$&$1.53053$&$0.08211$\\ \hline 
$1.99$&$\gopt$&$2.65991$&$0.48300$&$0.93912$&$0.69274$&$0.30011$&$0.74764$&$1.53169$&$0.08048$\\ \hline 
$2.12$&$\gopt$&$2.65925$&$0.48273$&$0.93728$&$0.69293$&$0.29727$&$0.74960$&$1.53117$&$0.08122$\\ \hline 
$3$&$\gopt$&$2.67687$&$0.48125$&$0.92365$&$0.69352$&$0.29824$&$0.76363$&$1.52790$&$0.08710$\\ \hline 
$5$&$\gopt$&$2.75611$&$0.47832$&$0.89074$&$0.69456$&$0.33467$&$0.79931$&$1.51873$&$0.10341$\\ \hline 
\end{tabular}
\end{center}
\caption{Demonstration to the evaluation of the upper bound $\max\limits_{L_0\le L\le L_1}C_1(\eps,\gamma,L)$ (see~\eqref{C1(L)def}) for the constant $\rozC(\eps,\gamma)$ in the Rozovskii-type inequality~\eqref{MainRozovskyIneq} for moderate values of $L\coloneqq\roz(\eps,\gamma)\in[L_0,L_1]=[0.03,0.65]$ and some $\eps,\gamma$: Extreme values of $C_1(\eps,\gamma,L)\le C_1(\eps,\gamma,L_*)$ on the interval $L\in[L_0,L_1]$ (column~3); maximizer $L_*$, rounded down (column~4); optimal values of the parameters  $T_0$ and $T_1$ in~\eqref{C1(L)def} multiplied  by $L$ and $L^3$ in the extremal point $L=L_*$, rounded down (columns 5, 6);  values of the normalized integrals $I_k/L_*^3,$ $k=1,2,3,4,$ rounded down (columns 7--10), so that $C_1(\eps,\gamma,L_*)=(I_1+I_2+I_3+I_4)/L_*^3$. Recall that $\gopt=0.5599\ldots\ .$}
\label{Tab:RozovskyC1(L)}
\end{table}

\end{document}